\documentclass[onefignum,onetabnum]{siamart171218}

\usepackage{lipsum}
\usepackage{amsfonts, amssymb}
\usepackage[mathscr]{eucal}
\usepackage{graphicx}
\usepackage{subfig}
\usepackage{float}
\usepackage{color}
\usepackage{hyperref}
\usepackage[utf8]{inputenc}
\usepackage[english]{babel}
\usepackage{algorithm}
\usepackage{algpseudocode}

\newtheorem{problem}{Problem}[section]
\newtheorem{Problem*}{Problem}

\newtheorem{remark}{Remark}[section]
\numberwithin{equation}{section}
\newcommand{\height}{4.5cm}

\newcommand{\x}{{\bf x}}

\headers{X-ray Tomography with incomplete data}{M. V. Klibanov  and L. H. Nguyen}

\title{PDE-based numerical method for a limited angle X-ray tomography\thanks{ {\bf Funding}: This work was partially supported by the US Army Research Laboratory and US
Army Research Office grant W911NF-15-1-0233 as well as by the Office of
Naval Research grant N00014-15-1-2330.}
}

\author{Michael V. Klibanov\thanks{Department of Mathematics and Statistics, University of North Carolina at
Charlotte, Charlotte, NC 28223, mklibanv@uncc.edu}
\and Loc H. Nguyen\thanks{Department of Mathematics and Statistics, University of North Carolina at
Charlotte, Charlotte, NC 28223, loc.nguyen@uncc.edu}
}

\begin{document}

\maketitle

\begin{abstract}
 A new numerical method for X-ray tomography for a specific case of
incomplete Radon data is proposed. Potential applications are in checking
out bulky luggage in airports. This method is based on the analysis of the
transport PDE governing the X-ray tomography rather than on the conventional
integral formulation. The quasi-reversibility method is applied. Convergence
analysis is performed using a new Carleman estimate. Numerical results are
presented and compared with the inversion of the Radon transform using the
well-known filtered back projection algorithm. In addition, it is shown how
to use our method to study the inversion of the attenuated X-ray transform
for the same case of incomplete data.

\end{abstract}

\begin{keywords}
  tomographic inverse problem, X-ray transform,
incomplete data, Carleman estimate
\end{keywords}

\begin{AMS}
 	35R30, 44A12
\end{AMS}

\section{Introduction}

\label{sec:1}

Computing a function from its Radon transform, which was first introduced by
Radon in 1917 \cite{Radon1917, Radon:IEEE1986}, is considered as the theory
behind the first commercial computed tomography (CT) scanner, invented by
Hounsfield. Due to this contribution, Hounsfield was awarded the Nobel prize
in 1979. In the current paper, we develop a new numerical method for this
inverse problem for the case of a special type of limited angle data for the
Radon transform. This type of limited angle data might find applications in
checking out baggages in airports as well as checking out interior
structures of walls.

The Radon transform of a function $f$ is the integral of $f$ on a set of
segments of straight lines. If that set includes all straight lines in the
plane, we say that the Radon transform transform data are completely given.
The exact reconstruction of the function $f$ from its complete Radon
transform data can be computed by the well-known filtered back projection
algorithm, see \cite{Natterer:cmsiam2001, PanSidkyVannier:ip2009}. The full
observation of the data is important since the filtered back projection
formula involves a non local operator. In some applications, due to some
technical reasons, the complete Radon transform cannot be collected. We
refer the reader to \cite{Borg,Borg:Meas2017, Borg:report2017} for some
circumstances about this incompleteness, e.g., when the $X$-rays are blocked
by metal bars, see \cite[Section 7]{Borg} for a detailed discussion. On the
other hand, in this paper, we design another experimental situation in which
a large amount of the Radon data is lost, see Section \ref{sec problem
statement}. Our goal is to image objects (or equivalently to determine a
function) when an interval of view angles is limited in a special way. This
situation includes \textit{Radon transform for the limited angle problem} 
\cite{Louis:nm1985}. The missing data leads to the instability of the
reconstruction. We cite some important papers \cite{Linh:ip2016, Borg,
Frikel:ip2013, Linh:SIAM2015, Linh:jfaa2017} and references therein that
characterize and (or) introduce the strategies to reduce the resulting 
\textit{artifacts}, which appear in the reconstructed image. To treat the
case of incompleteness, one might non-rigorously fill the missing data by
the number $0$, see figures 2-4 in \cite{Borg}.

In this paper, we propose a new numerical method which analytically and
numerically produces a good approximation of the inversion of the Radon
transform for a special case of limited angle data. Unlike the filtered back
projection algorithm, our approach is not based on the integral form of the
Radon transform and is not intended to derive an explicit inversion formula.
We use a boundary value problem for a linear partial differential system
whose solution directly yields the solution to the inverse problem. Instead
of using the conventional \textquotedblleft integral" approach, we consider
a well known PDE \cite{HasanogluRomanov:s2017} governing the propagation of
X-rays. In fact, this is the stationary transport PDE without absorption and
integral terms. That PDE involves two unknown functions: its solution $u$
and the target function of interest $f$. At each point of the boundary, one
of boundary conditions for $u$ is exactly the integral along a line segment,
which is considered in Radon transform. In fact, boundary conditions for
that PDE are over determined ones. Using one of ideas of the
Bukhgeim-Klibanov method \cite{BukhKlib}, we next differentiate that PDE
with respect to the source location to obtain another PDE, in which the
function $f$ is not involved, see \cite{Klibanov:jiipp2013} for a survey of
the method of \cite{BukhKlib} as well as books \cite{BK,BY,KT}. That new
equation contains the unknown function $u$ as well as its partial derivative
with respect to the source location. However, a theory on how to solve the
resulting over determined boundary value problem is not available yet. In
this paper, we only approximate the solution of this problem by the solution
of an overdetermined boundary value problem for a linear system of coupled
PDEs of the first order. The solution of the latter problem is used to
compute a partial sum of the Fourier series for the function $u$ with
respect to a special orthonormal basis. This \textquotedblleft cut-off"
technique was first introduced in \cite{Klibanov:jiip2017} for a class of
coefficient inverse problems. Then, it was successfully applied to
numerically solve some coefficient inverse problems \cite{EIT,KlibKol}.
Having that approximation for the function $u$ in hands, we compute the
corresponding approximation for the function $f$ directly.

As mentioned in the above paragraph, to obtain that system of PDEs, we
truncate the Fourier series with respect to a special orthonormal basis and
assume that the corresponding approximation of the function $u$ still
satisfies the above mentioned PDE. Let $N$ be the number of terms of that
truncated series. Even though the original series converges of course as $%
N\rightarrow \infty ,$ the question about the convergence of resulting
numerical solutions as $N\rightarrow \infty $ is a very challenging one. The
true \textquotedblleft hidden" reason of this challenge is the ill-posedness
of the originating problem. Thus, we do not provide here the proof of
convergence of those numerical solutions at $N\rightarrow \infty $. We
estimate an optimal number $N$ numerically, see Remark 5.3 in Section 5.4.
In other words, we consider an \emph{approximate mathematical model}, which
is a common place in numerical methods for ill-posed problems. Indeed, it is
well known that proofs of convergence of numerical solutions resulting from
truncations of a variety of Fourier series, as $N\rightarrow \infty ,$ are
quite challenging ones in many other inverse/ill-posed problems. Hence,
these proofs are usually omitted. Still, it is also well known that
approximate mathematical models based on truncated Fourier series work
successfully numerically even for coefficient inverse problems, which are
nonlinear, unlike the linear problem of this paper. As some examples of
those successes, we refer to, e.g. works of Kabanikhin with coauthors \cite%
{Kab1,Kab2,Kab3} for the 2D version of the Gelfand-Levitan-Krein method, as
well as to publications \cite{KlTh,EIT,KlibKol} of the first author with
coauthors.

The above mentioned overdetermined boundary value problem for a system of
PDEs of the first order is solved here by the quasi-reversibility method,
which is well known to be a perfect tool to solve overdetermined boundary
value problems for PDEs. This method was first introduced by Latt\`{e}s and
Lions \cite{LattesLions:e1969} for numerical solutions of ill-posed problems
for PDEs. It has been studied intensively since then, see e.g., \cite%
{Becacheelal:AIMS2015, Bourgeois:ip2006, BourgeoisDarde:ip2010,
ClasonKlibanov:sjsc2007, Dadre:ipi2016, KS, Klibanov:jiipp2013, Loc:Arxiv2018}. A recent
survey on this method can be found in \cite{Ksurvey}.

In the convergence analysis of this paper we consider a semi discrete form
of our system of PDEs, which is more realistic for computations than the
conventional continuous form. More precisely, we assume that partial
derivatives with respect to one of two variables are written in finite
differences, whereas derivatives with respect to the second variable are
written in the conventional continuous form. However, we do not allow the
step size of the grid $h\rightarrow 0,$ unlike many conventional well posed
problems for PDEs. Indeed, the analysis at $h\rightarrow 0$ is a very
challenging one due to the ill-posedness of the problem. As to the fully
discrete form, in which both partial derivatives are written via finite
differences, it is clear from, e.g. \cite{KS}, that, in the case of
ill-posed problems (as opposed to some conventional well posed problems for
PDEs), this case is far more complicated. Thus, it is not considered in this
first publication about our new method.

It is well known that proofs of convergence of regularized solutions of the
quasi-reversibility method are based on Carleman estimates, see, e.g. \cite%
{Klibanov:jiipp2013,Ksurvey}. Hence, first, we prove a new Carleman
estimate. Next, using this estimate, we prove the existence and uniqueness
of the minimizer (i.e., the regularized solution \cite{T}) for our semi
discrete version of the quasi reversibility method. Finally, using the same
Carleman estimate, we establish a convergence rate of regularized solutions
to the exact solution.

An important part of the paper is devoted to the numerical implementation of
our method. We present here some numerical results. In particular, we
compare performance of our method with the performance of the filtered back
projection algorithm in which the missed data are filled by zeros. We point
out that, at this early stage of the development, we are not interested in
treating fine details, such as artifacts, for example. Rather, we arrange a
simple post processing.

The paper is organized as follows. In Section \ref{sec problem statement} we
state the problem. In Section \ref{sec 4} we derive the above mentioned
overdetermined boundary value problem for a system of PDEs of the first
order, which does not contain the target function $f$. In Section \ref{sec 5}
we introduce first the quasi-reversibility method to solve that problem. \
Next, we prove a new Carleman estimate and use this Carleman estimate to
prove the existence and uniqueness of the minimizer and establish the
convergence rate of the minimizers to the exact solution as the level of the
measurement noise tends to zero. In Section \ref{sec imple} we discuss the
numerical implementation of our method. Numerical studies are described in
Section \ref{sec num}. We present concluding remarks in Section \ref{sec
concl}. In addition, we explain in Section \ref{sec concl} how to extend our
approach to solve the inverse attenuated tomographic problem. Below all
functions are real valued ones.

\section{Problem statement}

\label{sec problem statement}

Everywhere below all functions are real valued ones and $\mathbf{x}=\left(
x,y\right) $ denotes points in $\mathbb{R}^{2}.$ Let $b>a>0$ and $d,R>0$ be
some numbers. Consider the rectangle $\Omega \subset \mathbb{R}^{2}$ 
\begin{equation}
\Omega =(-R,R)\times (a,b).  \label{1}
\end{equation}%
Let $\Gamma _{d}\subset \mathbb{R}^{2}$ be the segment of the horizontal
line where our point sources are located, 
\begin{equation}
\Gamma _{d}=\left\{ \mathbf{x}=\left( x,y\right) :x\in (-d,d),y=0\right\} .
\label{2.10}
\end{equation}%
Let $f(\mathbf{x})$ be the unknown function whose support is contained in $%
\Omega $, i.e. 
\begin{equation}
f\left( \mathbf{x}\right) =0\text{ for }\mathbf{x}\in \mathbb{R}%
^{2}\setminus \Omega .  \label{2.1}
\end{equation}%
For the purpose of our theoretical analysis, we assume below that $f\in
C^{2}(\mathbb{R}^{2}).$ In the case of X-ray tomography the function $%
f\left( \mathbf{x}\right) $ represents the X-ray attenuation coefficient at
the point $\mathbf{x}$, see \cite{Natterer:cmsiam2001}. Consider point
sources $\mathbf{x}_{\alpha }=\left( \alpha ,0\right) \in \Gamma _{d}.$ We
define the function $u(\mathbf{x},\mathbf{x}_{\alpha })$ as 
\begin{equation}
u(\mathbf{x},\mathbf{x}_{\alpha })=\int_{L(\mathbf{x},\mathbf{x}%
_{\alpha })}f(\mathbf{\xi })d\sigma ,\quad  \label{2.2}
\end{equation}%
where $L(\mathbf{x},\mathbf{x}_{\alpha })$ is the line segment connecting
points $\mathbf{x}$ and $\mathbf{x}_{\alpha }$. We are interested in the
following problem:

\begin{problem}[Tomographic inverse problem with incomplete data]
Determine the function $f$ from the measurement of $Rf$, where 
\begin{equation}
Rf=u(\mathbf{x},\mathbf{x}_{\alpha })  \label{2.3}
\end{equation}%
for all $\mathbf{x}=(x,y)\in \partial \Omega $ and all $\mathbf{x}_{\alpha
}=(\alpha ,0)\in \Gamma _{d}.$ The function $Rf$ is known as the Radon
transform of the function $f$. \label{pro 2.1}
\end{problem}

\begin{remark}
The Radon transform, along with the inversion formula, was first introduced
by Radon in 1917 in his celebrated paper \cite{Radon1917}. We also refer the
reader to \cite{Radon:IEEE1986} for the translation of \cite{Radon1917} into
English .
\end{remark}

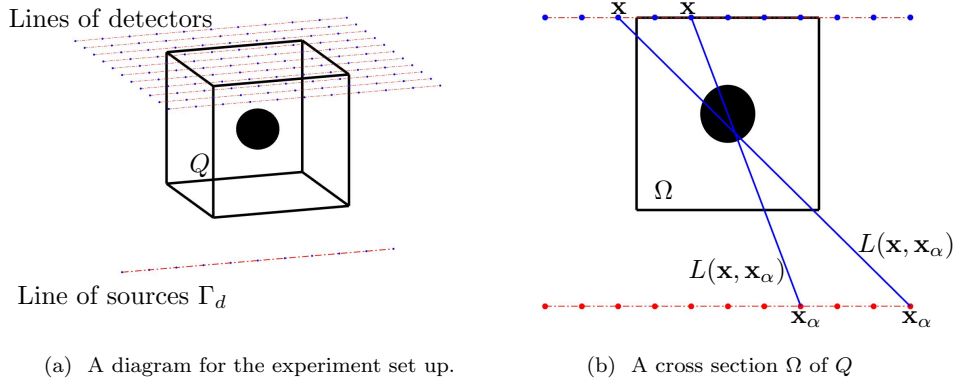
\begin{figure}[h]
\begin{center}
\subfloat[\label{fig 1a} A diagram for the experiment set up.]{\includegraphics[width=0.48\textwidth]{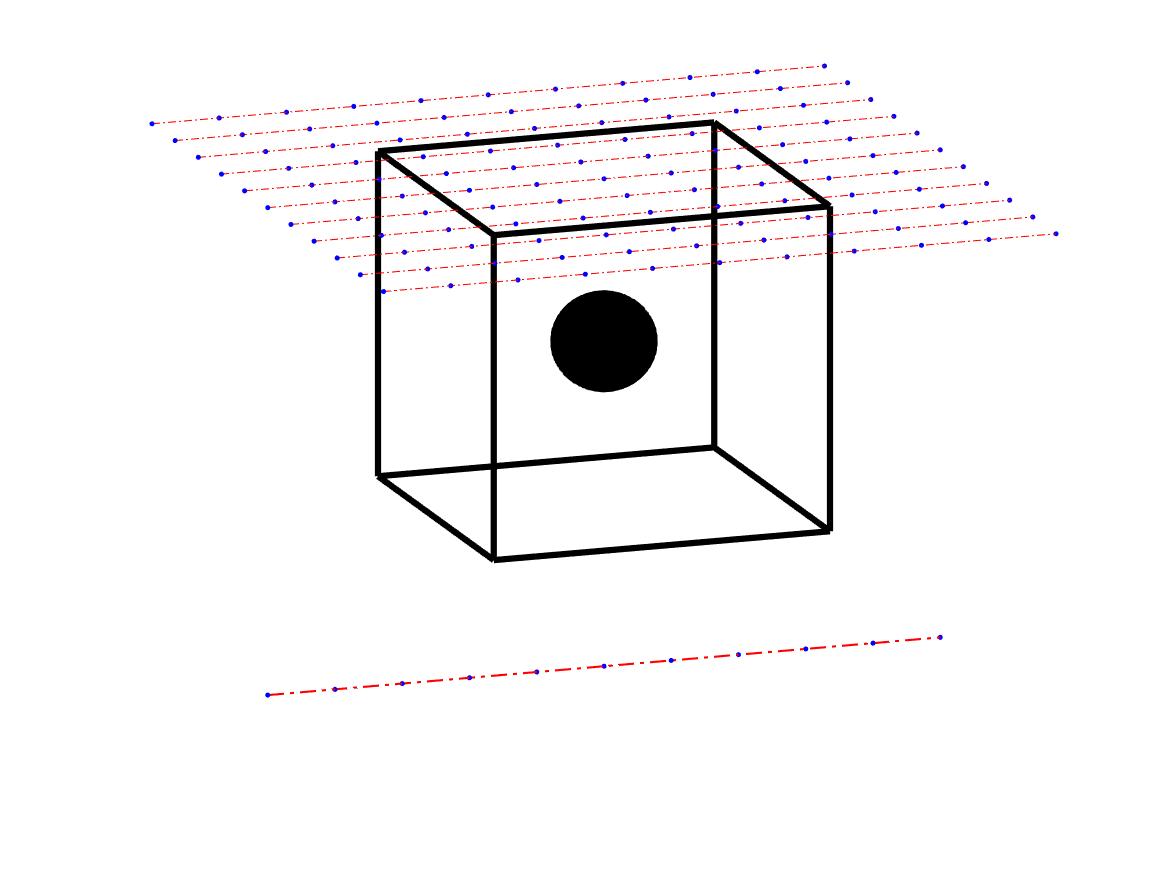} 
\put(-180,120){Lines of detectors}
\put(-180,15){ Line of sources $\Gamma_d$}
\put(-115,65){ $Q$}
} 
\subfloat[\label{fig 1b} A cross section $\Omega$ of $Q$]{\includegraphics[width=0.48\textwidth]{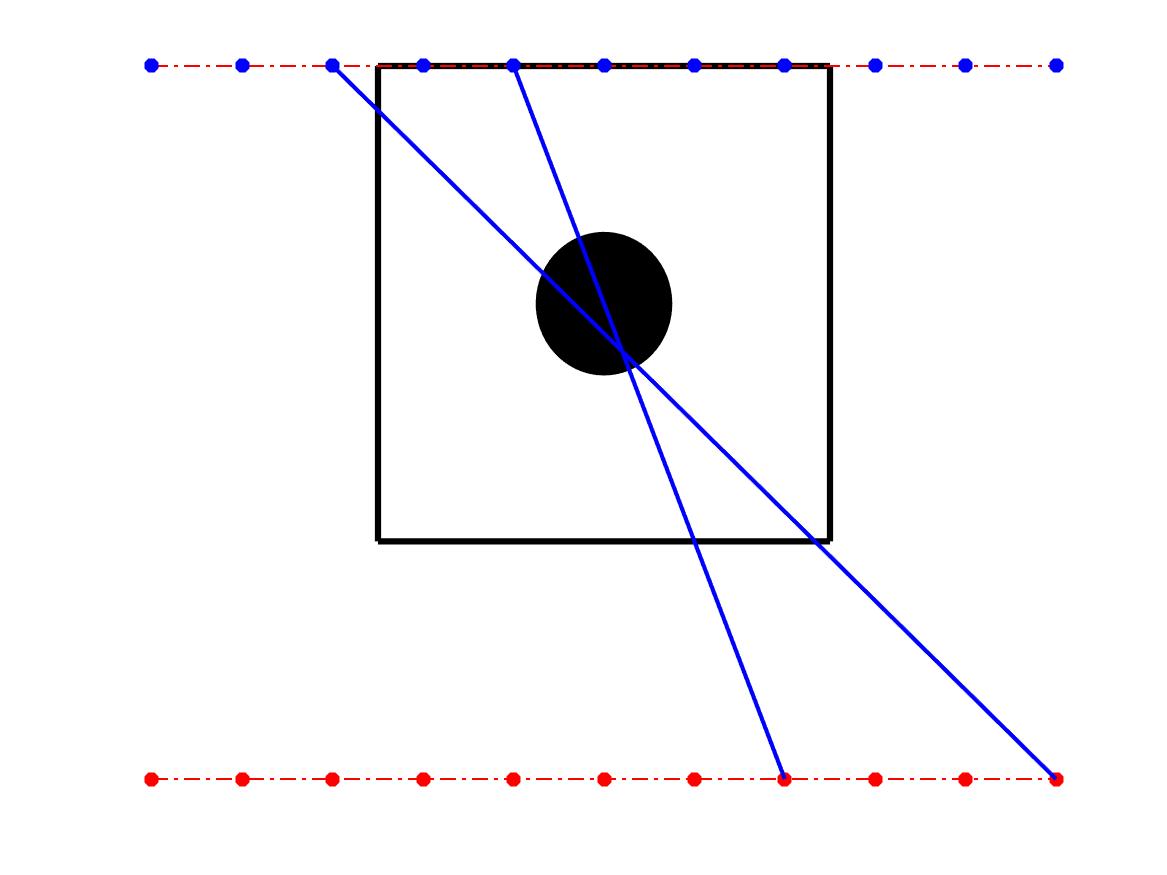}
\put(-62,8){$\x_\alpha$}
\put(-20,8){$\x_\alpha$}
\put(-130,125){$\x$}
\put(-104,125){$\x$}
\put(-40,35){ $L(\x, \x_\alpha)$}
\put(-104,25){ $L(\x, \x_\alpha)$}
\put(-117,55){ $\Omega$}
}
\end{center}
\caption{\textit{An illustration of a 3D tomographic experiment. One can
detect the 3D object in $Q$ by repeatedly solving Problem \protect\ref{pro
2.1} at each cross section $\Omega $ of $Q$ on the plane defined by the line
of source and each observation line. In (b), $\mathbf{x}$ represents the
location of detectors and $\mathbf{x}_{\protect\alpha }$ denotes the
locations of the source. In this tomographic setting, the line $L(\mathbf{x},%
\mathbf{x}_{\protect\alpha })$ in (b) is assumed to be the geodesic line
connecting $\mathbf{x}_{\protect\alpha }$ and $\mathbf{x}$. }}
\label{fig 1}
\end{figure}

The case when the data $Rf(\mathbf{x},\mathbf{x}_{\alpha })$ are available
for all $\mathbf{x}_{\alpha }$ and $\mathbf{x}$ such that the set of lines $%
L(\mathbf{x},\mathbf{x}_{\alpha })$ contains all possible lines intersecting 
$\Omega $, Problem \ref{pro 2.1} is known as the tomographic inverse problem
with \textit{complete} data. This inverse problem with complete data is
exactly solved by the filtered back projection formula \cite%
{Natterer:cmsiam2001}. Unlike this, in the current paper, the point source $%
\mathbf{x}_{\alpha }$ is allowed to \textquotedblleft move" only along the
line segment $\Gamma _{d}$, which is located below $\Omega ,$ rather than on
a curve surrounding $\Omega ,$ as illustrated in Figure \ref{fig 1}. In this
setting, one can easily find many straight lines intersecting $\Omega $ but
not belonging to our set of lines $L(\mathbf{x},\mathbf{x}_{\alpha }).$
Therefore, the data in Problem \ref{pro 2.1} is said to be \textit{incomplete%
}. See Figures \ref{fig Case 1 full}--\ref{fig Case 7 full} versus Figures %
\ref{fig Case 1 incomplete}--\ref{fig Case 7 incomplete} for the
illustrations of the amount of missing data.

Due to a large amount of missing data, the Radon inversion via the
well-known filtered back projection algorithm built in MATLAB does not work
well. In addition, this formula is not rigorously established for this case.
These motivate us to develop a new numerical method to solve Problem \ref%
{pro 2.1}. We use the well-known transport PDE that governs the function $Rf$%
. Next, we establish and solve an inverse source problem for this equation.
This is our PDE approach.

Problem \ref{pro 2.1} arises in X-ray tomography. Assume that we want to
image an object in a 3D domain $Q$, illustrated on Figure \ref{fig 1a}. A
source, located at each point $\mathbf{x}_{\alpha }$ on the line $\Gamma
_{d} $ in (\ref{2.10}) below $Q$, generates tomographic data that can be
measured at an array of detectors on a rectangle on the top of $Q$. One can
arrange such detectors on a set of \textquotedblleft observation lines" that
are parallel to $\Gamma _{d}$. Each observation line, together with $\Gamma
_{d}$, defines a plane. The cross section of $Q$ by that plane is our 2D
domain $\Omega $. Figure \ref{fig 1b} illustrates an example of such cross
section. Hence, we believe that results of this paper have potential
applications in, e.g. checking out a bulky baggage in airports.

We next discuss the issue of the data for Problem \ref{pro 2.1} on the
boundary of $\Omega $. The data on the top of $\Omega $ can be collected
directly. As to the data on two vertical sides of $\Omega ,$ one can easily
see from Figure \ref{fig 1b} that if the measurement line is sufficiently
long, then we can calculate the data on the sides of $\Omega $ using the
data on the measurement side as well as (\ref{2.1}) and (\ref{2.2}). Also,
by (\ref{2.1}) and (\ref{2.2}) the data on the bottom side of $\Omega $ is
identically zero. Solving problem \ref{pro 2.1} provides a knowledge of a
cross section of the desired object.

\begin{remark}[A non uniqueness example and the uniqueness of Problem 
\protect\ref{pro 2.1}]
It is not hard to verify that $f(x,y)=h(y)$, for some function $h(y)$
satisfying \[\displaystyle\int_{a}^{b}h(y)dy=0,\] is in the null space of the
\textquotedblleft incomplete" Radon transform whose domain is all pairs $(%
\mathbf{x},\mathbf{x}_{\alpha })$ where $\mathbf{x}$ is on the top of $%
\Omega $ and the line $L(\mathbf{x},\mathbf{x}_{\alpha })$ does not
intersect the vertical sides of $\Omega $. Hence, the knowledge of data for
Problem \ref{pro 2.1} on the vertical sides of $\Omega $ is crucial. The
uniqueness of Problem \ref{pro 2.1} is considered as an assumption in this
paper. On the other hand, we consider in this paper an approximate
mathematical model, which is obtained via the truncation of a certain
Fourier series. Uniqueness for the latter case follows immediately from our
convergence result (Theorem \ref{thm convergence f}).
\end{remark}

\begin{lemma}
Assume that the function $f\in C^{k}(\mathbb{R}^{2})$, $k\geq 1$ and
satisfies condition \eqref{2.1}. Then, the function $u(\mathbf{x},\mathbf{x}%
_{\alpha })=Rf(\mathbf{x},\mathbf{x}_{\alpha })$ is $k$ times continuously
differentiable with respect to both $\mathbf{x}\in \Omega $ and $\alpha \in
(-d,d)$. Moreover, those derivatives are bounded in $\overline{\Omega }%
\times \lbrack -d,d].$ \label{lem 2.1}
\end{lemma}

\textbf{Proof}. 
We have 
\begin{align*}
L(\mathbf{x},\mathbf{x}_{\alpha }) &=\left\{ \mathbf{x}\left( t\right)
=\left( 1-t\right) \mathbf{x}_{\alpha }+t\mathbf{x,} \quad t\in \left(
0,1\right) \right\} \\
&=\left\{ \left( x(t), y(t)\right) :x\left( t\right) =\alpha +t\left(
x-\alpha \right) ,y\left( t\right) =ty, \quad t\in \left( 0,1\right)
\right\}.
\end{align*}
Hence, by (\ref{2.2}) 
\begin{equation*}
u(\mathbf{x},\mathbf{x}_{\alpha })=\int_{L(\mathbf{x},\mathbf{x}%
_{\alpha })}f(\mathbf{\xi })d\sigma =\sqrt{\left( x-\alpha \right) ^{2}+y^{2}%
}\int_{0}^{1}f\left( \alpha +t\left( x-\alpha \right) ,ty\right) dt.%
\text{ }\square
\end{equation*}


\section{An approximation for the model governing the X-ray tomographic data}

\label{sec 4}

We establish in this section a system of first order PDEs that leads to our
numerical method to solve Problem \ref{pro 2.1}.

\subsection{The exact PDE governing the X-ray tomographic function}

\label{sec:3}

For each source $\mathbf{x}_{\alpha }=(\alpha ,0)$ in $\Gamma _{d}$ and $%
\mathbf{x}=(x,y)$ in $\Omega $, let $\varphi $ be the angle constituted by
the line $L(\mathbf{x},\mathbf{x}_{\alpha })$ and the $x-$axis. The
directional derivative of $u(\mathbf{x},\mathbf{x}_{\alpha })$ with respect
to the direction $(\cos \varphi ,\sin \varphi )$ of the line $L(\mathbf{x},%
\mathbf{x}_{0})$ is given by 
\begin{align*}
\cos \varphi u_{x}+\sin \varphi u_{y}& =\lim_{t\rightarrow 0}\frac{u(x+t\cos
\varphi ,y+t\sin (\varphi ),\mathbf{x}_{\alpha })-u(x,y,\mathbf{x}_{\alpha })%
}{t} \\
& =\lim_{t\rightarrow 0}\frac{1}{t}\int_{l_{t}}f(\sigma )d\sigma
\end{align*}%
where $l_{t}\subset L(\mathbf{x},\mathbf{x}_{0})$ is the line connecting the
point $\mathbf{x}$ and $(\mathbf{x}+t(\cos \varphi ,\sin \varphi ))$. Since
the length of $l_{t}$ is $t$ and the function $f$ is continuous, then the
above limit is $f(\mathbf{x})$. Since 
\begin{equation*}
\cos \varphi =\frac{x-\alpha }{|\mathbf{x}-\mathbf{x}_{\alpha }|}=\frac{%
x-\alpha }{\sqrt{\left( x-\alpha \right) ^{2}+y^{2}}}\text{ and }\sin
\varphi =\frac{y}{|\mathbf{x}-\mathbf{x}_{\alpha }|}=\frac{y}{\sqrt{\left(
x-\alpha \right) ^{2}+y^{2}}},
\end{equation*}%
for each $\mathbf{x}_{\alpha }=(\alpha ,0)$, $\alpha \in (-d,d)$, then the
function $u(\mathbf{x},\mathbf{x}_{\alpha })$ satisfies the following form
of the transport equation: 
\begin{equation}
\frac{x-\alpha }{|\mathbf{x}-\mathbf{x}_{\alpha }|}u_{x}+\frac{y}{|\mathbf{x}%
-\mathbf{x}_{\alpha }|}u_{y}=f(x,y).  \label{main eqn}
\end{equation}

\begin{remark}
Although equation \eqref{main eqn} is well known, see e.g., \cite{HasanogluRomanov:s2017}, we have briefly derived it as above. 
This is
because equation \eqref{main eqn} leads us to a PDE approach to
solve the tomographic inverse problem with incomplete data, Problem \ref{pro
2.1}. Equation \eqref{main eqn} is the exact mathematical model that governs
the function $u$.
\end{remark}

\subsection{An orthonormal basis in $L^{2}(-d,d)$}

\label{sec 4.1}

We start by recalling a special orthonormal basis of $L^{2}(-d,d)$ that is
different from the basis constructed from either standard orthonormal
polynomials or trigonometric functions. If one considers such a usual basis
of $L^{2}(-d,d)$, then one of the elements of that basis is a constant,
meaning that its derivative is identically zero. Unlike this, for our
approach, we need to construct an orthonormal basis $\left\{ \Psi _{n}\left(
\alpha \right) \right\} _{n=1}^{\infty }$ in $L_{2}(-d,d),$ which has the
following two properties:

\begin{enumerate}
\item $\Psi _{n}\in C^{1}\left[ -d,d\right] ,$ $\forall n= 1, 2, \dots$

\item Let $\left( ,\right) $ denotes the scalar product in $L^{2}\left(
-d,d\right) $ and let $a_{mn}=\left( \Psi _{n}^{\prime },\Psi _{m}\right) .$
Then the matrix $M_{N}=\left( a_{mn}\right) _{m,n=1}^{N}$ should be
invertible for any $N=1,2,\dots$
\end{enumerate}

Such a basis was first constructed in \cite{Klibanov:jiip2017}. We now
reproduce that construction for the convenience of the reader. For $\alpha
\in (-d,d)$, consider the set of functions $\{\alpha ^{n-1}e^{\alpha
}\}_{n=1}^{\infty }$. These functions are linearly independent and form a
complete set in $L^{2}(-d,d)$. Applying the classical Gram-Schmidt
orthonormalization procedure to this set, we obtain the orthonormal basis $%
\{\Psi _{n}\left( \alpha \right) \}_{n=1}^{\infty }$ of $L^{2}(-d,d)$. It is
obvious that for each $n\geq 1$, $\Psi _{n}(\alpha )=P_{n-1}(\alpha
)e^{\alpha }$, where $P_{n-1}$ is a polynomial of the degree $n-1$. The
following lemma holds true:

\begin{lemma}[Theorem 2.1 in \protect\cite{Klibanov:jiip2017}]
The function $\Psi _{n}^{\prime }$ is not identically zero for any $n\geq 1$%
. Moreover, we have 
\begin{equation*}
\phi _{mn}=\int_{-d}^{d}\Psi _{n}^{\prime }(\alpha )\Psi _{m}(\alpha
)d\alpha =\left\{ 
\begin{array}{ll}
1 & \mbox{if }n=m, \\ 
0 & \mbox{if }n<m.%
\end{array}%
\right.
\end{equation*}%
Consequently, for any integer $N\geq 1$, the matrix $M_{N}=(\phi
_{mn})_{m,n=1}^{N}$ has determinant 1 and is, therefore, invertible. \label%
{lem 3.1}
\end{lemma}


The function $u(x,y,\alpha )$ can be represented via the following Fourier
series, which converges in $L^{2}(-d,d)$ for every point $\left( x,y\right)
\in \overline{\Omega }:$%
\begin{equation*}
u(x,y,\alpha )=\sum_{n=1}^{\infty }u_{n}(x,y)\Psi _{n}(\alpha ),\quad
(x,y)\in \overline{\Omega },\alpha \in (-d,d).
\end{equation*}%
In order to introduce our approximate mathematical model mentioned in
Introduction, we approximate the function $u(\mathbf{x},\mathbf{x}_{\alpha
})=u(x,y,\alpha )$ as: 
\begin{equation}
u(x,y,\alpha )\approx \sum_{n=1}^{N}u_{n}(x,y)\Psi _{n}(\alpha ),\quad
(x,y)\in \overline{\Omega },\alpha \in (-d,d),  \label{4.2}
\end{equation}%
where $N\geq 1$ is a certain integer, which is chosen later numerically, and 
\begin{equation}
u_{n}(x,y)=\int_{-d}^{d}u(x,y,\alpha )\Psi _{n}(\alpha )d\alpha ,\quad
(x,y)\in \overline{\Omega },\alpha \in (-d,d).  \label{eqn find un}
\end{equation}%
Our approximate mathematical model mentioned in the Introduction amounts to
the replacement in (\ref{4.2}) \textquotedblleft $\approx $" with
\textquotedblleft $=$" as well as to the assumption that the resulting
function solves equation (\ref{main eqn}). Thus, everywhere below%
\begin{equation}
u(x,y,\alpha )=\sum_{n=1}^{N}u_{n}(x,y)\Psi _{n}(\alpha ),\text{ \ }(x,y)\in 
\overline{\Omega },\alpha \in (-d,d).  \label{4.3}
\end{equation}

\subsection{A system of first order PDEs}

\label{sec 4.2}

The goal of this section is to derive a system of linear coupled PDEs, whose
solution directly yields numerical solution to Problem \ref{pro 2.1}.
Differentiating equation \eqref{main eqn} with respect to $\alpha $ and
denoting $v=\partial _{\alpha }u$, we obtain 
\begin{equation}
\frac{x-\alpha }{\left\vert \mathbf{x}-\mathbf{x}_{\alpha }\right\vert }%
v_{x}+\frac{y}{\left\vert \mathbf{x}-\mathbf{x}_{\alpha }\right\vert }v_{y}-%
\frac{y^{2}}{\left\vert \mathbf{x}-\mathbf{x}_{\alpha }\right\vert ^{3}}%
u_{x}+\frac{(x-\alpha )y}{\left\vert \mathbf{x}-\mathbf{x}_{\alpha
}\right\vert ^{3}}u_{y}=0
 \label{4.40}
\end{equation}
for all $\x = (x, y) \in \Omega$ and $\alpha \in (-d, d).$
Equation (\ref{4.40}) is equivalent with 
\begin{equation}
v_{y}=-\frac{x-\alpha }{\left\vert \mathbf{x}-\mathbf{x}_{\alpha
}\right\vert ^{2}}u_{y}-\frac{x-\alpha }{y}v_{x}+\frac{y}{\left\vert \mathbf{%
x}-\mathbf{x}_{\alpha }\right\vert ^{2}}u_{x},\text{ \ }\mathbf{x}=(x,y)\in
\Omega ,\alpha \in (-d,d).  \label{4.4}
\end{equation}%
By (\ref{4.3}) the function $v(x,y,\alpha )$ can be written as 
\begin{equation}
v(\mathbf{x},\alpha )=\sum_{n=1}^{N}u_{n}(\mathbf{x})\Psi _{n}^{\prime
}(\alpha ),\quad \mathbf{x}=(x,y)\in \overline{\Omega },\alpha \in (-d,d).
\label{4.5}
\end{equation}%
Plugging the function $u$ and $v$ in (\ref{4.3}) and (\ref{4.5})
respectively into (\ref{4.4}), we obtain%
\begin{multline}
\sum_{n=1}^{N}\partial _{y}u_{n}(\mathbf{x})\Psi _{n}^{\prime }(\alpha )=-%
\frac{x-\alpha }{\left\vert \mathbf{x}-\mathbf{x}_{\alpha }\right\vert ^{2}}%
\sum_{n=1}^{N}\partial _{y}u_{n}(\mathbf{x})\Psi _{n}(\alpha )
\\
-\frac{x-\alpha }{y}\sum_{n=1}^{N}\partial _{x}u_{n}(\mathbf{x})\Psi
_{n}^{\prime }(\alpha )+\frac{y}{\left\vert \mathbf{x}-\mathbf{x}_{\alpha
}\right\vert ^{2}}\sum_{n=1}^{N}\partial _{x}u_{n}(\mathbf{x})\Psi
_{n}(\alpha )
\label{4.6}
\end{multline}
for all $\x = (x, y) \in \Omega$ and $\alpha \in (-d, d).$
Multiplying both sides of (\ref{4.6}) by $\Psi _{m}(\alpha )$, $m\in
\{1,\dots ,N\}$, and then integrating the resulting equation with respect to 
$\alpha \in (-d,d)$, we obtain 
\begin{multline*}
\sum_{n=1}^{N}\partial _{y}u_{n}(\mathbf{x})\int_{-d}^{d}\Psi _{m}(\alpha
)\Psi _{n}^{\prime }(\alpha )d\alpha =-\sum_{n=1}^{N}\partial _{y}u_{n}(%
\mathbf{x})\int_{-d}^{d}\frac{x-\alpha }{\left\vert \mathbf{x}-\mathbf{x}%
_{\alpha }\right\vert ^{2}}\Psi _{m}(\alpha )\Psi _{n}(\alpha )d\alpha  \\
+\sum_{n=1}^{N}\partial _{x}u_{n}(\mathbf{x})\int_{-d}^{d}\left( -\frac{%
x-\alpha }{y}\Psi _{m}(\alpha )\Psi _{n}^{\prime }(\alpha )+\frac{y}{%
\left\vert \mathbf{x}-\mathbf{x}_{\alpha }\right\vert ^{2}}\Psi _{m}(\alpha
)\Psi _{n}(\alpha )\right) d\alpha .
\end{multline*}%
Recalling Lemma \ref{lem 3.1}, we obtain 
\begin{equation}
M_{N}\mathbf{U}_{y}(\mathbf{x})=D_{1}\left( \mathbf{x}\right) \mathbf{U}_{y}(%
\mathbf{x})+D_{2}\left( \mathbf{x}\right) \mathbf{U}_{x}(\mathbf{x}),\quad 
\mathbf{x}=(x,y)\in \Omega ,  \label{4.7}
\end{equation}%
where the $N-$ dimensional vector valued function $\mathbf{U}\left( \mathbf{x%
}\right) $ is 
\begin{equation}
\mathbf{U}\left( \mathbf{x}\right) =\left( u_{1},\dots,u_{N}\right) ^{T}\left( 
\mathbf{x}\right)   \label{4.70}
\end{equation}%
and $D_{1}\left( \mathbf{x}\right) ,D_{2}\left( \mathbf{x}\right) $ are two $%
N\times N$ matrices whose $mn^{\mathrm{th}}$, $1\leq m,n\leq N$, entries 
\begin{align}
(D_{1})_{mn}&=\int_{-d}^{d}\frac{x-\alpha }{\left\vert \mathbf{x}-\mathbf{x}%
_{\alpha }\right\vert ^{2}}\Psi _{m}(\alpha )\Psi _{n}(\alpha )d\alpha ,
\label{4.71}
\\
(D_{2})_{mn}&=\int_{-d}^{d}\left( -\frac{x-\alpha }{y}\Psi _{m}(\alpha )\Psi
_{n}^{\prime }(\alpha )+\frac{y}{\left\vert \mathbf{x}-\mathbf{x}_{\alpha
}\right\vert ^{2}}\Psi _{m}(\alpha )\Psi _{n}(\alpha )\right) d\alpha 
\label{4.72}
\end{align}%
belong to $C^{\infty }\left( \overline{\Omega }\right) .$ The following
lemma follows immediately from (\ref{4.71}) and (\ref{4.72}):

\begin{lemma}
Suppose that in the definition (\ref{1}) of the domain $\Omega $ the number $%
a>1.$ Then the following estimates hold: 
\begin{equation*}
\max_{\mathbf{x}\in \overline{\Omega }}\left\Vert D_{1}\left( \mathbf{x}%
\right) \right\Vert \leq \frac{C_{1}}{a^{2}},\text{ }\max_{\mathbf{x}\in 
\overline{\Omega }}\left\Vert D_{2}\left( \mathbf{x}\right) \right\Vert \leq 
\frac{C_{1}}{a}.
\end{equation*}%
\label{lem 3.2}
\end{lemma}

\begin{remark}
The system of PDEs \eqref{4.7} for the $N$-dimensional vector valued
function $\mathbf{U}(\mathbf{x})$ is our approximate mathematical model for
the exact one \eqref{main eqn}. Our method to solve Problem \ref{pro 2.1} is
based on a numerical solver for \eqref{4.7}. \label{rem 3.2}
\end{remark}

Here and everywhere below the norm of a matrix is the square root of the sum
of squares of its entries. Also, in Lemma \ref{lem 3.3} and everywhere below 
$C_{1}=C_{1}\left( N,R,d\right) >0$ denotes different constants independent
on the number $a$. Rewrite (\ref{4.7}) as%
\begin{equation}
M_{N}\left( I-M_{N}^{-1}D_{1}\left( \mathbf{x}\right) \right) \mathbf{U}_{y}(%
\mathbf{x})+D_{2}\left( \mathbf{x}\right) \mathbf{U}_{x}(\mathbf{x})=0, \quad \mathbf{x}=(x,y)\in \Omega .  \label{4.8}
\end{equation}

Lemma \ref{lem 3.3} follows immediately from Lemmata \ref{lem 3.1} and \ref%
{lem 3.2}.

\begin{lemma}
For each $N\geq 1$, there exists a sufficiently large number $%
a_{0}=a_{0}(N,R,d)>1$ such that for any $a\geq a_{0},$ the matrix $%
M_{N}\left( I-M_{N}^{-1}D_{1}\left( \mathbf{x}\right) \right) $ is
invertible. Denote 
\[
	D\left( \mathbf{x}\right) =-\left[ M_{N}\left(
I-M_{N}^{-1}D_{1}\left( \mathbf{x}\right) \right) \right] ^{-1}D_{2}\left( 
\mathbf{x}\right) \quad \x \in \overline \Omega.
\]
We have 
\begin{equation}
\max_{\mathbf{x}\in \overline{\Omega }}\left\Vert D\left( \mathbf{x}\right)
\right\Vert \leq C_{1}.  \label{4.80}
\end{equation}%
Moreover, equation (\ref{4.8}) is equivalent to 
\begin{equation}
\mathbf{U}_{y}(\mathbf{x})+D\left( \mathbf{x}\right) \mathbf{U}_{x}(\mathbf{x%
})=0,\text{ }\mathbf{x}=(x,y)\in \Omega .  \label{4.9}
\end{equation}%
\label{lem 3.3}
\end{lemma}

In addition to (\ref{4.9}), the following vector function $\mathbf{g}\left( 
\mathbf{x}\right) $ of boundary conditions is known%
\begin{equation}
\mathbf{U}(\mathbf{x})=\mathbf{g}(\mathbf{x}),\quad \mathbf{x}\in \partial
\Omega  \label{4.10}
\end{equation}%
via using \eqref{eqn find un} and \eqref{4.70} for $\mathbf{x}=(x,y)\in
\partial \Omega $. In particular 
\begin{equation}
\mathbf{g}\left( \mathbf{x}\right) =0\text{ for }\mathbf{x}=(x,a).
\label{4.100}
\end{equation}
Thus, we solve below boundary value problem (\ref{4.9}), (\ref{4.10}).
Suppose that we have obtained its approximate solution. Then the
corresponding approximation for the target function $f\left( \mathbf{x}%
\right) $ should be obtained via the substitution of (\ref{4.3}) in (\ref%
{main eqn}), see (\ref{5.32}).

\begin{remark}
\label{rem 4.1} As it was mentioned in Introduction, the number $N$ should
be chosen numerically, also see Remark \ref{rem 5.1}.
\end{remark}

\section{The Quasi-Reversibility Method for the first order system of PDEs 
\eqref{4.9}--\eqref{4.10}}

\label{sec 5}

The boundary value problem (\ref{4.9}), (\ref{4.10}) is overdetermined since
the boundary data (\ref{4.10}) for the system (\ref{4.9}) of PDEs of the
first order are given on the whole boundary $\partial \Omega $ rather than
on its part. Therefore, to find an approximate solution of problem (\ref{4.9}%
), (\ref{4.10}), we use the quasi-reversibility method, which, in general,
works properly for overdetermined problems.

For vector functions $\mathbf{U\in }H^{1}\left( \Omega \right)^N,$ consider
the functional $J_{\varepsilon }\left( \mathbf{U}\right) $%
\begin{equation}
J_{\varepsilon }\left( \mathbf{U}\right) =\int_{\Omega }\left| 
\mathbf{U}_{y}(\mathbf{x})+D\left( \mathbf{x}\right) \mathbf{U}_{x}(\mathbf{x%
})\right|^{2}d\mathbf{x} + \varepsilon \left\Vert \mathbf{U}\right\Vert
_{H^{1}(\Omega)^N}^{2},  \label{5.1}
\end{equation}%
where $\varepsilon \in \left( 0,1\right) $ is the regularization parameter.
The quasi-reversibility method for problem (\ref{4.9})--(\ref{4.10}) amounts
to the following minimization problem:

\begin{problem*}[Solving \eqref{4.9}--\eqref{4.10} by the
quasi-reversibility method in the continuous form]
Minimize functional (\ref{5.1}) on the set of vector functions $\mathbf{U}%
\in H^{1}\left( \Omega \right) ^{N},$ subject to boundary condition %
\eqref{4.10}. The resulting minimizer is called the regularized solution of %
\eqref{4.9}--\eqref{4.10}.
\end{problem*}

Conventionally, the convergence analysis of the quasi-reversibility method
is performed on the basis of Carleman estimates \cite{Ksurvey}. However,
since $\mathbf{U}\left( \mathbf{x}\right) $ is a vector function rather than
a 1D function and also since the matrix $D\left( \mathbf{x}\right) $ is
likely not self adjoint, we cannot currently derive a proper Carleman
estimate for the differential operator in the integrand of the right hand
side of (\ref{5.1}). Hence, we consider this operator in its semi discrete
form, assuming the finite differences in the $x-$direction. However, we do
not \textquotedblleft allow" the step size $h$ of the finite difference tend
to zero and, do not estimate the distance between the finite difference and
continuous solutions. We observe that the semi discrete form is more
realistic for computations than the continuous form. In our numerical
realization we consider the fully discrete form, see Section \ref{sec imple}%
. As it is often the case in the field of ill-posed and inverse problems,
the theory for the fully discrete case is more complicated, see, e.g. \cite%
{KS}. Thus, it is outside of the scope of this first publication about our
method. is not yet developed.

\subsection{ Semi discrete formulation of the quasi-reversibility method}

\label{sec 5.1}

Let the number $h_{0}\in \left( 0,1\right) .$ We assume that there exists a
number $h\in \left[ h_{0},1\right) $ such that the number $K=2R/h$ is an
integer. When saying below \textquotedblleft for all $h\in \left[
h_{0},1\right) ",$ we mean only those number $h$ for which the number $2R/h$
is an integer. In any case, let $h$ be one of such numbers. In the interval $%
x\in \left[ -R,R\right] $, consider the grid of the finite difference scheme
with the step size $h$,%
\[
	x_0 = -R < x_1 = -R + h < \dots < x_i = -R + i h < \dots  < x_K -R + Kh = R.
\]
We define the domain $\Omega ^{h}$ as 
\begin{equation}
\Omega ^{h}=\left\{ \mathbf{x}=\left( x,y\right)
:x=x_{i}=-R+ih,i=1,\dots,\left( K-1\right) ;y\in \left( a,b\right) \right\} .
\label{2}
\end{equation}%
For any $N-$dimensional vector function $Q\left( \mathbf{x}\right) \in
C\left( \overline{\Omega }\right) $,  denote 
\begin{align}
\mathbf{Q}_{i}^{h}\left( y\right) &=\mathbf{Q}\left( -R+ih,y\right) \quad
i=0,\dots, K,y\in (a, b) ,  \label{3}
\\
\mathbf{Q}^{h}\left( y\right) &=\left( \mathbf{Q}_{1}^{h}\left( y\right) ,\dots,%
\mathbf{Q}_{K-1}^{h}\left( y\right) \right) ^{T} \quad y\in  (a, b),  \label{5.2}
\\
\widetilde{\mathbf{Q}}^{h}\left( y\right) &=\left( \mathbf{Q}_{0}^{h}\left(
y\right) ,\mathbf{Q}_{1}^{h}\left( y\right) ,\dots,\mathbf{Q}_{K-1}^{h}\left(
y\right) ,\mathbf{Q}_{K}^{h}\left( y\right) \right) ^{T} \nonumber
\\
&=\left( \mathbf{Q}%
_{0}^{h}\left( y\right), \mathbf{Q}^{h}\left( y\right) , \dots, 
\mathbf{Q}_{K}^{h}\left( y\right) \right) ^{T} \quad y\in (a, b).
\label{30}
\end{align}
Note that, unlike $\widetilde{\mathbf{Q}}^{h}\left( y\right) ,$ vector
functions in (\ref{5.2}) do not include boundary terms 
\begin{equation}
\mathbf{Q}_{0}^{h}\left( y\right) =\mathbf{Q}\left( -R,y\right) ,\mathbf{Q}%
_{K}^{h}\left( y\right) =\mathbf{Q}\left( R,y\right), \quad y\in (a, b)
\label{6}
\end{equation}%
at the vertical sides of the rectangle $\Omega $ in (\ref{1}). Since $%
\mathbf{Q}\left( \mathbf{x}\right) $ is an $N-$dimensional vector valued
function, then $\mathbf{Q}^{h}\left( y\right) $ and $\widetilde{\mathbf{Q}}%
^{h}\left( y\right) $ are $N\times \left( K-1\right) $ and $N\times \left(
K+1\right) $ respectively matrix valued functions of the variable $y$. Let $%
D^{h}\left( y\right) $ be the block diagonal matrix, whose block matrices on the diagonal are $%
K-1$ sub-matrices of the form $D\left( -R+h,y\right) ,\dots,D\left( -R+\left(
K-1\right) h,y\right) ,y\in \left[ a,b\right] .$ It follows from (\ref{4.80}%
) that%
\begin{equation}
\max_{y\in \left[ a,b\right] }\left\Vert D^{h}\left( y\right) \right\Vert
\leq C_{1}.  \label{5.200}
\end{equation}

We now come back to our vector function $\mathbf{U}\left( \mathbf{x}\right)
. $ Using definition (\ref{5.2}), we approximate the derivative $\mathbf{U}%
_{x} $ at the point $\left( -R+jh,y\right) \in \Omega ^{h}$ by the central
finite difference as%
\begin{equation}
\mathbf{U}_{jx}^{h}\left( y\right) =\frac{\mathbf{U}_{j+1}^{h}\left(
y\right) -\mathbf{U}_{j-1}^{h}\left( y\right) }{2h} \quad j=1,\dots,K-1.
\label{5.3}
\end{equation}%
By Lemma \ref{lem 2.1}, $\mathbf{U}\in C^{2}\left( \overline{\Omega }\right) 
$. Hence, it follows from (\ref{5.3}) that 
\begin{equation*}
\mathbf{U}_{jx}^{h}\left( y\right) =\mathbf{U}_{jx}\left( -R+jh,y\right)
+O\left( h\right) \text{ as }h\rightarrow 0;\text{ }j=1,\dots,K-1.
\end{equation*}%
Denote 
\begin{equation*}
\mathbf{U}_{x}^{h}\left( y\right) =\left( \mathbf{U}_{1x}^{h}\left(
-R+h,y\right) ,\dots,\mathbf{U}_{\left( K-1\right) x}^{h}\left( -R+\left(
K-1\right) h,y\right) \right) ^{T}.
\end{equation*}%
Hence, dropping $O\left( h\right) ,$ we obtain the following finite
difference analog of problem (\ref{4.9}), (\ref{4.10})%
\begin{equation}
\left\{ 
\begin{array}{ll}
\mathbf{U}_{y}^{h}(y)+D^{h}\left( y\right) \mathbf{U}_{x}^{h}(y)=\mathbf{0},
&\text{ in }\Omega ^{h}, \\ 
\mathbf{U}^{h}\left( a\right) =\mathbf{g}^{h}(a)=\mathbf{0},
\mathbf{U}%
^{h}\left( b\right) =\mathbf{g}^{h}(b), 
\\
\mathbf{U}_{0}^{h}\left( y\right) =\mathbf{g}\left( -R,y\right) ,\mathbf{U}%
_{K}^{h}\left( y\right) =\mathbf{g}\left( R,y\right) &y\in \left[ a,b\right]
.%
\end{array}%
\right.  \label{5.5}
\end{equation}
where the boundary matrix $\mathbf{g}^{h}\left( b\right) $ is known and is
defined using the vector function $\mathbf{g}\left( \mathbf{x}\right) ,%
\mathbf{x}\in \partial \Omega $ in the obvious manner, also see (\ref{4.100}%
) and (\ref{3})-(\ref{6}).

We now introduce semi discrete functional spaces for matrices $\mathbf{Q}^{h}
$, $\widetilde{\mathbf{Q}}^{h}$defined in (\ref{3})-(\ref{6}). We set 
\begin{align*}
L^{2,h}\left( \Omega ^{h}\right) & =\left\{ \mathbf{U}^{h}\left( y\right)
:\left\Vert \mathbf{U}^{h}\left( y\right) \right\Vert _{L^{2,h}\left( \Omega
^{h}\right) }^{2}=\sum_{i=1}^{K-1}h\int_{a}^{b}\left[ \mathbf{U}%
_{i}^{h}\left( y\right) \right] ^{2}dy<\infty \right\} , \\
H^{1,h}\left( \Omega ^{h}\right) & =\Big\{\mathbf{U}^{h}\left( y\right)
:\left\Vert \mathbf{U}^{h}\left( y\right) \right\Vert _{H^{1,h}\left( \Omega
^{h}\right) }^{2} \\
& =\sum_{j=1}^{K-1}h\int_{a}^{b}\big[\left( \mathbf{U}_{jx}^{h}\left(
-R+jh,y\right) \right) ^{2}+\left( \partial _{y}\mathbf{U}_{j}^{h}\left(
y\right) \right) ^{2}+\left( \mathbf{U}_{j}^{h}\left( y\right) \right) ^{2}%
\big]dy<\infty \Big\}.
\end{align*}%
\begin{equation}
\widetilde{H}^{1,h}\left( \Omega ^{h}\right) =\left\{ 
\begin{array}{ll}
\widetilde{\mathbf{P}}^{h}\left( y\right) =\left( \mathbf{P}_{0}^{h}\left(
y\right) ,\mathbf{P}_{1}^{h}\left( y\right) ,\dots,\mathbf{P}_{K-1}^{h}\left(
y\right) ,\mathbf{P}_{K}^{h}\left( y\right) \right) ^{T}: \\ 
\mathbf{P}^{h}\left( y\right) =\left( \mathbf{P}_{1}^{h}\left( y\right) ,\dots,%
\mathbf{P}_{K-1}^{h}\left( y\right) \right) ^{T}\in H^{1,h}\left( \Omega
^{h}\right) , \\ 
\mathbf{P}^{h}\left( a\right) =\mathbf{P}^{h}\left( b\right) =\mathbf{0,P}%
_{0}^{h}\left( y\right) =\mathbf{P}_{K}^{h}\left( y\right) =\mathbf{0,} 
\\
\left\Vert \widetilde{\mathbf{P}}^{h}\right\Vert _{\widetilde{H}^{1,h}\left(
\Omega ^{h}\right) }=\left\Vert \mathbf{P}^{h}\right\Vert _{H^{1,h}\left(
\Omega ^{h}\right) }.%
\end{array}%
\right.   
y \in [a, b]
\label{5.50}
\end{equation}%
also see (\ref{6}). Scalar products in these spaces are defined in the
obvious manner. We denote the scalar product in the space $L^{2,h}(\Omega
^{h})$ as $\left( ,\right) ^{h}$ and the one for the latter two Sobolev
spaces as $\left[ ,\right] ^{h}$. Below we fix the number $h_{0}\in \left(
0,1\right) .$ It follows from (\ref{5.3}) that there exists a constant $%
B_{h_{0}}=B_{h_{0}}\left( h_{0}\right) >0$ depending only on $h_{0}$ such
that%
\begin{equation}
\left\Vert \mathbf{Q}_{x}^{h}\left( y\right) \right\Vert _{L^{2,h}\left(
\Omega ^{h}\right) }^{2}\leq B_{h_{0}}\left\Vert \mathbf{Q}^{h}\right\Vert
_{L^{2,h}\left( \Omega ^{h}\right) }^{2},\text{ }\forall \mathbf{Q}^{h}:%
\widetilde{\mathbf{Q}}^{h}\in \widetilde{H}^{1,h}\left( \Omega ^{h}\right)
,\forall h\in \left[ h_{0},1\right] .  \label{7}
\end{equation}%

\begin{remark}
 Thus, according to (\ref{30}), (\ref{6}) and (\ref{7}), if a matrix 
$\mathbf{Z}^{h}$ is defined on the set $\Omega ^{h},$ then $\widetilde{%
\mathbf{Z}}^{h}$ means that this matrix is complemented by boundary
conditions at $x=-R$, $x=R$, $y=a,y=b$. In particular,  $\widetilde{\mathbf{Z%
}}^{h}\in \widetilde{H}^{1,h}\left( \Omega ^{h}\right) $ means that $\mathbf{%
Z}^{h}\in H^{1,h}\left( \Omega ^{h}\right) $ and those boundary conditions
are zeros.
\label{rem bdry}
\end{remark}

The semi discrete quasi-reverisibility method applied to problem \eqref{5.5}
is:

\begin{problem*}[Solving \eqref{5.5} by the quasi-reversibility method in
the semi discrete form]
Let $\varepsilon \in \left[ 0,1\right) $ be the regularization parameter.
Minimize the functional $J_{\varepsilon }^{h}\left( \mathbf{U}^{h}\right) ,$ 
\begin{equation}
J_{\varepsilon }^{h}\left( \widetilde{\mathbf{U}}^{h}\right) =\left\Vert 
\mathbf{U}_{y}^{h}\left( y\right) +D^{h}\left( y\right) \mathbf{U}%
_{x}^{h}\left( y\right) \right\Vert _{L^{2,h}\left( \Omega ^{h}\right)
}^{2}+\varepsilon \left\Vert \mathbf{U}^{h}\right\Vert _{H^{1,h}\left(
\Omega ^{h}\right) }^{2}  \label{5}
\end{equation}%
on the set of matrices $\widetilde{\mathbf{U}}^{h}$ such that boundary
conditions of \eqref{5.5} are satisfied, also see Remark \ref{rem bdry} for $\mathbf{U}^{h}$ and $\widetilde{\mathbf{U}}^{h}$.
 \label{pro 2}
\end{problem*}

\subsection{Existence and uniqueness of the solution of Problem \protect\ref%
{pro 2}}

\label{sec 5.2}

First, we prove a new Carleman estimate:

\begin{lemma}[Carleman estimate]
Let the parameter $\lambda >0.$ The following Carleman estimate holds true 
\begin{equation}
\int_{a}^{b}\left( w^{\prime }\right) ^{2}e^{2\lambda y}dy\geq \frac{1}{2}%
\int_{a}^{b}\left( w^{\prime }\right) ^{2}e^{2\lambda y}dy+\frac{1}{2%
}\lambda ^{2}\int_{a}^{b}w^{2}e^{2\lambda y}dy,\text{ }\forall w\in 
\widetilde{H}^{1}\left( a,b\right) ,\forall \lambda >0.  \label{4}
\end{equation}%
Here, $\widetilde{H}^{1}\left( a,b\right) $ is the subspace of functions $%
w\in H^{1}\left( a,b\right) $ satisfying $w\left( b\right) =0.$ \label{lem
Carleman}
\end{lemma}

Note that usually a generic constant $C>0$ is used in Carleman estimates,
see, e.g. Chapter 4 in \cite{LRS}. In (\ref{4}), however, we have a specific
value $C=1/2.$

\textbf{Proof of Lemma \ref{lem Carleman}}. Introduce a new function $%
p\left( y\right) =w\left( y\right) e^{\lambda y}.$ Then $w\left( y\right)
=p\left( y\right) e^{-\lambda y}.$ Hence, $w^{\prime }=p^{\prime
}e^{-\lambda y}-\lambda pe^{-\lambda y}.$ We have 
\begin{align*}
\left( w^{\prime }\right) ^{2}e^{2\lambda y}&=\left( p^{\prime }-\lambda
p\right) ^{2}=\left( p^{\prime }\right) ^{2}-2\lambda p^{\prime }p+\lambda
^{2}p^{2}
\\
&\geq -2\lambda p^{\prime }p+\lambda ^{2}p^{2}=\left( -\lambda p^{2}\right)
^{\prime }+\lambda ^{2}p^{2}=\left( -\lambda w^{2}e^{2\lambda y}\right)
^{\prime }+\lambda ^{2}w^{2}e^{2\lambda y}.
\end{align*}%
Hence,%
\begin{equation*}
\int_{a}^{b}\left( w^{\prime }\right) ^{2}e^{2\lambda y}dy\geq -\lambda
w^{2}\left( b\right) e^{2\lambda b}+\lambda w^{2}\left( a\right) e^{2\lambda
a}+\lambda ^{2}\int_{a}^{b}w^{2}e^{2\lambda y}dy\geq \lambda
^{2}\int_{a}^{b}w^{2}e^{2\lambda y}dy.
\end{equation*}%
%
%
%
Therefore, 
\begin{equation*}
2\int_{a}^{b}\left( w^{\prime }\right) ^{2}e^{2\lambda y}dy\geq
\int_{a}^{b}\left( w^{\prime }\right) ^{2}e^{2\lambda y}dy+\lambda
^{2}\int_{a}^{b}w^{2}e^{2\lambda y}dy.
\end{equation*}%
Dividing this inequality by 2, we obtain (\ref{4}). $\square $

\begin{theorem}
Assume that $a\geq a_{0}=a_{0}\left( N,R,d\right) >1,$ where $a_{0}\left(
N,R,d\right) $ is the number defined in Lemma \ref{lem 3.3}. Also, assume
that functions $\mathbf{g}\left( -R,y\right) ,\mathbf{g}\left( R,y\right)
\in C^{1}\left[ a,b\right] .$ Suppose that there exists an $N\times \left(
K+1\right) $ matrix 
\[\widetilde{\mathbf{F}}^{h}\left( y\right) =\left( 
\mathbf{F}_{0}^{h}\left( y\right) ,\mathbf{F}_{1}^{h}\left( y\right) ,\dots,%
\mathbf{F}_{K-1}^{h}\left( y\right) ,\mathbf{F}_{K}^{h}\left( y\right)
\right) ^{T}
\] 
such that $\mathbf{F}^{h}\left( y\right) =\left( \mathbf{F}%
_{1}^{h}\left( y\right) ,\dots,\mathbf{F}_{K-1}^{h}\left( y\right) \right)
^{T}\in H^{1,h}\left( \Omega ^{h}\right) $ and 
\begin{equation*}
\mathbf{F}^{h}\left( a\right) =\mathbf{0},\mathbf{F}^{h}\left( b\right) =%
\mathbf{g}^{h}(b);\text{ }\mathbf{F}_{0}^{h}\left( y\right) =\mathbf{g}%
\left( -R,y\right) ,\mathbf{F}_{K}^{h}\left( y\right) =\mathbf{g}\left(
R,y\right) ,y\in \left[ a,b\right] ,
\end{equation*}%
see (\ref{5.5}). Then for each number $\varepsilon \in \left[ 0,1\right) $
and for each $h\in \left[ h_{0},1\right) $ there exists unique solution $%
\widetilde{\mathbf{U}}_{\min }^{h}\left( y\right) $ with $\mathbf{U}_{\min
}^{h}\left( y\right) \in H^{1,h}\left( \Omega ^{h}\right) $ of the Problem %
\ref{pro 2} (see (\ref{3})-(\ref{6})). Furthermore, there exists a constant $%
C_{h_{0}}=C_{h_{0}}\left( N,\Omega ,d,h_{0}\right) >0$ depending only on
listed parameters such that the following estimate holds: 
\begin{equation}
\left\Vert \mathbf{U}_{\min }^{h}\right\Vert _{H^{1,h}(\Omega ^{h})}\leq
C_{h_{0}}\left\Vert \mathbf{F}^{h}\right\Vert _{H^{1,h}(\Omega ^{h})}.
\label{5.11}
\end{equation}%
\label{thm existence}
\end{theorem}

\textbf{Proof}. Everywhere below $C_{h_{0}}=C_{h_{0}}\left( N,\Omega
,d,h_{0}\right) >0$\emph{\ }denotes different constants depending only on
listed parameters.\emph{\ }Consider the matrix $\widetilde{\mathbf{V}}%
^{h}\in \widetilde{H}^{1,h}\left( \Omega ^{h}\right) $ defined as $%
\widetilde{\mathbf{V}}^{h}=\widetilde{\mathbf{U}}^{h}-\widetilde{\mathbf{F}}%
^{h}.$ Hence, the functional functional $J_{\varepsilon }^{h}$ defined in (\ref{5}) becomes the functional $I_{\varepsilon }^{h}( 
\mathbf{V}^{h}),$ where 
\begin{multline}
I_{\varepsilon }^{h}(\widetilde{\mathbf{V}}^{h}) =J_{\varepsilon
}^{h}(\widetilde{\mathbf{V}}^{h}+\widetilde{\mathbf{F}}^{h})  
=\left\Vert \mathbf{V}_{y}^{h}\left( y\right) +D^{h}\left( y\right) \mathbf{V%
}_{x}^{h}\left( y\right) +\mathbf{G}^{h}\left( y\right) \right\Vert
_{L^{2,h}\left( \Omega ^{h}\right) }^{2}
\\
+\varepsilon \left\Vert \mathbf{V}%
^{h}+\mathbf{F}^{h}\right\Vert _{H^{1,h}\left( \Omega ^{h}\right) }^{2},
\label{5.120}
\end{multline}
for all $\mathbf{V}^{h}$ such that $\widetilde{\mathbf{V}}^{h}\in \widetilde{H}
^{1,h}\left( \Omega ^{h}\right)$, see Remark \ref{rem bdry}. In (\ref%
{5.120}) 
\begin{equation}
\mathbf{G}^{h}\left( y\right) =\mathbf{F}_{y}^{h}\left( y\right)
+D^{h}\left( y\right) \mathbf{F}_{x}^{h}( y), \quad y\in \left(
a,b\right) . 
 \label{5.13}
\end{equation}%
The matrix $\widetilde{\mathbf{V}}^{h}\in \widetilde{H}^{1,h}\left( \Omega
^{h}\right) $ minimizes the functional (\ref{5.120}) if and only if the
matrix $\widetilde{\mathbf{U}}^{h}=\widetilde{\mathbf{V}}^{h}+\widetilde{%
\mathbf{F}}^{h}$ solves Problem 2$.$ Let $\widetilde{\mathbf{V}}_{\min
}^{h}\in \widetilde{H}^{1,h}\left( \Omega ^{h}\right) $ be a minimizer of
the functional (\ref{5.120}). Then by the variational principle%
\begin{multline}
\left( \mathbf{V}_{\min y}^{h}\left( y\right) +D\left( y\right) \mathbf{V}%
_{\min x}^{h}\left( y\right) ,\mathbf{W}_{y}^{h}\left( y\right) +D^{h}\left(
y\right) \mathbf{W}_{x}^{h}\left( y\right) \right) ^{h}+\varepsilon \left[ 
\mathbf{V}^{h},\mathbf{W}^{h}\right] ^{h}
\\
=-\left( \mathbf{G}^{h}\left( y\right) ,\mathbf{W}_{y}^{h}\left( y\right)
+D^{h}\left( y\right) \mathbf{W}_{x}^{h}\left( y\right) \right)
^{h}-\varepsilon \left[ \mathbf{F}^{h},\mathbf{W}^{h}\right] ^{h},
\label{5.14}
\end{multline}%
for all $\mathbf{W}^{h}$ such that $\widetilde{\mathbf{W}}^{h}\in \widetilde{H}^{1,h}\left( \Omega ^{h}\right)$, see Remark \ref{rem bdry}.
Using (\ref{5.200}), (\ref{7}) and the Cauchy-Schwarz inequality, we obtain 
\begin{align}
\big\Vert\mathbf{W}_{y}^{h}\left( y\right) &+ D^{h}\left( y\right) \mathbf{W}%
_{x}^{h}\left( y\right) \big\Vert_{L^{2,h}\left( \Omega ^{h}\right)
}^{2}=\left\Vert \left[ \mathbf{W}_{y}^{h}\left( y\right) +D^{h}\left(
y\right) \mathbf{W}_{x}^{h}\left( y\right) \right] e^{2\lambda
y}e^{-2\lambda y}\right\Vert _{L^{2,h}\left( \Omega ^{h}\right) }^{2}  \notag
\\
&\geq e^{-2\lambda b}\left\Vert \left[ \mathbf{W}_{y}^{h}\left( y\right)
+D^{h}\left( y\right) \mathbf{W}_{x}^{h}\left( y\right) \right] e^{2\lambda
y}\right\Vert _{L^{2,h}\left( \Omega ^{h}\right) }^{2} \notag
\\
&\geq \frac{1}{2}e^{-2\lambda b}\left\Vert \mathbf{W}_{y}^{h}\left( y\right)
e^{2\lambda y}\right\Vert _{L^{2,h}\left( \Omega ^{h}\right)
}^{2}-C_{h_{0}}e^{-2\lambda b}\left\Vert \mathbf{W}^{h}e^{2\lambda
y}\right\Vert _{L^{2,h}\left( \Omega ^{h}\right) }^{2}, \notag
\end{align}
for all $\mathbf{W}^{h}$ such that $\widetilde{\mathbf{W}}^{h}\in \widetilde{H}^{1,h}(
\Omega ^{h}),$ see Remark \ref{rem bdry}.
Hence, by Lemma \ref{lem Carleman} 
\begin{multline}
\big\Vert\mathbf{W}_{y}^{h}\left( y\right) + D^{h}\left( y\right) \mathbf{W}%
_{x}^{h}\left( y\right) \big\Vert_{L^{2,h}\left( \Omega ^{h}\right) }^{2}
\geq \frac{1}{4}e^{-2\lambda b}\left\Vert \mathbf{W}_{y}^{h}\left( y\right)
e^{2\lambda y}\right\Vert _{L^{2,h}\left( \Omega ^{h}\right) }^{2}
\\+\frac{1}{4%
}e^{-2\lambda b}\lambda ^{2}\left\Vert \mathbf{W}^{h}e^{2\lambda
y}\right\Vert _{L^{2,h}\left( \Omega ^{h}\right) }^{2}  
-C_{h_{0}}e^{-2\lambda b}\left\Vert \mathbf{W}^{h}e^{2\lambda y}\right\Vert
_{L^{2,h}\left( \Omega ^{h}\right) }^{2},
\label{5.140}
\end{multline}
for all $\mathbf{W}^{h}$ such that $\widetilde{\mathbf{W}}^{h}\in \widetilde{H}^{1,h}(
\Omega ^{h}),$ see Remark \ref{rem bdry}.
Recalling (\ref{5.200}), fix a sufficiently large number $\lambda =\lambda
\left( C_{1},h_{0}\right) >\sqrt{8C_{h_{0}}}.$ Then (\ref{5.140}) implies
that 
\begin{align*}
\big\Vert\mathbf{W}_{y}^{h}\left( y\right) &+D^{h}\left( y\right) \mathbf{W}%
_{x}^{h}\left( y\right) \big\Vert_{L^{2,h}\left( \Omega ^{h}\right) }^{2}
\\
&\geq \frac{1}{4}e^{-2\lambda b}\left\Vert \mathbf{W}_{y}^{h}\left( y\right)
e^{2\lambda y}\right\Vert _{L^{2,h}\left( \Omega ^{h}\right) }^{2}+\frac{1}{8%
}e^{-2\lambda b}\lambda ^{2}\left\Vert \mathbf{W}^{h}e^{2\lambda
y}\right\Vert _{L^{2,h}\left( \Omega ^{h}\right) }^{2}
\\
&\geq \frac{1}{4}e^{-2\lambda \left( b-a\right) }\left\Vert \mathbf{W}%
_{y}^{h}\left( y\right) \right\Vert _{L^{2,h}\left( \Omega ^{h}\right) }^{2}+%
\frac{1}{8}e^{-2\lambda \left( b-a\right) }\lambda ^{2}\left\Vert \mathbf{W}%
^{h}\right\Vert _{L^{2,h}\left( \Omega ^{h}\right) }^{2},
\end{align*}
for all $\mathbf{W}^{h}$ such that $\widetilde{\mathbf{W}}^{h}\in \widetilde{H}^{1,h}(
\Omega ^{h}),$ see Remark \ref{rem bdry}. 
Therefore, 
\begin{equation}
\left\Vert \mathbf{W}_{y}^{h}\left( y\right) +D^{h}\left( y\right) \mathbf{W}%
_{x}^{h}\left( y\right) \right\Vert _{L^{2,h}\left( \Omega ^{h}\right)
}^{2}\geq C_{h_{0}}\left\Vert \mathbf{W}^{h}\right\Vert _{H^{1,h}\left(
\Omega ^{h}\right) }^{2},
 \label{5.15}
\end{equation}
for all $\mathbf{W}^{h}$ such that $\widetilde{\mathbf{W}}^{h}\in \widetilde{H}^{1,h}(
\Omega ^{h}),$ see Remark \ref{rem bdry}, 
and $h\in \left[ h_{0},1\right] .$ It follows from (\ref{5.15}) that we can
define a new scalar product in the space $\widetilde{H}^{1,h}\left( \Omega
^{h}\right) $ as 
\begin{equation}
\left\{ \widetilde{\mathbf{P}}^{h},\widetilde{\mathbf{Q}}^{h}\right\} 
=
\left( \mathbf{P}_{y}^{h}\left( y\right) +D^{h}\left( y\right) \mathbf{P}%
_{x}^{h}\left( y\right) ,\mathbf{Q}_{y}^{h}\left( y\right) +D^{h}\left(
y\right) \mathbf{Q}_{x}^{h}\left( y\right) \right) ^{h}
+\varepsilon \left[ 
\mathbf{P}^{h},\mathbf{Q}^{h}\right] ^{h}\label{5.16}
\end{equation}
for all $\mathbf{P}^{h},\mathbf{Q}
^{h}$ such that $\widetilde{\mathbf{P}}^{h},\widetilde{\mathbf{Q}}^{h}\in \widetilde{H}
^{1,h}( \Omega ^{h}).$
By (\ref{5.15}) and (\ref{5.16}) the corresponding norm $\left\{ \widetilde{%
\mathbf{P}}^{h}\right\} =\sqrt{\left\{ \widetilde{\mathbf{P}}^{h},\widetilde{%
\mathbf{P}}^{h}\right\} }$ satisfies the following inequalities: 
\begin{equation}
B_{h_{0}}^{\left( 2\right) }\left\Vert \widetilde{\mathbf{P}}^{h}\right\Vert
_{H^{1,h}\left( \Omega ^{h}\right) }^{2}\geq \left\{ \widetilde{\mathbf{P}}%
^{h}\right\} ^{2}\geq B_{h_{0}}^{\left( 1\right) }\left\Vert \widetilde{%
\mathbf{P}}^{h}\right\Vert _{H^{1,h}\left( \Omega ^{h}\right) }^{2} \quad 
\widetilde{\mathbf{P}}^{h}\in \widetilde{H}^{1,h}\left( \Omega ^{h}\right) ,
\label{5.17}
\end{equation}%
for all $h\in \left[ h_{0},1\right] ,$ where constants $B_{h_{0}}^{\left(
1\right) },B_{h_{0}}^{\left( 2\right) }>0$ depend only on the number $h_{0}$
and the matrix $D^{h}\left( y\right) $ (recall that $\varepsilon \in \left[
0,1\right) ).$ Hence, the new norm $\left\{ \widetilde{\mathbf{P}}%
^{h}\right\} $ in $\widetilde{H}^{1,h}\left( \Omega ^{h}\right) $ is
equivalent with the previous norm $\left\Vert \widetilde{\mathbf{P}}%
^{h}\right\Vert _{H^{1,h}\left( \Omega ^{h}\right) }$ for $\widetilde{%
\mathbf{P}}^{h}\in \widetilde{H}^{1,h}\left( \Omega ^{h}\right) $. Hence,
using (\ref{5.14}), we obtain 
\begin{equation}
\left\{ \widetilde{\mathbf{V}}_{\min }^{h},\widetilde{\mathbf{W}}%
^{h}\right\} =  
-\left( \mathbf{G}^{h}\left( y\right) ,\mathbf{W}_{y}^{h}\left( y\right)
+D^{h}\left( y\right) \mathbf{W}_{x}^{h}\left( y\right) \right)
^{h}-\varepsilon \left[ \mathbf{F}^{h},\mathbf{W}^{h}\right] ^{h},
\label{5.170}
\end{equation}
for all $\mathbf{W}^{h}$ such that $\widetilde{\mathbf{W}}^{h}\in \widetilde{H}^{1,h}(
\Omega ^{h}),$ see Remark \ref{rem bdry}.
Next, (\ref{5.200}), (\ref{5.13}) and (\ref{5.16}) imply that the right hand
side of (\ref{5.170}) can be estimated from the above as%
\begin{multline}
\left\vert \left( \mathbf{G}^{h}\left( y\right) ,\mathbf{W}_{y}^{h}\left(
y\right) +D^{h}\left( y\right) \mathbf{W}_{x}^{h}\left( y\right) \right)
^{h}-\varepsilon \left[ \mathbf{F}^{h},\mathbf{W}^{h}\right] ^{h}\right\vert
\\
\leq C_{h_{0}}\left\Vert \mathbf{F}^{h}\right\Vert _{H^{1,h}\left( \Omega
^{h}\right) }\left\Vert \mathbf{W}^{h}\right\Vert _{H^{1,h}\left( \Omega
^{h}\right)},  \label{5.171}
\end{multline}%
for all $\mathbf{W}^{h}$ such that $\widetilde{\mathbf{W}}^{h}\in \widetilde{H}%
^{1,h}\left( \Omega ^{h}\right),$ see Remark \ref{rem bdry}. Hence, Riesz theorem and (\ref{5.17})
imply that there exists unique matrix $\widetilde{\mathbf{P}}^{h}\in \in 
\widetilde{H}^{1,h}\left( \Omega ^{h}\right) $ such that the right hand side
of (\ref{5.170}) can be represented as%
\begin{equation}
-\left( \mathbf{G}^{h}\left( y\right) ,\mathbf{W}_{y}^{h}\left( y\right)
+D^{h}\left( y\right) \mathbf{W}_{x}^{h}\left( y\right) \right)
^{h}-\varepsilon \left[ \mathbf{F}^{h},\mathbf{W}^{h}\right] ^{h}=\left\{ 
\widetilde{\mathbf{P}}^{h},\widetilde{\mathbf{W}}^{h}\right\},  \label{5.18}
\end{equation}%
for all $\mathbf{W}^{h}$ such that $\widetilde{\mathbf{W}}^{h}\in \widetilde{%
H}^{1,h}\left( \Omega ^{h}\right),$ see Remark \ref{rem bdry}. Furthermore, Riesz theorem and (\ref%
{5.171}) also imply that 
\begin{equation}
\left\Vert \widetilde{\mathbf{P}}^{h}\right\Vert _{\widetilde{H}^{1,h}\left(
\Omega ^{h}\right) }=\left\Vert \mathbf{P}^{h}\right\Vert _{H^{1,h}\left(
\Omega ^{h}\right) }\leq C_{h_{0}}\left\Vert \mathbf{F}^{h}\right\Vert
_{H^{1,h}\left( \Omega ^{h}\right) }.  \label{5.19}
\end{equation}%
Hence, (\ref{5.170}) and (\ref{5.18}) imply that%
\begin{equation*}
\left\{ \widetilde{\mathbf{V}}_{\min }^{h},\widetilde{\mathbf{W}}%
^{h}\right\} =\left\{ \widetilde{\mathbf{P}}^{h},\widetilde{\mathbf{W}}%
^{h}\right\} ,\quad \forall \widetilde{\mathbf{W}}^{h}\in \widetilde{H}%
^{1,h}\left( \Omega ^{h}\right) .
\end{equation*}%
This means that the minimizer $\widetilde{\mathbf{V}}_{\min }^{h}$ of the
functional $I_{\alpha }^{h}\left( \widetilde{\mathbf{V}}^{h}\right) $
exists, it is unique, and $\widetilde{\mathbf{V}}_{\min }^{h}=\widetilde{%
\mathbf{P}}^{h}\in \widetilde{H}^{1,h}\left( \Omega ^{h}\right) $ $.$
Therefore, the unique solution of Problem \ref{pro 2} is the matrix $%
\widetilde{\mathbf{U}}_{\min }^{h}=\widetilde{\mathbf{P}}^{h}+\widetilde{%
\mathbf{F}}^{h}.$

To prove (\ref{5.11}), we use the fourth line of (\ref{5.50}) and (\ref{5.19}%
) 
\begin{align}
\left\Vert \widetilde{\mathbf{U}}_{\min }^{h}-\widetilde{\mathbf{F}}%
^{h}\right\Vert _{\widetilde{H}^{1,h}\left( \Omega ^{h}\right) }&=\left\Vert 
\widetilde{\mathbf{P}}^{h}\right\Vert _{\widetilde{H}^{1,h}\left( \Omega
^{h}\right) }=\left\Vert \mathbf{P}^{h}\right\Vert _{H^{1,h}\left( \Omega
^{h}\right) } \notag
\\
&=\left\Vert \mathbf{U}_{\min }^{h}-\mathbf{F}^{h}\right\Vert _{H^{1,h}\left(
\Omega ^{h}\right) }\leq C_{h_{0}}\left\Vert \mathbf{F}^{h}\right\Vert
_{H^{1,h}\left( \Omega ^{h}\right) }.  \label{5.190}
\end{align}%
By the triangle inequality and (\ref{5.190}) 
\begin{equation*}
\left\Vert \mathbf{U}_{\min }^{h}\right\Vert _{H^{1,h}\left( \Omega
^{h}\right) }-\left\Vert \mathbf{F}^{h}\right\Vert _{H^{1,h}\left( \Omega
^{h}\right) }\leq \left\Vert \mathbf{U}_{\min }^{h}-\mathbf{F}%
^{h}\right\Vert _{H^{1,h}\left( \Omega ^{h}\right) }\leq C_{h_{0}}\left\Vert 
\mathbf{F}^{h}\right\Vert _{H^{1,h}\left( \Omega ^{h}\right)}.
\end{equation*}
Hence, $\left\Vert \mathbf{U}_{\min }^{h}\right\Vert _{H^{1,h}\left( \Omega
^{h}\right) }\leq \left( C_{h_{0}}+1\right) \left\Vert \mathbf{F}%
^{h}\right\Vert _{H^{1,h}\left( \Omega ^{h}\right) }.$ $\ \square $

\subsection{Convergence rate of regularized solutions}

\label{sec 5.3}

Let $\widetilde{\mathbf{U}}_{\min }^{h}$ be the minimizer of the functional $%
J_{\alpha }^{h}\left( \widetilde{\mathbf{U}}^{h}\right) ,$ which was found
in Theorem \ref{thm existence}. Then $\widetilde{\mathbf{U}}_{\min }^{h}$ is
called the \textquotedblleft regularized solution" in the regularization
theory \cite{T}. Naturally, it is important to prove convergence of
regularized solutions to the exact solution of the overdetermined system of
PDEs \eqref{5.5}, as long as the level of the noise in the data of second
and third lines of (\ref{5.5}) tends to zero. Recall that, according to the
regularization theory, one needs to assume the existence of the
\textquotedblleft idealized" exact solution, i.e. the solution which
corresponds to the noiseless data \cite{T}.

Let $\delta \in \left( 0,1\right) $ be the level of noise in the data. Let $%
\widetilde{\mathbf{U}}^{\ast ,h}\left( \mathbf{x}\right) $ be the exact
solution of problem (\ref{5.5}) with noiseless data $\mathbf{g}^{\ast
,h}=\left( \mathbf{g}^{\ast ,h}(a)=\mathbf{0},\mathbf{g}^{\ast ,h}(b),%
\mathbf{g}^{\ast }\left( -R,y\right) ,\mathbf{g}^{\ast }\left( R,y\right)
\right) $. Suppose that there exists a matrix 
\[
\widetilde{\mathbf{F}}^{\ast
,h}\left( y\right) =\left( \mathbf{F}_{0}^{\ast ,h}\left( y\right) ,\mathbf{F%
}_{1}^{\ast ,h}\left( y\right) ,\dots,\mathbf{F}_{K-1}^{\ast ,h}\left(
y\right) ,\mathbf{F}_{K}^{\ast ,h}\left( y\right) \right) ^{T}
\]
 such that $%
\mathbf{F}^{\ast ,h}\left( y\right) =\left( \mathbf{F}_{1}^{\ast ,h}\left(
y\right) ,\dots,\mathbf{F}_{K-1}^{\ast ,h}\left( y\right) \right) ^{T}\in
H^{1,h}\left( \Omega ^{h}\right) $
\begin{equation*}
\mathbf{F}^{\ast ,h}\left( a\right) =\mathbf{0},\mathbf{F}^{\ast ,h}\left(
b\right) =\mathbf{g}^{\ast ,h}(b);\text{ }\mathbf{F}_{0}^{\ast ,h}\left(
y\right) =\mathbf{g}^{\ast }\left( -R,y\right) ,\mathbf{F}_{K}^{\ast
,h}\left( y\right) =\mathbf{g}^{\ast }\left( R,y\right),
\end{equation*}%
 for all $y \in (a, b)$.
Also, let $\mathbf{g}_{\delta }^{h}=\left( \mathbf{g}_{\delta
}^{h}\left( a\right) =\mathbf{0},\mathbf{g}_{\delta }^{h}(b),\mathbf{g}%
_{\delta }\left( -R,y\right) ,\mathbf{g}_{\delta }\left( R,y\right) \right) $
be the noisy data in (\ref{5.5}) and assume that there exists a matrix $%
\widetilde{\mathbf{F}}_{\delta }^{h}\left( y\right) $\emph{\ }satisfying the
same conditions as ones for $\widetilde{\mathbf{F}}^{\ast ,h}\left( y\right) 
$ with the replacement of $\mathbf{g}^{\ast ,h}$ by $\mathbf{g}_{\delta }^{h}
$.\emph{\ }We assume that the following error estimate holds: 
\begin{equation}
\left\Vert \mathbf{F}_{\delta }^{h}-\mathbf{F}^{\ast ,h}\right\Vert
_{H^{1,h}\left( \Omega ^{h}\right) }\leq \delta .  \label{5.20}
\end{equation}

\begin{theorem}[The convergence of the regularized solution to the exact one]
Let $\widetilde{\mathbf{U}}^{\ast ,h}\left( \mathbf{x}\right) $ be the exact
solution of the problem (\ref{5.5}) with noiseless data $\mathbf{g}^{\ast
,h},$ which replace $\mathbf{g}^{h}$ in \eqref{5.5}. Let $\widetilde{\mathbf{%
U}}_{\delta }^{h}\left( \mathbf{x}\right) $ be the solution of (\ref{5.5})
with noisy data $\mathbf{g}_{\delta }^{h}$, which was found in Theorem \ref%
{thm existence}. Assume that conditions of Theorem \ref{thm existence} hold
true and that the error estimate (\ref{5.20}) is valid. Then for all $h\in %
\left[ h_{0},1\right) $ the following convergence rate is valid: 
\begin{equation}
\left\Vert \mathbf{U}_{\delta }^{h}-\mathbf{U}^{\ast ,h}\right\Vert
_{H^{1,h}\left( \Omega ^{h}\right) }\leq C_{h_{0}}\left( \delta +\sqrt{%
\varepsilon }\left\Vert \mathbf{U}^{\ast ,h}\right\Vert _{H^{1,h}\left(
\Omega ^{h}\right) }\right) .  \label{5.21}
\end{equation}%
In particular, choosing $\varepsilon \in \left[ 0,\delta ^{2}\right] ,$ we
obtain $\left\Vert \mathbf{U}_{\delta }^{h}-\mathbf{U}^{\ast ,h}\right\Vert
_{H^{1,h}\left( \Omega ^{h}\right) }\leq C_{h_{0}}\delta .$ \label{thm
convergence u}
\end{theorem}

\textbf{Proof}. Denote $\widetilde{\mathbf{V}}_{\delta }^{h}=\widetilde{%
\mathbf{U}}_{\delta }^{h}-\widetilde{\mathbf{F}}_{\delta }^{h}\in $ and $%
\widetilde{\mathbf{V}}^{\ast ,h}=\widetilde{\mathbf{U}}^{\ast ,h}-\widetilde{%
\mathbf{F}}^{\ast ,h}$. Similarly with (\ref{5.13}) and (\ref{5.14}), we
obtain 
\begin{multline}
\left( \mathbf{V}_{\delta y}^{h}\left( y\right) +D^{h}\left( y\right) 
\mathbf{V}_{\delta x}^{h}\left( y\right) ,\mathbf{W}_{y}^{h}\left( y\right)
+D^{h}\left( y\right) \mathbf{W}_{x}^{h}\left( y\right) \right)
^{h}+\varepsilon \left[ \mathbf{V}_{\delta }^{h},\mathbf{W}^{h}\right] ^{h}
\\
=-\left( \mathbf{G}_{\delta }^{h}\left( y\right) ,\mathbf{W}_{y}^{h}\left(
y\right) +D^{h}\left( y\right) \mathbf{W}_{x}^{h}\left( y\right) \right)
^{h}-\varepsilon \left[ \mathbf{F}_{\delta }^{h},\mathbf{W}^{h}\right] ^{h},
\label{5.22}
\end{multline}
for all $\mathbf{W}^{h}$ such that $\widetilde{\mathbf{W}}^{h}\in \widetilde{%
H}^{1,h}\left( \Omega ^{h}\right)$, see Remark \ref{rem bdry}, where
\begin{equation}
\mathbf{G}_{\delta }^{h}\left( y\right) =\mathbf{F}_{\delta y}^{h}\left(
y\right) +D^{h}\left( y\right) \mathbf{F}_{\delta x}^{h}\left( y\right)
,\quad y\in \left( a,b\right) .  \label{5.23}
\end{equation}%
 Also, by the same arguments, we have 
\begin{multline}
\left( \mathbf{V}_{y}^{\ast ,h}\left( y\right) +D^{h}\left( y\right) \mathbf{%
V}_{x}^{\ast ,h}\left( y\right) ,\mathbf{W}_{y}^{h}\left( y\right)
+D^{h}\left( y\right) \mathbf{W}_{x}^{h}\left( y\right) \right)
^{h}+\varepsilon \left[ \mathbf{V}^{\ast h},\mathbf{W}^{h}\right] ^{h}
\\
=-\left( \mathbf{G}^{\ast ,h}\left( y\right) ,\mathbf{W}_{y}^{h}\left(
y\right) +D^{h}\left( y\right) \mathbf{W}_{x}^{h}\left( y\right) \right)
^{h}+\varepsilon \left[ \mathbf{V}^{\ast ,h},\mathbf{W}^{h}\right] ^{h},
\label{5.25}
\end{multline}
for all $\mathbf{W}^{h}$ such that $\widetilde{\mathbf{W}}^{h}\in \widetilde{%
H}^{1,h}\left( \Omega ^{h}\right)$, see Remark \ref{rem bdry}, where
\begin{equation}
\mathbf{G}^{\ast ,h}\left( y\right) =\mathbf{F}_{y}^{\ast ,h}\left( y\right)
+D^{h}\left( y\right) \mathbf{F}_{x}^{\ast ,h}\left( y\right) ,\quad y\in
\left( a,b\right) .  \label{5.26}
\end{equation}%
Denote $\widetilde{\mathbf{X}}_{\delta }^{h}=\widetilde{\mathbf{V}}_{\delta
}^{h}-\widetilde{\mathbf{V}}^{\ast ,h}$ and 
\begin{equation}
\mathbf{Z}_{\delta }^{h}=\mathbf{G}_{\delta }^{h}-\mathbf{G}^{\ast ,h}.
\label{5.27}
\end{equation}%
Obviously, $\widetilde{\mathbf{X}}_{\delta }^{h}\in \widetilde{H}%
^{1,h}\left( \Omega ^{h}\right) .$ Subtracting (\ref{5.25}) from (\ref{5.22}%
) and using (\ref{5.23}) and (\ref{5.26}), we obtain%
\begin{multline}
\left( \mathbf{X}_{\delta y}^{h}\left( y\right) +D^{h}\left( y\right) 
\widetilde{\mathbf{X}}_{\delta x}^{h}\left( y\right) ,\mathbf{W}%
_{y}^{h}\left( y\right) +D^{h}\left( y\right) \mathbf{W}_{x}^{h}\left(
y\right) \right) ^{h}+\varepsilon \left[ \mathbf{X}_{\delta }^{h},\mathbf{W}%
^{h}\right] ^{h}
\\
=-\left( \mathbf{Z}_{\delta }^{h}\left( y\right) ,\mathbf{W}_{y}^{h}\left(
y\right) +D^{h}\left( y\right) \mathbf{W}_{x}^{h}\left( y\right) \right)
^{h}-\varepsilon \left[ \mathbf{V}^{\ast ,h}+\mathbf{F}_{\delta }^{h},%
\mathbf{W}^{h}\right]^{h},
\label{5.28}
\end{multline}
for all $\mathbf{W}^{h}$ such that $\widetilde{\mathbf{W}}^{h}\in \widetilde{%
H}^{1,h}\left( \Omega ^{h}\right),$ see Remark \ref{rem bdry}.
By (\ref{5.27}) 
\begin{equation}
\Vert \mathbf{Z}_{\delta }^{h}\Vert _{L^{2,h}\left( \Omega ^{h}\right) }\leq
C_{h_{0}}\delta .  \label{5.29}
\end{equation}%
Setting in (\ref{5.28}) $\widetilde{\mathbf{W}}^{h}=$ $\widetilde{\mathbf{X}}%
_{\delta }^{h}$, noting that 
\begin{equation*}
\left\Vert \mathbf{V}^{\ast ,h}+\mathbf{F}_{\delta }^{h}\right\Vert
_{H^{1,h}\left( \Omega ^{h}\right) }=\left\Vert \mathbf{U}^{\ast ,h}+\left( 
\mathbf{F}_{\delta }^{h}-\mathbf{F}^{\ast ,h}\right) \right\Vert
_{H^{1,h}\left( \Omega ^{h}\right) }\leq \Vert \mathbf{U}^{\ast ,h}\Vert
_{H^{1,h}\left( \Omega ^{h}\right) }+\delta ,
\end{equation*}%
and using (\ref{5.29}), we obtain%
\begin{equation}
\Vert \mathbf{X}_{\delta y}^{h}\left( y\right) +D^{h}\left( y\right) \mathbf{%
X}_{\delta x}^{h}\left( y\right) \Vert _{L^{2,h}\left( \Omega ^{h}\right)
}^{2}+\varepsilon \Vert \mathbf{X}_{\delta }^{h}\Vert _{H^{1,h}\left( \Omega
^{h}\right) }^{2}\leq C_{h_{0}}\delta ^{2}+\varepsilon \Vert \mathbf{U}%
^{\ast ,h}\Vert _{H^{1,h}\left( \Omega ^{h}\right) }^{2}.  \label{5.30}
\end{equation}%
Ignoring in (\ref{5.30}) the positive term $\varepsilon \Vert \mathbf{X}%
_{\delta }^{h}\Vert _{H^{1,h}\left( \Omega ^{h}\right) }^{2},$ recalling
that $\widetilde{\mathbf{X}}_{\delta }^{h}\in \widetilde{H}^{1,h}\left(
\Omega ^{h}\right) $ and applying (\ref{5.15}) to the rest of the left hand
side of (\ref{5.30}) we obtain%
\begin{equation}
\left\Vert \mathbf{X}_{\delta }^{h}\right\Vert _{H^{1,h}\left( \Omega
^{h}\right) }\leq C_{h_{0}}\left( \delta +\sqrt{\varepsilon }\left\Vert 
\mathbf{U}^{\ast ,h}\right\Vert _{H^{1,h}\left( \Omega ^{h}\right) }\right) .
\label{5.31}
\end{equation}%
Since 
\begin{equation*}
\mathbf{X}_{\delta }^{h}=\mathbf{V}_{\delta }^{h}-\mathbf{V}^{\ast
,h}=\left( \mathbf{U}_{\delta }^{h}-\mathbf{F}_{\delta }^{h}\right) -\left( 
\mathbf{U}^{\ast ,h}-\mathbf{F}^{\ast ,h}\right) =\left( \mathbf{U}_{\delta
}^{h}-\mathbf{U}^{\ast ,h}\right) -\left( \mathbf{F}_{\delta }^{h}-\mathbf{F}%
^{\ast ,h}\right) ,
\end{equation*}%
then by the triangle inequality and (\ref{5.20}) 
\begin{align}
\left\Vert \mathbf{X}_{\delta }^{h}\right\Vert _{H^{1,h}\left( \Omega
^{h}\right) } &\geq \left\Vert \mathbf{U}_{\delta }^{h}-\mathbf{U}^{\ast
,h}\right\Vert _{H^{1,h}\left( \Omega ^{h}\right) }-\left\Vert \mathbf{F}%
_{\delta }^{h}-\mathbf{F}^{\ast ,h}\right\Vert _{H^{1,h}\left( \Omega
^{h}\right) } \notag
\\
&=\left\Vert \mathbf{U}_{\delta }^{h}-\mathbf{U}^{\ast
,h}\right\Vert _{H^{1,h}\left( \Omega ^{h}\right) }-\delta .  \label{5.310}
\end{align}%
Thus, (\ref{5.21}) follows from (\ref{5.31}) and (\ref{5.310}). $\square$

\subsection{Reconstruction formula and its accuracy}

\label{sec 5.4}

We now estimate the accuracy of the reconstruction of the target function $%
f\left( \mathbf{x}\right) $. Recalling (\ref{2}), let $\mathbf{f}^{h}\left(
y\right) =\left( f\left( -R+h,y\right) ,\dots,f\left( -R+K-1\right) h,y\right)
,y\in \left( a,b\right) $ be the discrete analog of the function $f\left( 
\mathbf{x}\right) .$ By (\ref{main eqn}), (\ref{4.2}), (\ref{4.3}) and (\ref%
{5.3}) we have the following reconstruction formula for $y\in \left(
a,b\right) $ and $j=1,\dots,K-1:$  
\begin{multline}
	f(-R + jh) = \frac{1}{2d} \sum_{n = 1}^N\int_{-d}^d \Big[
		\frac{\left(
-R+jh\right) -\alpha }{\sqrt{\left( \left( -R+jh\right) -\alpha \right)
^{2}+y^{2}}}u_{n,jx}(-R+jh,y) \\
	+ \frac{y}{\sqrt{\left( \left( -R+jh\right) -\alpha \right) ^{2}+y^{2}}}%
\sum_{n=0}^{N-1}u_{n,y}(-R+jh,y)
	\Big] \Psi_n(\alpha)d\alpha
	\label{5.32}
\end{multline}

We have taken the average value with respect to $\alpha $ since the
integrand in (\ref{5.32}) depends on $\alpha $ in practical computations
whereas the function $f\left( ih,y\right) $ does not depend on $\alpha .$
That dependence on $\alpha $ is due to the approximate nature of our method.

To obtain the desired accuracy estimate, we note that in the case of the
noisy data discussed in Section \ref{sec 5.3} functions $u_{n}(-R+jh,y)$ in (%
\ref{5.32}) should be replaced with the components $u_{\delta
,n}^{h}(-R+jh,y)$ of the matrix $\mathbf{U}_{\delta }^{h},$ and in the case
of noiseless data they should be replaced with the components $u_{n}^{\ast
,h}(-R + jh,y)$ of the matrix $\mathbf{U}^{\ast ,h}$. Let $\mathbf{f}_{\delta
}^{h}\left( \mathbf{x}\right) $ and $\mathbf{f}^{\ast ,h}\left( \mathbf{x}%
\right) $ be the right hand sides of corresponding analogs of formula (\ref%
{5.32}). Subtracting these analogs and using Theorem \ref{thm convergence u}%
, we easily prove the following theorem:

\begin{theorem}[The convergence of the computed tomographic function to the
true one in our approximate context]
Assume that all conditions of Theorem \ref{thm convergence u} hold true.
Then for all $h\in \left[ h_{0},1\right) $ the following analog of the
convergence rate (\ref{5.21}) is valid: 
\begin{equation*}
\Vert \mathbf{f}_{\delta }^{h}-\mathbf{f}^{\ast ,h}\Vert _{L^{2,h}\left(
\Omega ^{h}\right) }\leq C_{h_{0}}\left( \delta +\sqrt{\varepsilon }\Vert 
\mathbf{f}^{\ast ,h}\Vert _{L^{2,h}\left( \Omega ^{h}\right) }\right) .
\end{equation*}%
In particular, choosing $\varepsilon \in \left[ 0,\delta ^{2}\right], $we
obtain $\Vert \mathbf{f}_{\delta }^{h}-\mathbf{f}^{\ast ,h}\Vert
_{L^{2,h}\left( \Omega ^{h}\right) }\leq C_{h_{0}}\delta .$ \label{thm
convergence f}
\end{theorem}

%
%
%

\section{Numerical implementation}

\label{sec imple}

In this section, we present some details of our computational implementation
for the numerical solution of Problem \ref{pro 2.1}. Recall that the $\Omega 
$ is defined in (\ref{1}), where numbers $R$, $a$ and $b$ will be chosen
later in each test of Section \ref{sec num}. In all our tests the line
segment with the sources $\Gamma _{d}$ in (\ref{2.10}) is the same, 
\begin{equation}
\Gamma _{d}=\left\{ \mathbf{x}=\left( x,y\right) :x\in
(-3.5,3.5),y=0\right\} =\left\{ \left( \alpha ,0\right) :\alpha \in \left(
-d,d\right) \right\} .  \label{6.1}
\end{equation}%
We calculate derivatives using finite differences. To do so, we fix the
number $T_{\mathbf{x}}=150$ and then consider grid points in the rectangle $%
\Omega ,$ 
\begin{equation}
(x_{i},y_{j})=(-R+(i-1)h _{x},a+(j-1)h _{y}),\quad 1\leq i,j\leq
T_{\mathbf{x}}+1.
\label{grid points}
\end{equation}%
where $h _{x}=2R/T_{\mathbf{x}}$ and $h _{y}=\left( b-a\right) /T_{%
\mathbf{x}}$ are the grid step sizes in $x$ and $y$ directions respectively.

By (\ref{6.1}) $d=3.5$ and the length of the line with sources is $7$. We
uniformly split the source interval $(-d,d)$ into $T_{\alpha }=100$
subintervals whose edge points are 
\begin{equation}
\alpha _{i}=-d+(i-1)\frac{2d}{T_{\alpha }},\quad i=1,\dots ,T_{\alpha }+1.
\label{6.2}
\end{equation}

\subsection{The forward problem and the noisy data}

\label{sec 6.1}

We solve the forward problem by calculating the Riemannian sum in the
integral in (\ref{2.2}) 
\begin{equation*}
u(\mathbf{x},\mathbf{x}_{\alpha })=\int_{L(\mathbf{x},\mathbf{x}%
_{\alpha })}f(\mathbf{\xi })d\sigma 
\end{equation*}%
for each $\mathbf{x}\in \partial \Omega ,\mathbf{x}_{\alpha _{i}}=(\alpha
_{i},0)\in \Gamma _{d}$. The step size of this sum depends on the pair $%
\mathbf{x},\mathbf{x}_{\alpha _{i}}$ and is chosen in such a way that there
are 150 grid points along the part $\widetilde{L}(\mathbf{x},\mathbf{x}%
_{\alpha _{i}})$ of the line $L(\mathbf{x},\mathbf{x}_{\alpha _{i}})$ which
lies inside of $\Omega :$ since $f\left( \mathbf{x}\right) =0$ outside of $%
\Omega .$ 

We generate random noise in our data for Problem \ref{pro 2.1} as 
\begin{equation}
Rf(\mathbf{x},\mathbf{x}_{\alpha _{i}})=u(\mathbf{x},\mathbf{x}_{\alpha
_{i}})(1+\delta (2\mathrm{rand}(\mathbf{x})-1)),\quad \mathbf{x}\in \partial
_{\alpha _{i}}\Omega ,i=1,\dots ,T_{\alpha },  \label{6.30}
\end{equation}%
where $\partial _{\alpha _{i}}\Omega =\partial \Omega \cap L(\mathbf{x},%
\mathbf{x}_{\alpha _{i}}),\delta >0$ is the noise level and $\mathrm{rand}$
is the function that generates uniformly distributed random numbers in the
interval $[0,1]$. In this paper, we choose two noise levels $\delta =0.05$
and $\delta =0.15,$ which correspond to $5\%$ and $15\%$ noise respectively.

The boundary data (\ref{4.10}) are read as%
\begin{equation}
\mathbf{g}(\mathbf{x})=\left( g_{1},\dots,g_{N}\right) ^{T}\left( \mathbf{x}%
\right) ,\quad g_{n}\left( \mathbf{x}\right) =\int_{-d}^{d}Rf(\mathbf{x},%
\mathbf{x}_{\alpha })\Psi _{n}(\alpha )d\alpha .  \label{eqn boundary data}
\end{equation}

\subsection{Calculating the vector function $\mathbf{U}$ and computing the
target function $f_{\mathrm{comp}}$}

\label{sec our approach}

Equation (\ref{4.9}) is obtained from equation (\ref{4.4}) and then from (%
\ref{4.7}) via singling out the $y-$derivative $\mathbf{U}_{y}.$ The latter
is done using the inverse of the matrix $M_{N}\left( I-M_{N}^{-1}D_{1}\left( 
\mathbf{x}\right) \right) $ (Lemma \ref{lem 3.3}). While equation (\ref{4.9}%
) is convenient for the theoretical analysis of Section \ref{sec 5}, our
computational experience tells us that in computations better not to invert
the matrix $M_{N}\left( I-M_{N}^{-1}D_{1}\left( \mathbf{x}\right) \right) $.
Thus, we work with an equivalent equation, in which the $y-$derivative $%
\mathbf{U}_{y}$ is not singled out. Denote 
\begin{equation}
A(\mathbf{x})=-D_{2}\left( \mathbf{x}\right) ,\quad B(\mathbf{x}%
)=M_{N}-D_{1}\left( \mathbf{x}\right) .  \label{eqn AB}
\end{equation}%
So, this equation together with the boundary condition (\ref{4.10}) becomes 
\begin{equation}
\left\{ 
\begin{array}{rcll}
A(\mathbf{x})\partial _{x}\mathbf{U}(\mathbf{x})+B(\mathbf{x})\partial _{y}%
\mathbf{U}(\mathbf{x}) & = & 0 & \mathbf{x}\in \Omega , \\ 
\mathbf{U}(\mathbf{x}) & = & \mathbf{g}(\mathbf{x}) & \mathbf{x}\in \partial
\Omega .%
\end{array}%
\right.   \label{6.3}
\end{equation}%
It follows from (\ref{4.8}) and Lemma \ref{lem 3.3} that for $a\geq
a_{0}\left( N,R,d\right) $ equation (\ref{6.3}) is equivalent with equation (%
\ref{4.9}).

We modify the objective functional (\ref{5.1}) as 
\begin{equation}
\mathcal{J}_{\epsilon _{1},\epsilon _{2}}(\mathbf{U})=\int_{\Omega }%
\left[ A(x,y)\partial _{x}\mathbf{U}(x,y)+B(x,y)\partial _{y}\mathbf{U}(x,y)%
\right] ^{2}dxdy  \label{new J}
\end{equation}%
\begin{equation*}
+\varepsilon _{1}\Vert \mathbf{U}\Vert _{L^{2}(\Omega )}^{2}+\varepsilon
_{2}\Vert \nabla \mathbf{U}\Vert _{L^{2}(\Omega )}^{2}.
\end{equation*}

\begin{remark}
In the original definition of this functional in (\ref{5.1}), we use only
one regularization parameter $\epsilon $. However, our computational
experience tells us that using two differential regularization parameters $%
\epsilon _{1}$ and $\epsilon _{2}$ yields better reconstructed results. In
this paper, we take $\epsilon _{1}=0.1$ and $\epsilon _{2}=0.01.$
\end{remark}

We consider the finite difference version of the functional $\mathcal{J}%
_{\epsilon _{1},\epsilon _{2}}(\mathbf{U}^{h}),$ 
\begin{multline*}
\mathcal{J}_{\epsilon _{1},\epsilon _{2}}^{h}(\mathbf{U}^{h})=h
_{x}h _{y}\sum_{i,j=2}^{T_{\mathbf{x}}}\Big|A(x_{i},y_{j})\frac{\mathbf{%
U}^{h}(x_{i+1},y_{j})-\mathbf{U}^{h}(x_{i},y_{j})}{h _{x}} \\
+B(x_{i},y_{j})\frac{\mathbf{U}^{h}(x_{i},y_{j+1})-\mathbf{U}%
^{h}(x_{i},y_{j})}{h _{y}}\Big|^{2}+\epsilon _{1}h _{x}h
_{y}\sum_{i,j=1}^{T_{\mathbf{x}}+1}|\mathbf{U}^{h}(x_{i},y_{j})|^{2} \\
+\epsilon _{2}h _{x}h _{y}\sum_{i,j=1}^{T_{\mathbf{x}}}\Big(\frac{|%
\mathbf{U}^{h}(x_{i+1},y_{j})-\mathbf{U}^{h}(x_{i},y_{j})|^{2}}{h _{x}}+%
\frac{|\mathbf{U}^{h}(x_{i},y_{j+1})-\mathbf{U}^{h}(x_{i},y_{j})|^{2}}{%
h _{y}}\Big),
\end{multline*}%
in which the integral in \eqref{new J} is approximated by its Riemann sum
and the derivative of $\mathbf{U}^{h}$ is in the finite difference form. Recall that the grid points $(x_i, y_j)$, $1 \leq i, j \leq T_\x + 1$, are defined in \eqref{grid points}. 
In this Rieman sum, by letting the indices $i$ and $j$ start from $2$ rather than 1, we ignore the boundary value of the integrand on $\partial \Omega$. This is acceptable since the measure of $\partial \Omega$ is zero.
As to the boundary conditions, see \eqref{5.12}  and \eqref{matrix C}
The functional $\mathcal{J}_{\epsilon _{1},\epsilon _{2}}^{h}(\mathbf{U}^{h})$
is written, with some suitable arrangement, in terms of entries of $\mathbf{U%
}^{h}$ as follows 
\begin{multline}
\mathcal{J}_{\epsilon _{1},\epsilon _{2}}^{h}(\mathbf{U}^{h})=h
_{x}h _{y}\sum_{i,j=2}^{T_{\mathbf{x}}}\sum_{n=1}^{N}\Big|\frac{%
A(x_{i},y_{j})}{h _{x}}u_{n}(x_{i+1},y_{j})+\frac{B(x_{i},y_{j})}{%
h _{y}}u_{n}(x_{i},y_{j-1})  \notag \\
-\Big(\frac{A(x_{i},y_{j})}{h _{x}}+\frac{B(x_{i},y_{j})}{h _{y}}%
\Big)u_{n}(x_{i},y_{j})\Big|^{2}+\epsilon _{1}h _{x}h
_{y}\sum_{i,j=1}^{T_{\mathbf{x}}+1}\sum_{n=1}^{N}|u_{n}(x_{i},y_{j})|^{2}
\label{5.7} \\
+\epsilon _{2}h _{x}h _{y}\sum_{i,j=1}^{T_{\mathbf{x}%
}}\sum_{n=1}^{N}\Big(\frac{|u_{n}(x_{i+1},y_{j})-u_{n}(x_{i},y_{j})|^{2}}{%
h _{x}}+\frac{|u_{n}(x_{i},y_{j+1})-u_{n}(x_{i},y_{j})|^{2}}{h _{y}%
}\Big).  \notag
\end{multline}

We next identify 
\begin{equation*}
\{\mathbf{U}^{h}(x_{i},y_{j})=(u_{1}(x_{i},y_{j}),u_{2}(x_{i},y_{j}),\dots
,u_{n}(x_{i},y_{j})):1\leq i,j,\leq T_{\mathbf{x}}+1\}
\end{equation*}%
by a column vector 
\begin{equation}
\mathfrak{U}=(\mathfrak{u}_{1},\mathfrak{u}_{2},\dots ,\mathfrak{u}_{(T_{%
\mathbf{x}}+1)^{2}N})^{T}  \label{5.8}
\end{equation}%
where 
\begin{equation}
\mathfrak{u}_{\mathfrak{i}}=u_{n}(x_{i},y_{j})  \label{5.9}
\end{equation}%
with 
\begin{equation}
\mathfrak{i}=(i-1)(T_{\mathbf{x}}+1)N+(j-1)N+n.  \label{eqn index}
\end{equation}

\begin{remark}
The map 
\[\{1, \dots, T_\mathbf{x}+1\} \times \{1, \dots, T_\mathbf{x}+1\} \times \{1, \dots, N\} \to \{1,
\dots, (T_\mathbf{x} + 1)^2N\} 
\] that sends $(i, j, n)$ to $\mathfrak{i}$ as
in \eqref{eqn index} is onto and one to one.
\end{remark}

Define the $(T_\mathbf{x} + 1)^2N \times (T_\mathbf{x} + 1)^2N$ matrix $%
\mathcal{M} = (\mathfrak{m}_{\mathfrak{i} \mathfrak{j}})_{1 \leq \mathfrak{i}%
, \mathfrak{j} \leq (T_\mathbf{x} + 1)^2N} $ as follows. For any $\mathfrak{i%
} = (i-1)(T_{\mathbf{x}} + 1)N + (j - 1)N + m$, $i, j \in \{2, \dots, T_%
\mathbf{x}\}$ and $m \in \{1, \dots, N\}$, set

\begin{enumerate}
\item $\mathfrak{m}_{\mathfrak{i} \mathfrak{j}} = - \displaystyle\Big(\frac{%
A(x_i, y_j)}{h_x} + \frac{B(x_i, y_{j})}{h_y}\Big)$ if $\mathfrak{j%
}$ is identical with $(i, j, n)$ in the sense of \eqref{eqn index} for any $%
n \in \{1, \dots, N\}.$

\item $\displaystyle \mathfrak{m}_{\mathfrak{i} \mathfrak{j}} = \frac{A(x_i,
y_j)}{h_x} $ if $\mathfrak{j}$ is identical with $(i + 1, j, n)$ in the
sense of \eqref{eqn index} for any $n \in \{1, \dots, N\}.$

\item $\displaystyle \mathfrak{m}_{\mathfrak{i} \mathfrak{j}} = \frac{B(x_i,
y_j)}{h_y} $ if $\mathfrak{j}$ is identical with $(i, j + 1, n)$ in the
sense of \eqref{eqn index} for any $n \in \{1, \dots, N\}.$

\item $\displaystyle \mathfrak{m}_{\mathfrak{i} \mathfrak{j}} = 0$ for other
pair $\mathfrak{i}, \mathfrak{j}$.
\end{enumerate}

Using the matrix $\mathcal{M},$ we can shorten the function $\mathcal{J}%
_{\epsilon _{1},\epsilon _{2}}^{h}(\mathbf{U}^{h})$ in \eqref{5.7} as 
\begin{equation*}
\mathfrak{J}_{\epsilon _{1},\epsilon _{2}}(\mathfrak{U})=h _{x}h
_{y}\Big(|\mathcal{M}\mathfrak{U}|^{2}+\epsilon _{1}|\mathfrak{U}%
|^{2}+\epsilon _{2}|D_{x}\mathfrak{U}|^{2}+\epsilon _{2}|D_{y}\mathfrak{U}%
|^{2}\Big),
\end{equation*}%
where $D_{x}$ and $D_{y}$ are the matrix that provide the finite difference
approximations of the partial derivatives of $\mathfrak{U}$ with respect to $%
x$ and $y$. Computationally, we solve the following minimization problem:

\begin{problem*}[Solving \eqref{6.3} by the quasi-reversibility method]
Minimize the functional $\mathfrak{J}_{\epsilon_1, \epsilon_2}(\mathfrak{U})$%
, subject to the finite difference analog of the boundary condition %
\eqref{eqn boundary data} 
\begin{equation}
\mathfrak{u}_{(i - 1)(T_\mathbf{x} + 1)N + (j - 1)N + n} = g_n(x_i, y_j)
\label{5.12}
\end{equation}
for all $i, j$ such that $(x_i, y_j)$ is on $\partial \Omega$ and $n \in
\{1, \dots, N\}.$ \label{pro 3}
\end{problem*}

It is convenient to employ the Linear Algebra package of MATLAB to solve
Problem \ref{pro 3}. Denote 
\begin{equation}
\mathcal{C=M}^{T}\mathcal{M}+\epsilon _{1}\text{Id}+\epsilon _{2}{\mathcal{D}%
}_{x}^{T}{\mathcal{D}}_{x}+\epsilon _{2}{\mathcal{D}}_{y}^{T}{\mathcal{D}}%
_{y}.
\label{matrix C}
\end{equation}
It is obvious that the minimizer of $\mathfrak{J}_{\epsilon_1, \epsilon_2}$
satisfies the equation $\mathcal{C} \mathfrak{U}=0$ subject to the
constraint \eqref{5.12}. We use the command lsqlin of MATLAB to compute such
vector $\mathfrak{U}$. The knowledge of $\mathfrak{U}$ yields that of $%
\mathbf{U}^h$ via \eqref{5.8}--\eqref{eqn index}. Denote the result obtained
by the procedure of this section as $\mathbf{U}^h_{\mathrm{comp}}=(u_{1}^{%
\mathrm{comp}},\dots ,u_{N}^{\mathrm{comp}})^{T}.$ The knowledge of this
vector function directly provides the knowledge of the function $u_h^{%
\mathrm{comp}}(\mathbf{x},\mathbf{x}_{\alpha })$ via (\ref{4.2}). The
reconstructed function $f_{\mathrm{comp}}$ is determined using the
reconstruction formula (\ref{5.32}) in which functions $u_{n}^{h}$ are
replaced with $u_{n}^{\text{comp}}.$

\subsection{The reconstruction via the filtered back projection algorithm}

\label{sec Radon}

We wish to compare our computational results with the results of the Radon
transform inversion which is widely used in the scientific community. To do
this, we employ the built-in function \textquotedblleft iradon" in MATLAB to
reconstruct the function $f$ from our data. In MATLAB, one can use the
function \textquotedblleft radon" to compute the Radon transform $Rf$ of a
function $f$ and then use the command \textquotedblleft iradon" for the
function $Rf$ to compute $f$. The command \textquotedblleft iradon" is based
on the \textit{filtered back projection formula} which is very well-known in
the scientific community \cite{Natterer:cmsiam2001}. In the case of complete
data, the filtered back projection formula provides a perfect reconstruction
of the function $f$. However, in the case of our incomplete data, the
filtered back projection formula does not work. Thus, we simply set that the
data to be zero for all those angles which are not involved in the data (\ref%
{2.3}), see Figures 2c-5c. It is clear from a visual comparison of these
figures with figures 2-4 of \cite{Borg} that we complement the missing data
similarly with \cite{Borg}. Of course, assigning zero to missing data is not
rigorous. But we are doing so just to have a crude comparison of our method
with the filtered back projection method. We point out that a detailed study
of the comparison issue of the filtered back projection method with our
method is outside of the scope of this publication. All what we want here is
to compare our reconstructions with a version of the filtered back
projection algorithm in which the missing data are set to zero.

It is well known that the arguments of the Radon transform $Rf(r,\theta )$
for the filtered back projection formula are a \textquotedblleft signed"
radius $r\in \lbrack -l/2,l/2]$ where $l$ is the length of the diagonal of $%
\Omega $ and an angle $\theta \in (0^{\circ },180^{\circ })$. For each $%
\theta \in (0^{\circ },180^{\circ }),$ let the $r$-axis be the line passing
through the center of $\Omega $ with its positive direction is the direction
of $(\cos \theta ,\sin \theta ).$ Then, the function $Rf(r,\theta )$ is
given by%
\begin{equation*}
Rf(r,\theta )=\int_{L(\mathbf{x},\mathbf{x}_{\alpha })}fd\sigma .
\end{equation*}%
Here points $\mathbf{x}\in \partial \Omega $ and $\mathbf{x}_{\alpha }\in
\Gamma _{d}$ are such that the line $L(\mathbf{x},\mathbf{x}_{\alpha })$ is
perpendicular to the $r-$axis and the intersection of $L(\mathbf{x},\mathbf{x%
}_{\alpha })$ with the $r-$axis is at the point $r$ on the $r$-axis. For
each pair $(r,\theta )$ we have $(r,\theta )\in (-l/2,l/2)\times (0^{\circ
},180^{\circ })$. Thus we have two cases:

\begin{enumerate}
\item \textbf{Case 1}. If there exists a corresponding pair $(\mathbf{x},%
\mathbf{x}_{\alpha })\in \partial \Omega \times \Gamma _{d}$ as above, then
we set $Rf(r,\theta )=u(\mathbf{x},\mathbf{x}_{\alpha }).$

\item \textbf{Case 2}. Otherwise, we set $Rf(r,\theta )=0.$
\end{enumerate}

In our computations, discrete values for the variable $r$ are: $%
-l/2+(i-1)l/216$, $i=0,\dots ,216.$ And discrete values for the variable $%
\theta $ are: $\{0^{\circ },1^{\circ },\dots ,179^{\circ }\}$.

After computing the incomplete $Rf$ from our data, we use the command
\textquotedblleft iradon" of MATLAB to reconstruct $f$. The discrete
function $f$ \ computed by the procedure in this section is denoted as $f_{%
\mathrm{comp}}^{\mathrm{iradon}}.$

\subsection{Post processing}

\label{sec post}

We need to \textquotedblleft clean up" the obtained results. To do this, we
perform the following two post processing steps:

\begin{enumerate}
\item \textbf{Step 1}. Let $f\left( \mathbf{x}\right) $ be either $f_{%
\mathrm{comp}}$ or $f_{\mathrm{comp}}^{\mathrm{iradon}}.$ We observe that
the image of $f\left( \mathbf{x}\right) $ contains unwanted artifacts. We
remove these artifacts by a simple procedure. Let $m=\displaystyle\max_{%
\mathbf{x}\in \overline{\Omega }}\{|f(\mathbf{x})|\}.$ We set 
\begin{equation}
\widetilde{f}(\mathbf{x})=\left\{ 
\begin{array}{ll}
0 & \text{ if }|f(\mathbf{x})|<0.2m, \\ 
f(\mathbf{x}) & \text{ otherwise.}%
\end{array}%
\right.  \label{eqn filter}
\end{equation}%
Next, we smooth out the function $\widetilde{f}(\mathbf{x})$ as in Step 2.
For brevity, we keep below the same notation $f\left( \mathbf{x}\right) $
for $\widetilde{f}(\mathbf{x})$. In Section \ref{sec num}, we display the
computed functions $f$ before and after using this 20\% artifact remover.

\item \textbf{Step 2}. Due to the presence of noise in the data , we have to
smooth the computed vector $\mathbf{U}_{\mathrm{comp}}^{h}$, the computed
function $f_{\mathrm{comp}}$ and $f_{\mathrm{comp}}^{iradon}$. More
precisely, for each $n\in \{0,\dots ,N\}$ and for each grid point $\mathbf{x}%
\in \Omega $, the number $u_{n}^{\mathrm{comp}}(\mathbf{x})$ is replaced by
the mean value of $u_{n}^{\mathrm{comp}}$ over neighboring grid points
located in the closed rectangle of the size $7h _{x}\times 7h _{y}$
centered at $\mathbf{x.}$ Only those grid points are counted which are
located in $\overline{\Omega }$. This smoothing step is applied to both
functions $f_{\mathrm{comp}}$ and $f_{\mathrm{comp}}^{\mathrm{iradon}}$. The
number $7$ here is chosen by a trial and error process.
\end{enumerate}

The numerical implementation of our approach to compute the functions $f_{%
\mathrm{comp}}$ is summarized in Algorithm \ref{alg 2}.

\begin{algorithm}
\caption{The numerical implementation to solve Problem \ref{pro 2.1}.} \label%
{alg 2} 
\begin{algorithmic}[1]
		\State{Choose $N = 15$.}
		\State{Calculate ${\bf g}(x, y) = (u_1(x, y), \dots, u_n(x, y))$ for $(x, y) \in \partial \Omega$ via formula \eqref{eqn boundary data}}.	
		\State{Compute the matrices $A(x, y)$ and $B(x, y)$ for $(x, y) \in \Omega$ by \eqref{eqn AB}.}			
	\State{Solve Problem \ref{pro 3} by the command ``linsolve" of MATLAB.}
	\State{Compute $u_{\rm comp}(x, y)$ using \eqref{4.2}.}
	\State{Compute $f_{\rm comp}(x, y)$ using \eqref{5.32}.}
	\State{Apply the above
post processing procedure.}
	\end{algorithmic}
\end{algorithm}

\begin{remark}
We have computationally observed that if the chosen number $N$ is small,
then resulting images are of not a good quality. On the other hand, if $N$
is too large, then our Algorithm \ref{alg 2} is time consuming. Here, we
have chosen $N=15$ in our computations by a trial and error procedure. We
have observed that with this choice of $N,$ the numerical results are
stable. Furthermore, numerical results change insignificantly when $N$
grows. \label{rem 5.1}
\end{remark}

\section{Numerical tests}

\label{sec num}

We test two cases (Tests 1,2) in which the true function $f^{\ast }$
consists of inclusions of the circular shape. In addition, we test two more
cases, in which the true functions $f^{\ast }$ are the characteristic
functions of some non-convex domains. To work with the first two cases, we
use a template which is a circular inclusion of the radius $1$ centered at $0
$ and described by the function $\varphi (\mathbf{x})$ 
\begin{equation}
\varphi (\mathbf{x})=\left\{ 
\begin{array}{ll}
\exp (-|\mathbf{x}|^{2}/(1-|\mathbf{x}|^{2})) & |\mathbf{x}|<1, \\ 
0 & \mbox{otherwise}%
\end{array}%
\right.   \label{7.1}
\end{equation}%
So, in the first two cases the function $f^{\ast }$ is generated by the
function $\varphi (\mathbf{x})$ in (\ref{7.1}) in which some parameters are
involved. These parameters provide a linear combination, scale and/or a
translation of the above circular inclusion. The functions $f^{\ast }$ in
the latter two cases are the characteristic functions of subsets of a
rhombus centered at the center of $\Omega $.

We consider four numerical tests listed below. In the case of the filtered
back projection method we add 5\% noise to the data, as in \eqref{6.30}.
However, in the case of our method we first add 5\% and then 15\% noise. In
all tests $d$, in (\ref{2.10}) is set to be $3.5$.

\begin{enumerate}
\item \textbf{Test 1}. The true function $f^{\ast }$ is given by 
\begin{equation*}
f^{\ast }(\mathbf{x})=\varphi ((\mathbf{x}-\mathbf{x}_{0})/r)
\end{equation*}%
where $r=0.2$ and $\mathbf{x}_{0}=(0,2)$ is the the center of $\Omega
=(-1,1)\times (1,3)$. In this setting, the distance between the source line $%
\Gamma _{d}$ in (\ref{2.10}) and the domain $\Omega $ is 1. Keeping in mind
comparison with the filtered back projection method, we number this as
\textquotedblleft inclusion number 1". The numerical result is displayed in
Figure \ref{fig case 1}.

\begin{figure}
	\begin{center}
		\subfloat[The true function $f^*$]{\includegraphics[width=.3\textwidth]{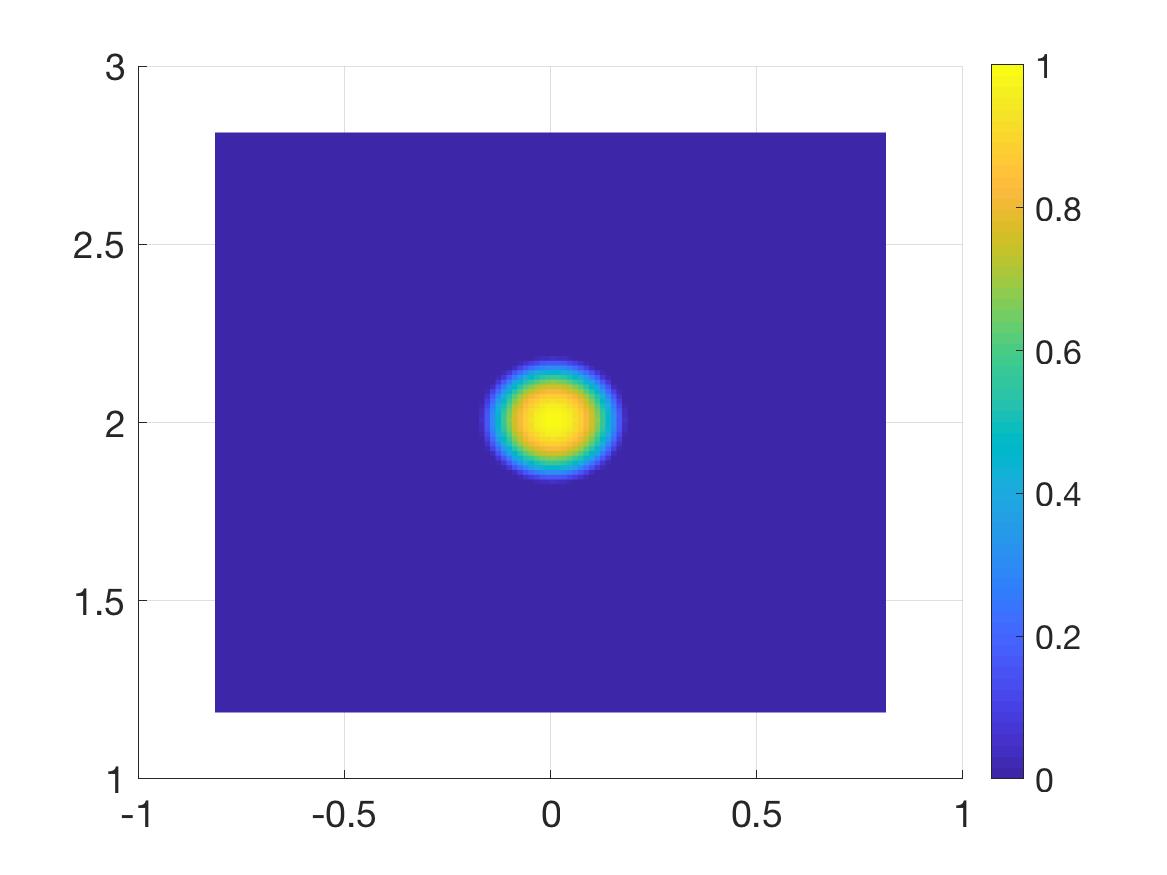}} \hspace{.2cm}
		\subfloat[\label{fig Case 1 full}The Radon transform of $f^*$ computed by the function ``radon" of Matlab]{\includegraphics[width=.3\textwidth]{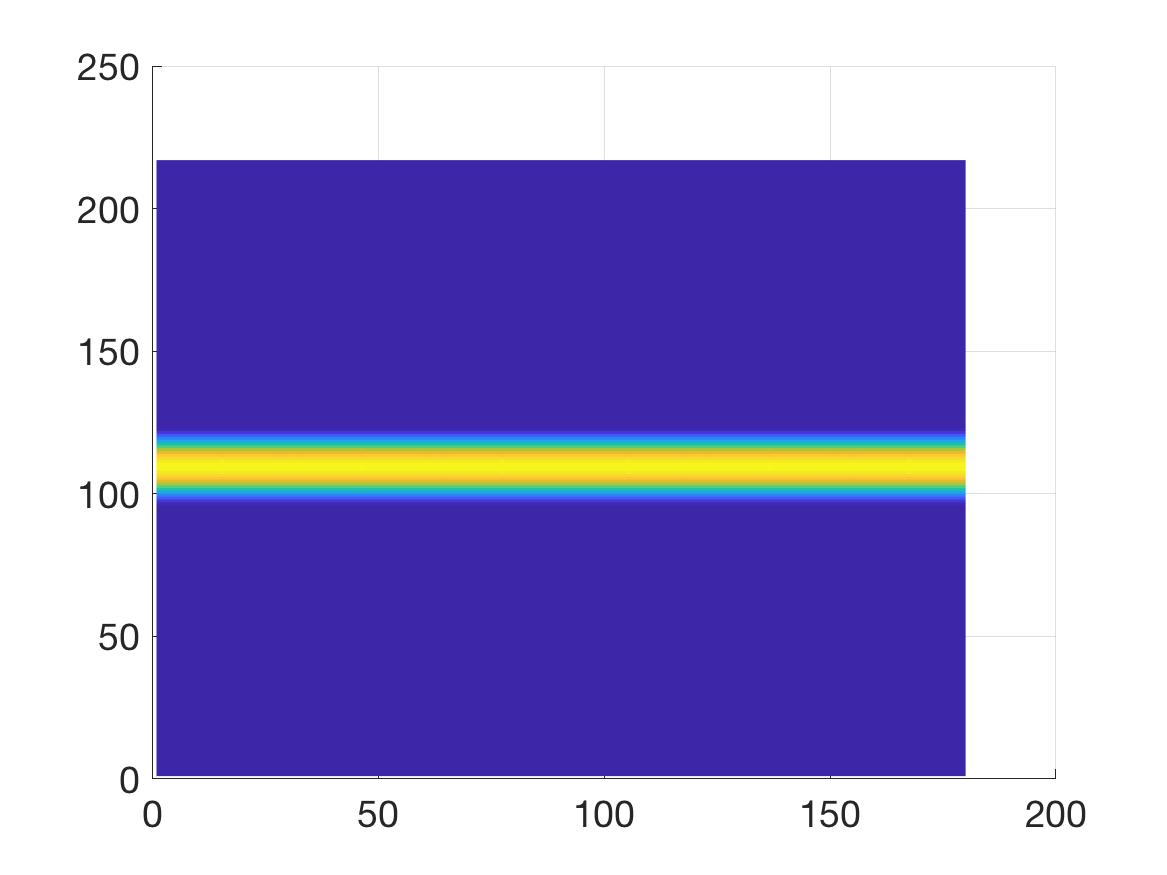}} \hspace{.2cm}
		\subfloat[\label{fig Case 1 incomplete} The incomplete tomographic data with $5\%$ noise]{\includegraphics[width=.3\textwidth]{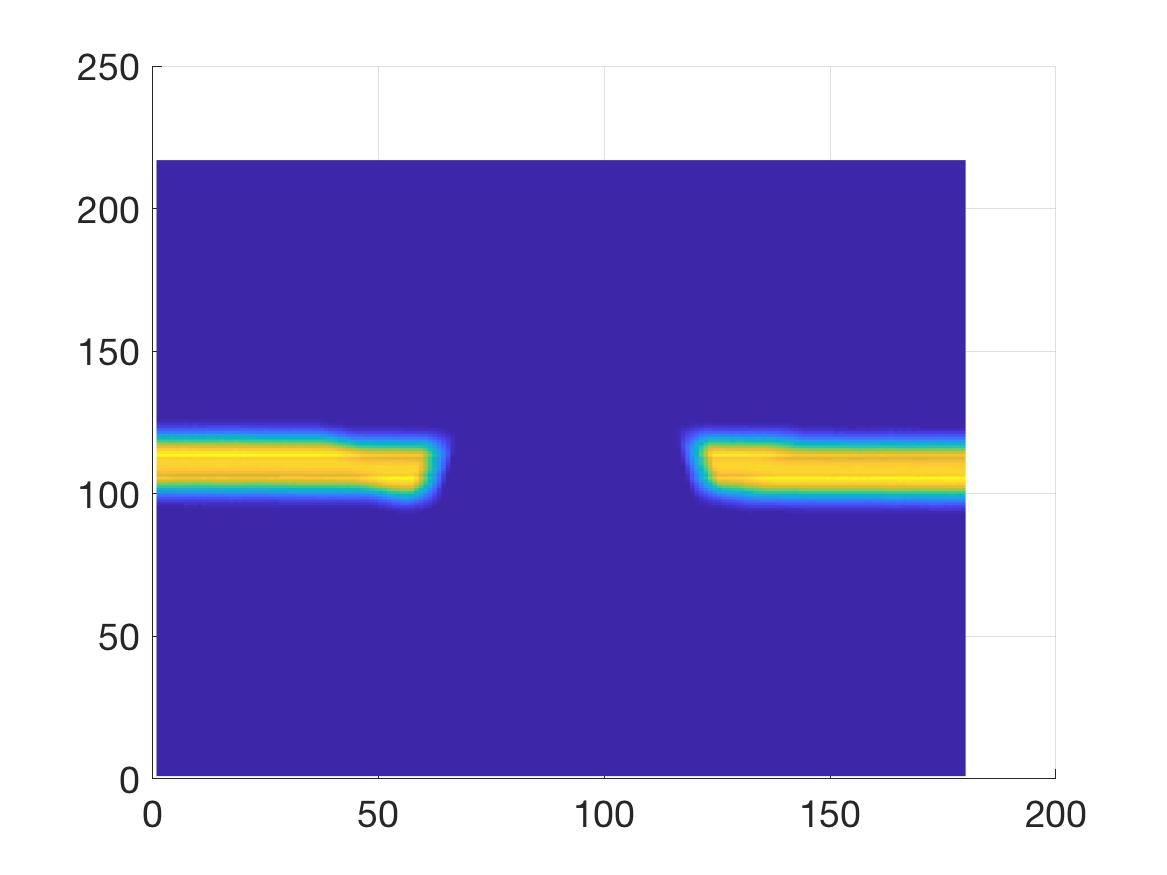}} 
		
		\subfloat[\label{fig Case 1 Radon No Clean}The function $f_{\rm comp}^{\rm iradon}$ computed by the filtered back projection algorithm, noise level $5\%$]{\includegraphics[width=.3\textwidth]{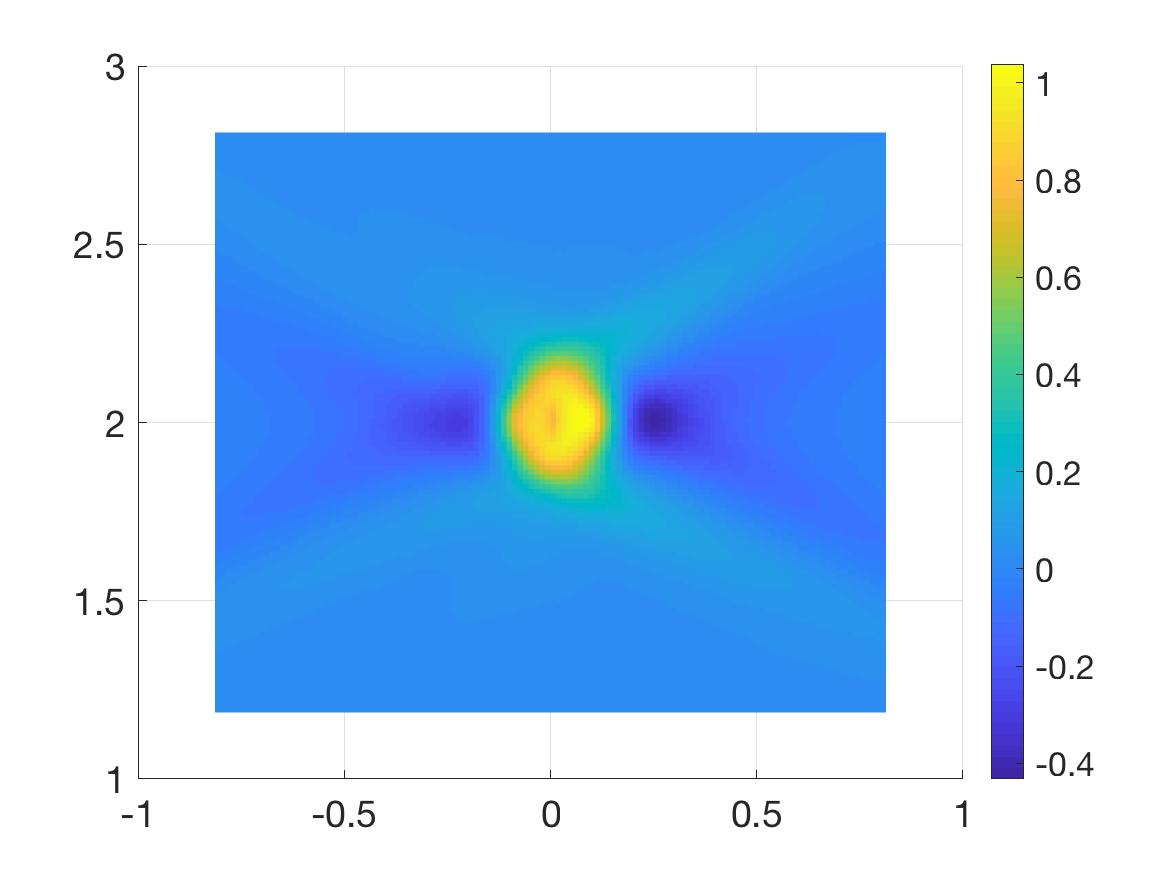}} \hspace{.2cm}
		\subfloat[\label{fig Case 1 Radon Clean}The function $f_{\rm comp}^{\rm iradon}$ computed by the filtered back projection algorithm, noise level $5\%$, together with the post processing of Section \ref{sec post}]{\includegraphics[width=.3\textwidth]{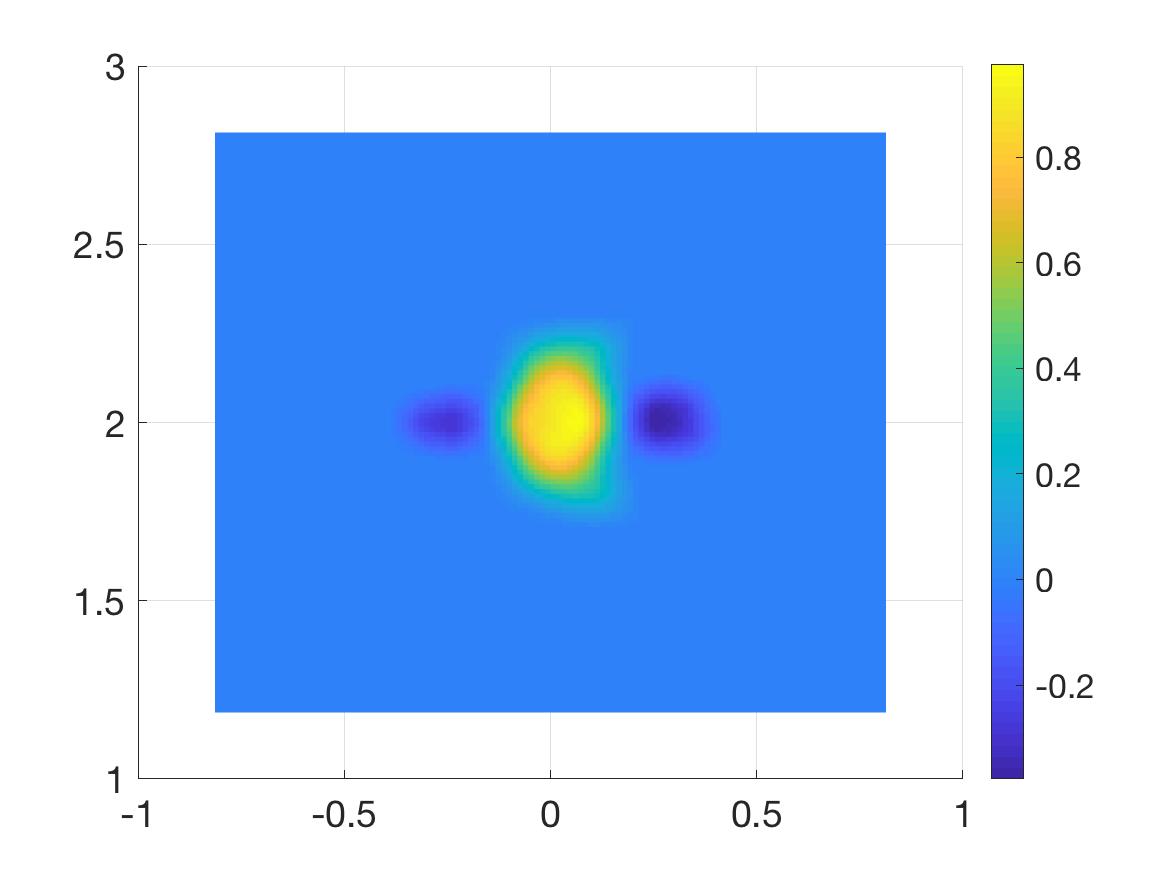}} \hspace{.2cm}
		\subfloat[\label{fig Case 1 with 5 percent noise}The function $f_{\rm comp}$ by our method in Section \ref{sec our approach}, noise level $5\%$]{\includegraphics[width=.3\textwidth]{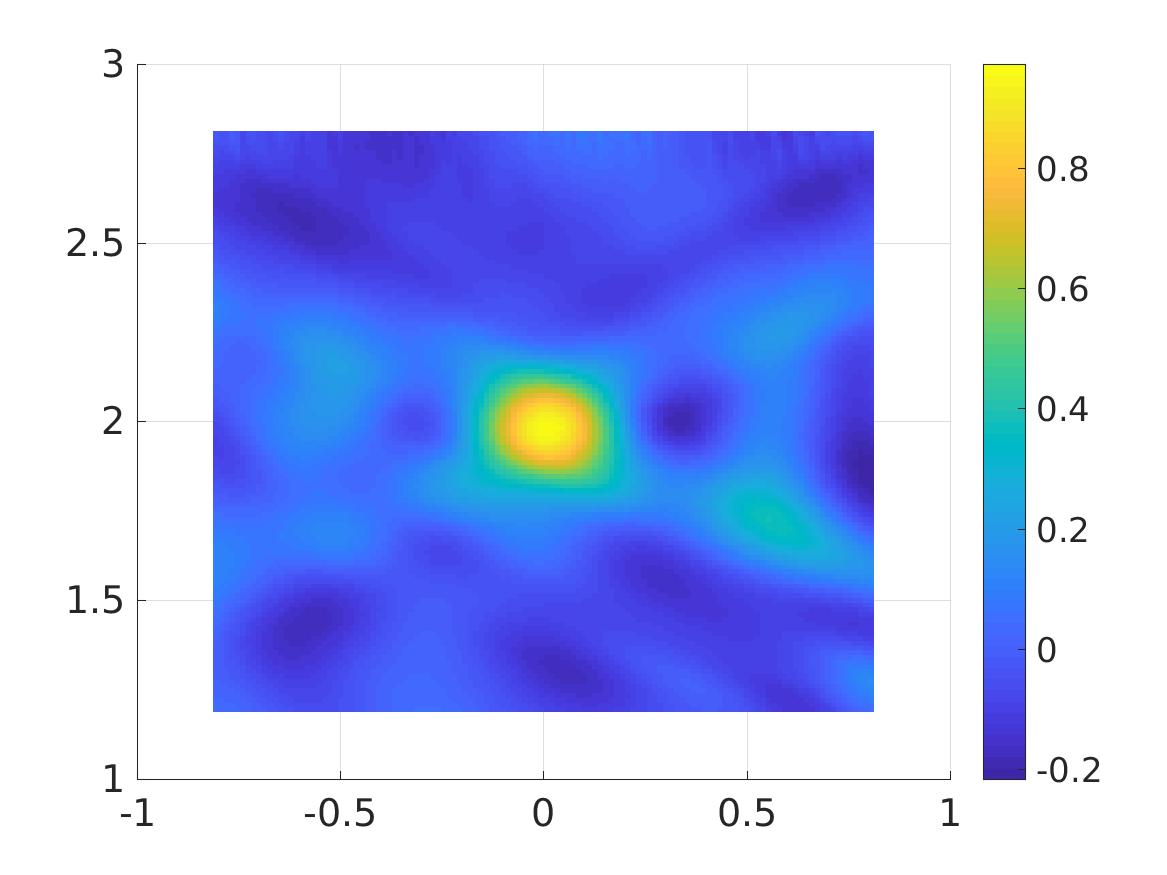}} 
		
		\subfloat[\label{fig Case 1 with 5 percent noise Clean} The function $f_{\rm comp}$ by our method in Section \ref{sec our approach}, noise level $5\%$, together withpost processing of Section \ref{sec post}]{\includegraphics[width=.3\textwidth]{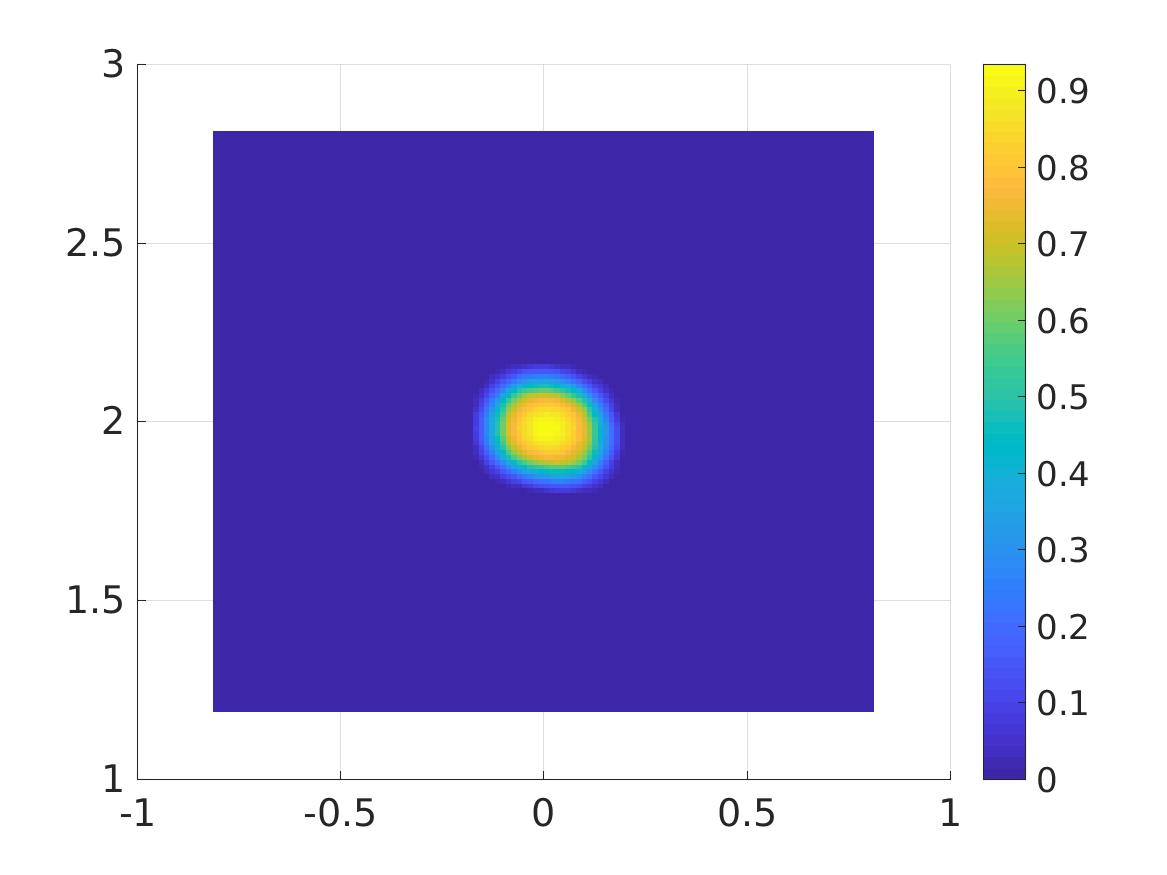}}  \hspace{.2cm}
		\subfloat[\label{fig Case 1 with 15 percent noise}The function $f_{\rm comp}$ by our method in Section \ref{sec our approach}, noise level $15\%$]{\includegraphics[width=.3\textwidth]{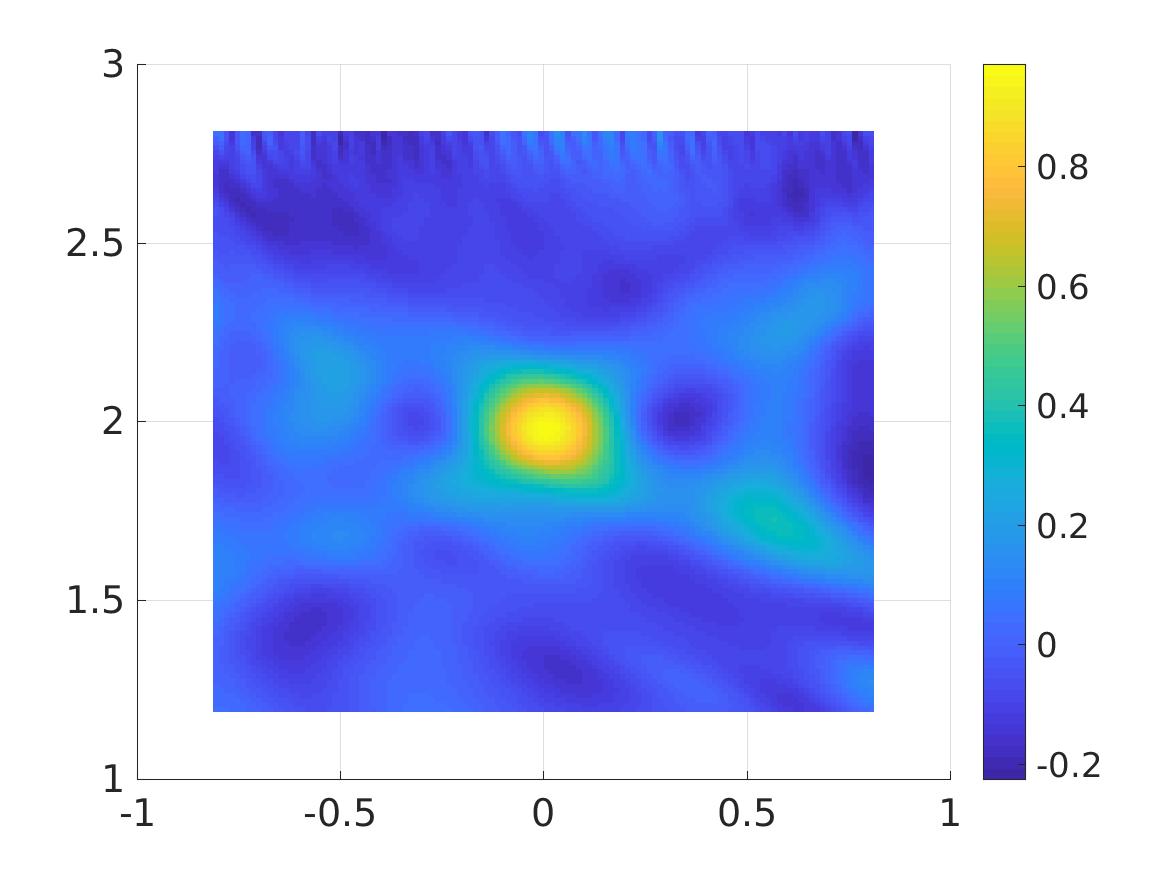}}  \hspace{.2cm}
		\subfloat[\label{fig Case 1 with 15 percent noise Clean}The function $f_{\rm comp}$ by our method in Section \ref{sec our approach}, noise level $15\%$, post processing of Section \ref{sec post}]{\includegraphics[width=.3\textwidth]{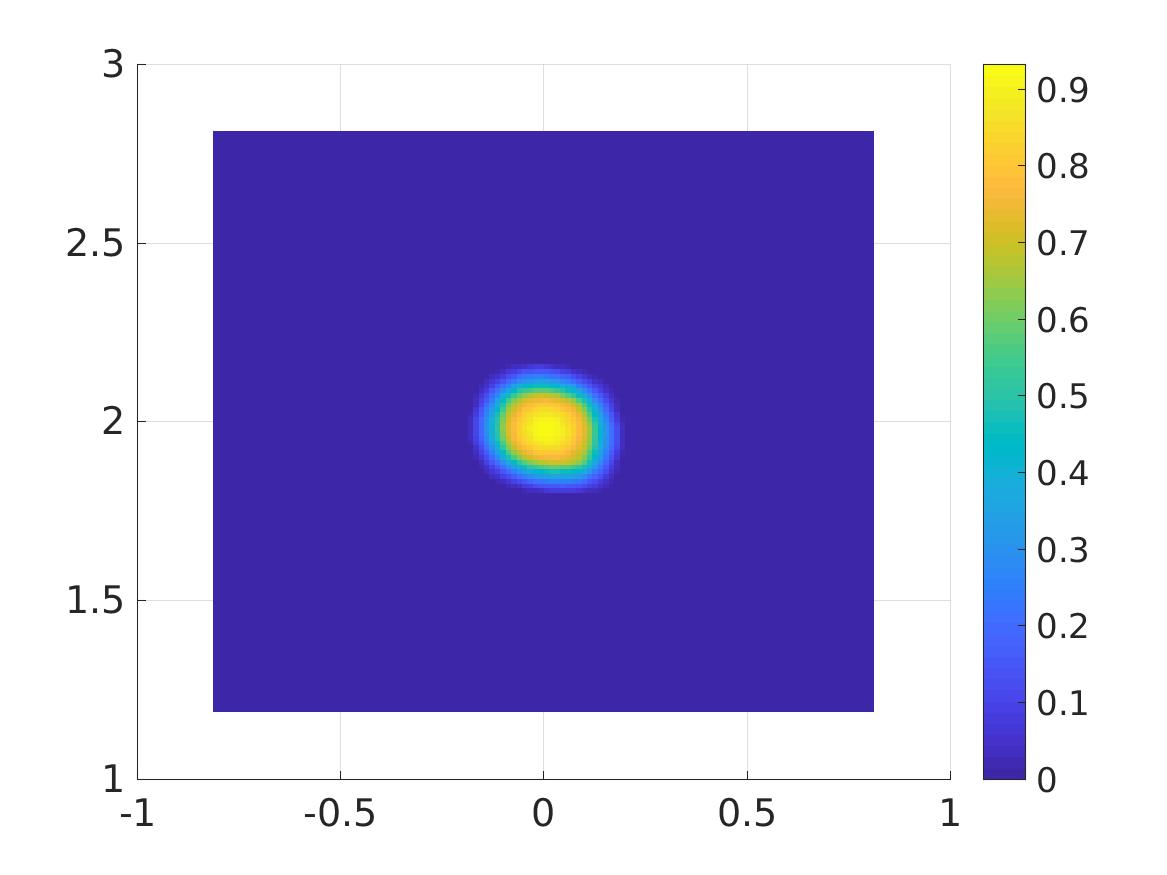}}  \hspace{.2cm}

\caption{\it  \label{fig case 1}
Test 1, inclusion number 1. The data and the reconstructions of the function
$f^{\ast }.$ One can see from (e),(g),(i) that the image quality
provided by our method is slightly better than that of the filtered back
projection method.}		
\end{center}
\end{figure}

\item \textbf{Test 2}. In this case, we set $\Omega =(-1,1)\times (3,5)$.
The distance between the source line $\Gamma _{d}$ and the domain $\Omega $
is 3, which is three times larger than the distance in test 1. The true
function $f^{\ast }$ is 
\begin{equation*}
f^{\ast }(\mathbf{x})=-6\varphi ((\mathbf{x}-\mathbf{x}_{1})/r_{1})+5\varphi
((\mathbf{x}-\mathbf{x}_{2})/r_{2})+6\varphi ((\mathbf{x}-\mathbf{x}%
_{3})/r_{3})
\end{equation*}%
where $\mathbf{x}_{1}=(-0.4,4),r_{1}=0.2$, $\mathbf{x}%
_{2}=(-0.1,3.57),r_{2}=0.23$, $\mathbf{x}_{3}=(0.4,4)$ and $r_{3}=0.18.$
Hence, we have here three different radii of circles varying between 0.18
and 0.23. We note that $f\left( \mathbf{x}\right) \leq 0$ inside of the
first circle, and $f\left( \mathbf{x}\right) \geq 0$ inside of second and
third circles. We number these inclusions as \textquotedblleft inclusions
number 1, 2 and 3" respectively. The true and reconstructed functions $f$
are displayed in Figure \ref{fig case 3}.

\begin{figure}
	\begin{center}
		\subfloat[The true function $f^*$]{\includegraphics[width=.3\textwidth]{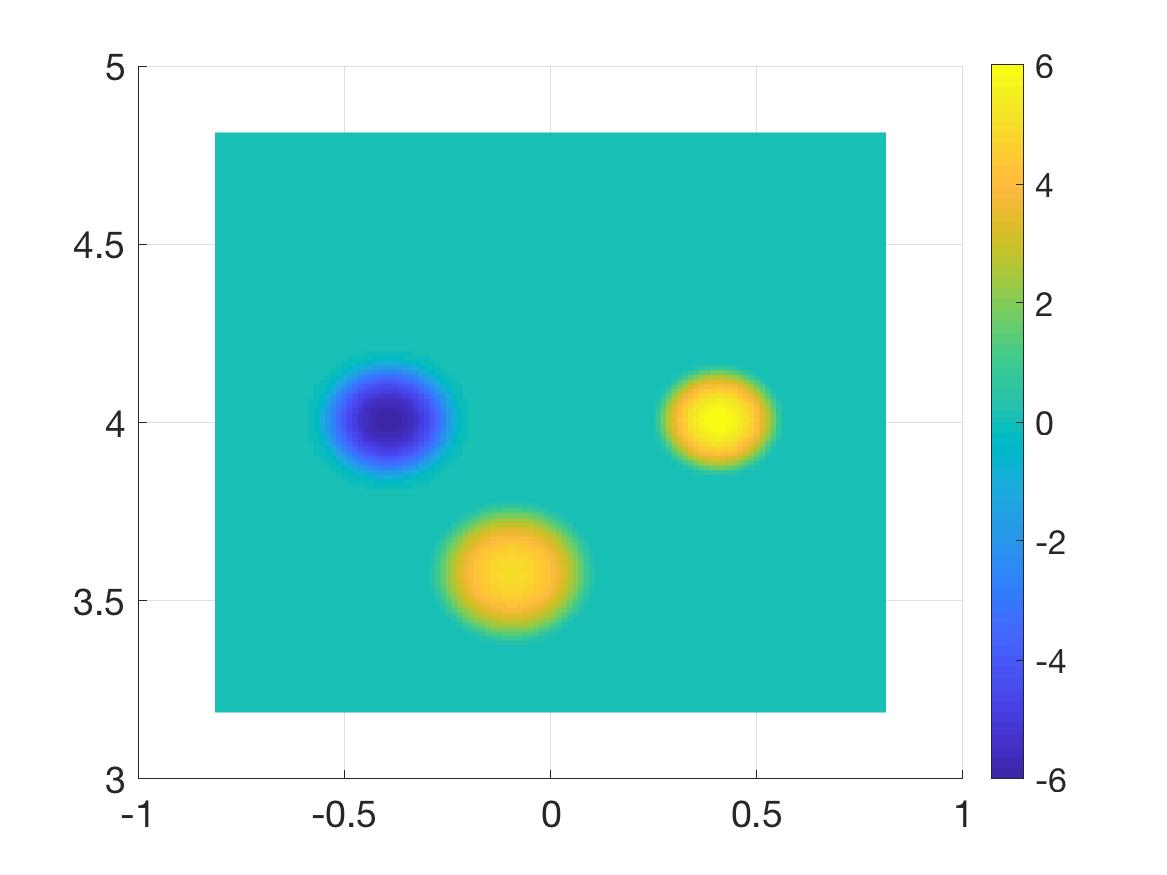}} \hspace{.2cm}
		\subfloat[\label{fig Case 3 full}The Radon transform of $f^*$ computed by the function ``radon" of Matlab]{\includegraphics[width=.3\textwidth]{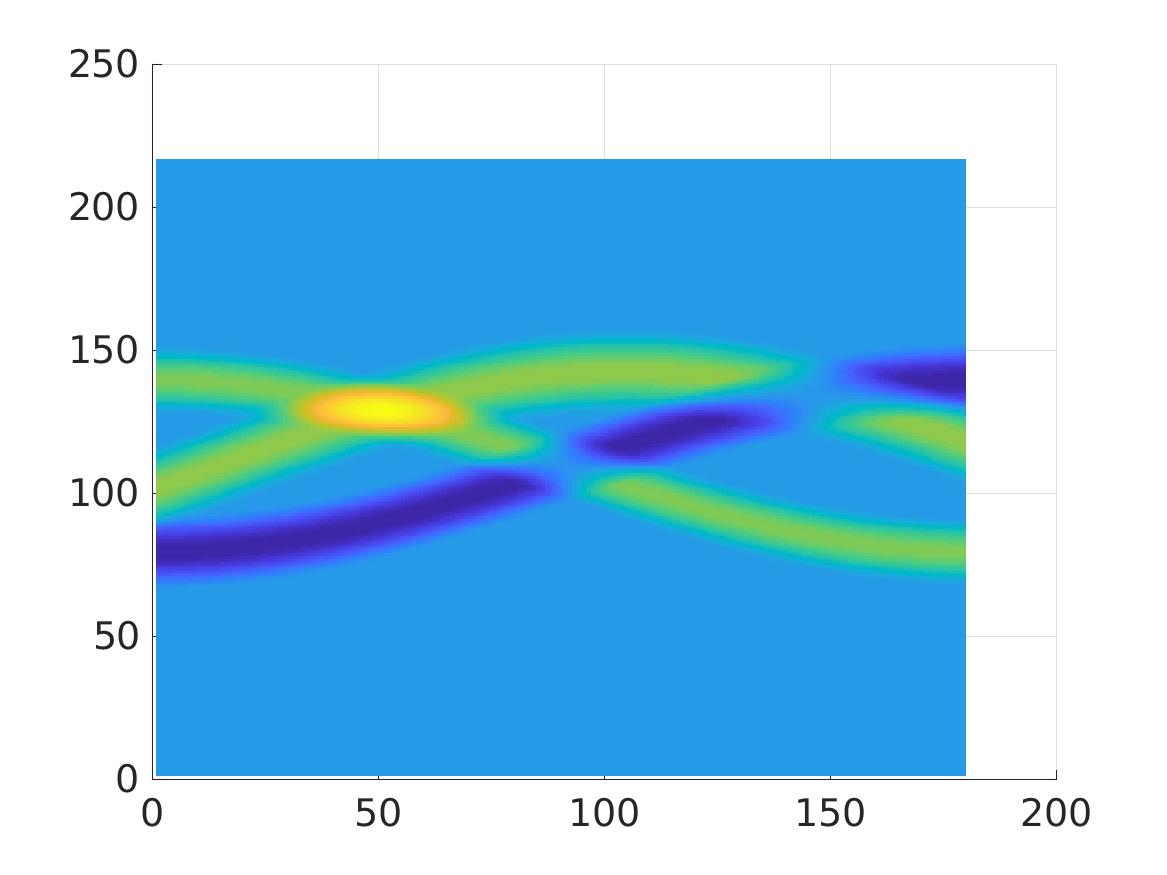}} \hspace{.2cm}
		\subfloat[\label{fig Case 3 incomplete} The incomplete tomographic data with $5\%$ noise]{\includegraphics[width=.3\textwidth]{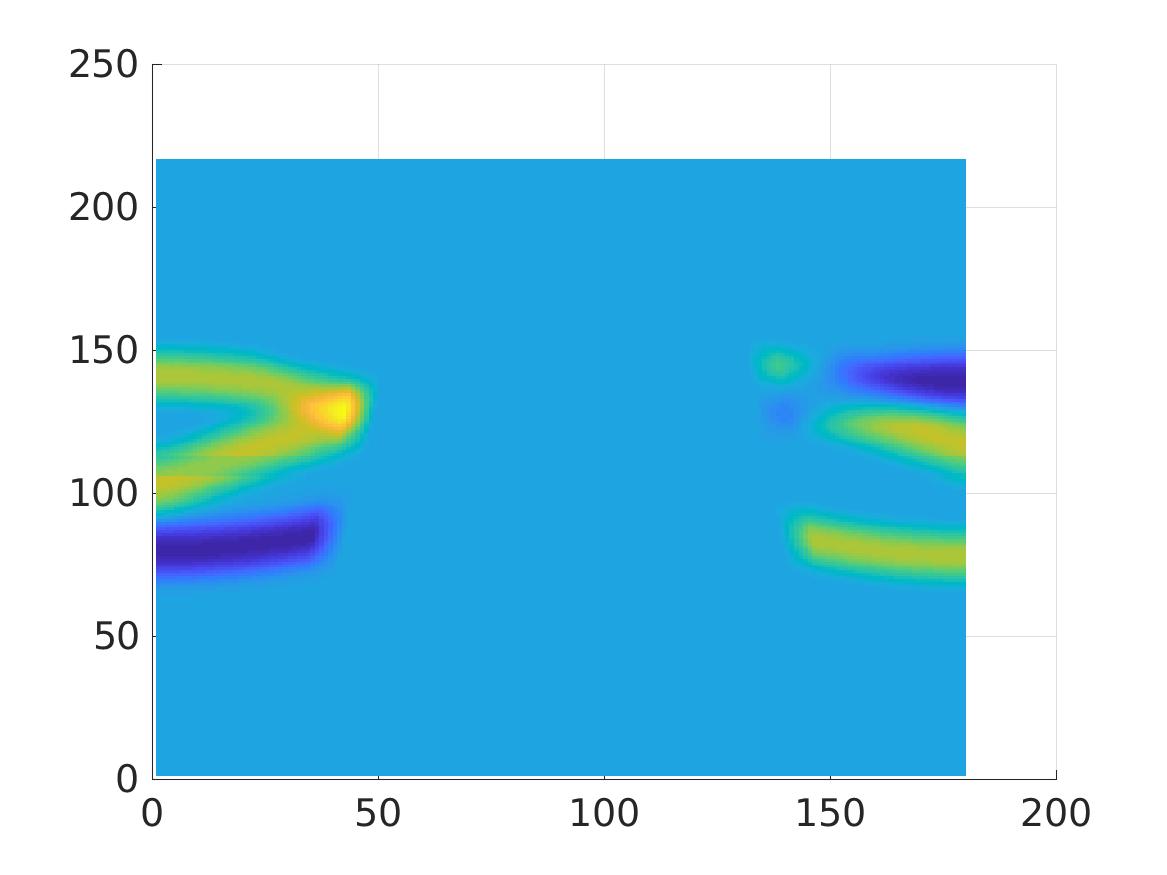}} 
		
		\subfloat[\label{fig Case 3 Radon No Clean}The function $f_{\rm comp}^{\rm iradon}$ computed by the filtered back projection algorithm, noise level $5\%$]{\includegraphics[width=.3\textwidth]{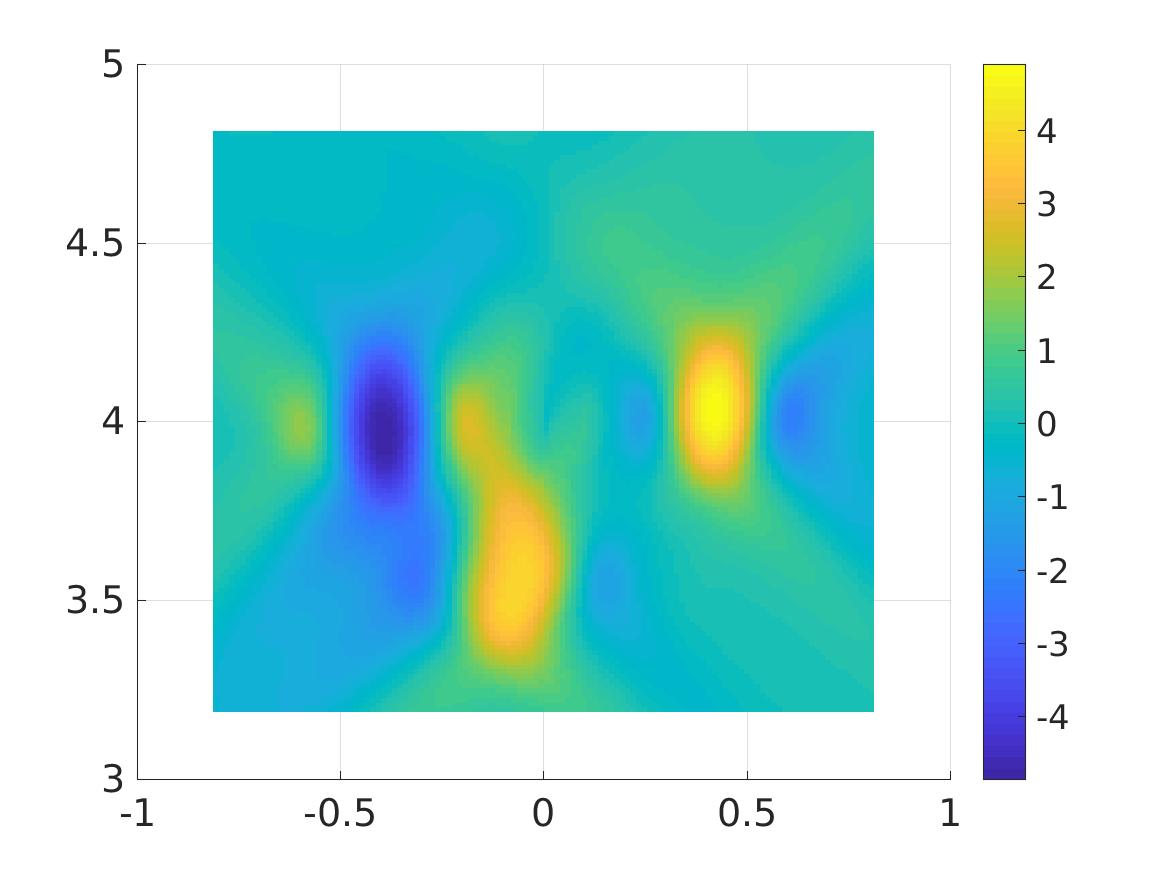}} \hspace{.2cm}
		\subfloat[\label{fig Case 3 Radon Clean}The function $f_{\rm comp}^{\rm iradon}$ computed by the filtered back projection algorithm, noise level $5\%$, together with the post processing of Section \ref{sec post}]{\includegraphics[width=.3\textwidth]{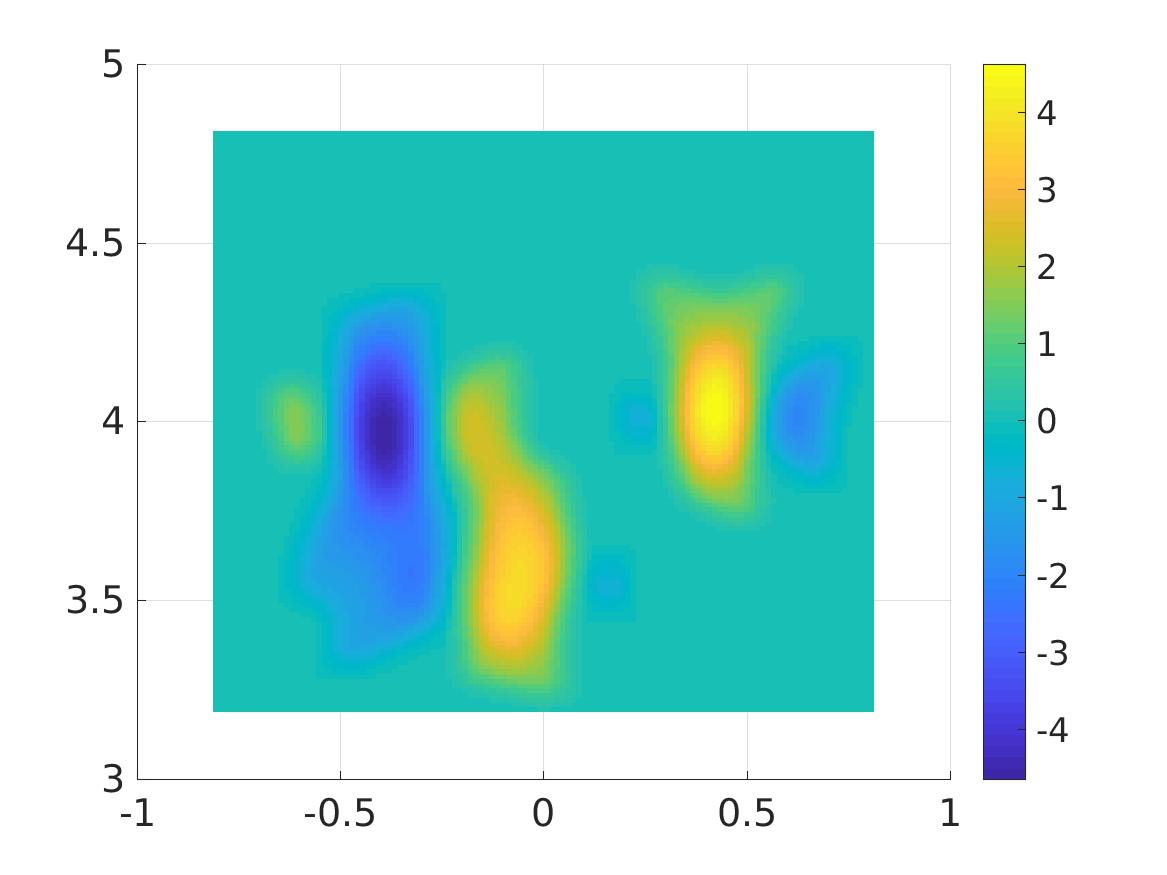}} \hspace{.2cm}
		\subfloat[\label{fig Case 3 with 5 percent noise}The function $f_{\rm comp}$ by our method in Section \ref{sec our approach}, noise level $5\%$]{\includegraphics[width=.3\textwidth]{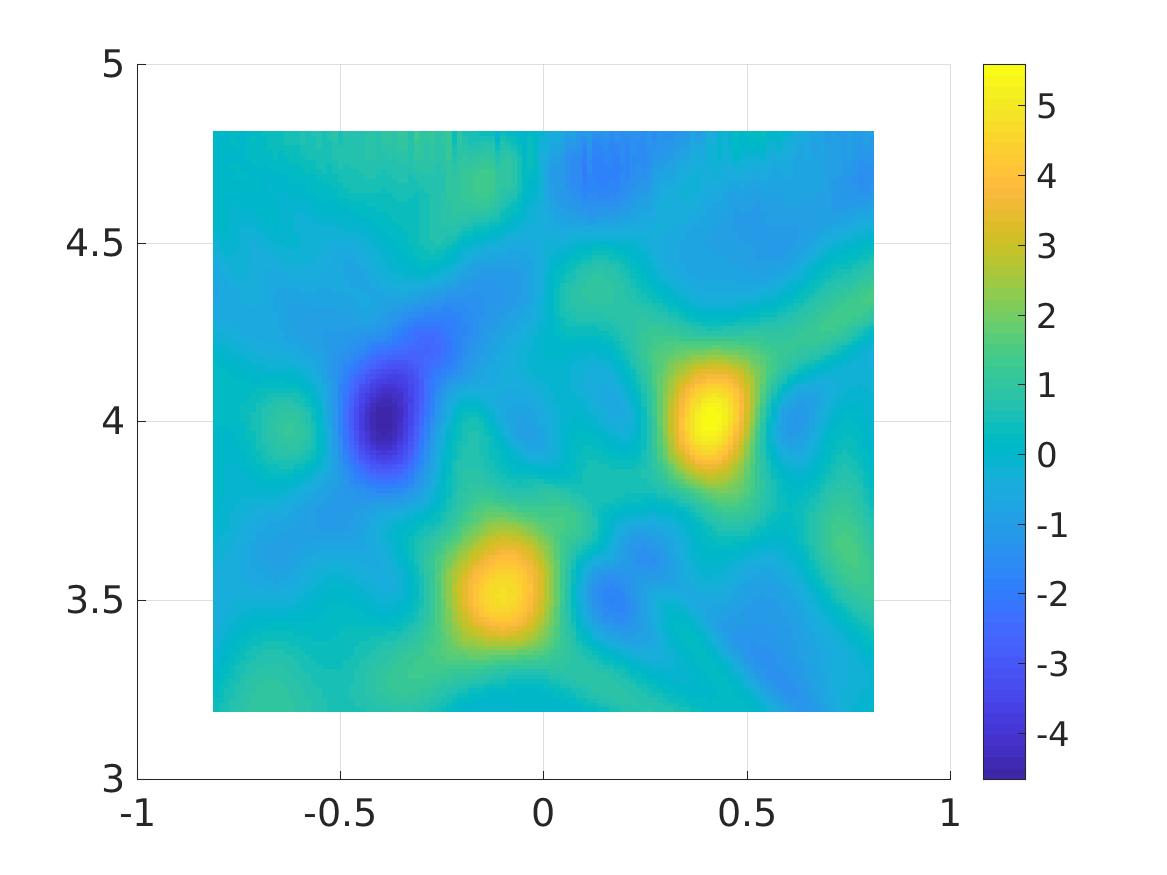}} 
		
		\subfloat[\label{fig Case 3 with 5 percent noise Clean} The function $f_{\rm comp}$ by our method in Section \ref{sec our approach}, noise level $5\%$, together withpost processing of Section \ref{sec post}]{\includegraphics[width=.3\textwidth]{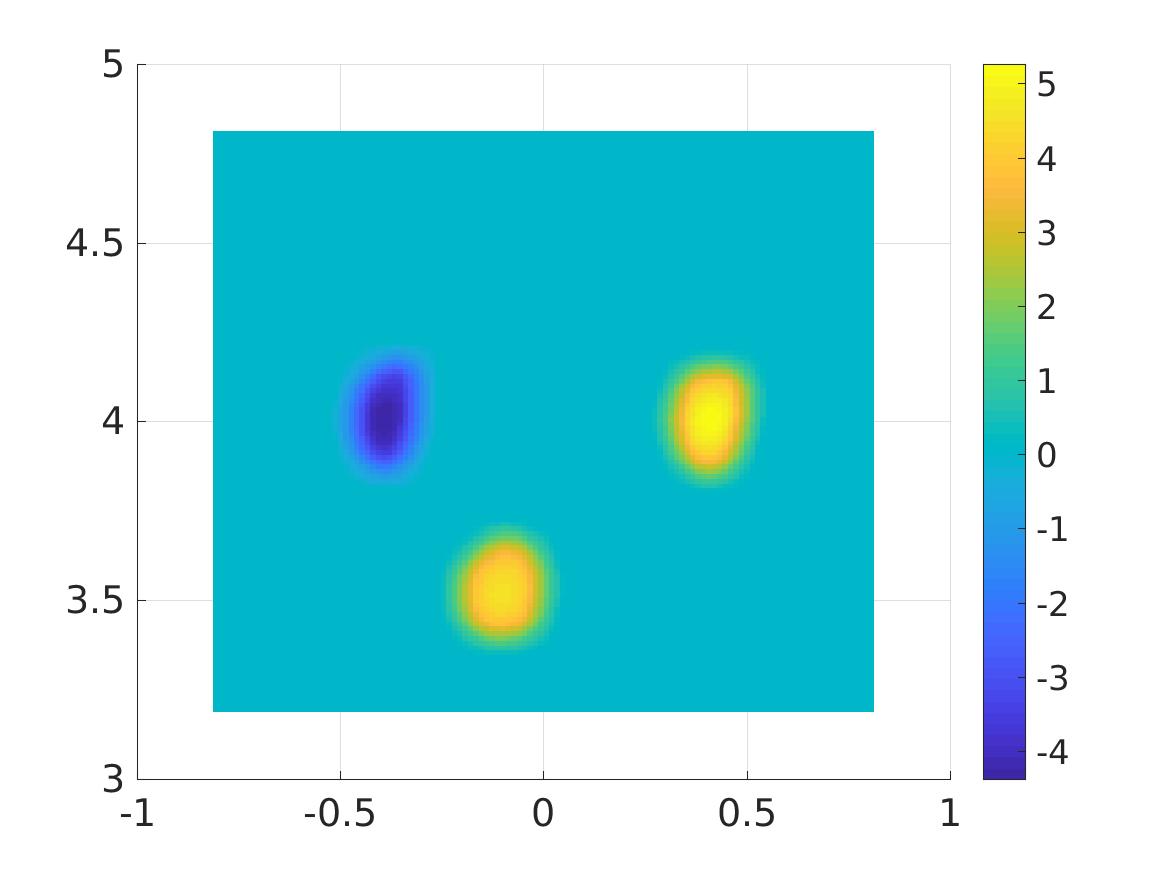}}  \hspace{.2cm}
		\subfloat[\label{fig Case 3 with 15 percent noise}The function $f_{\rm comp}$ by our method in Section \ref{sec our approach}, noise level $15\%$]{\includegraphics[width=.3\textwidth]{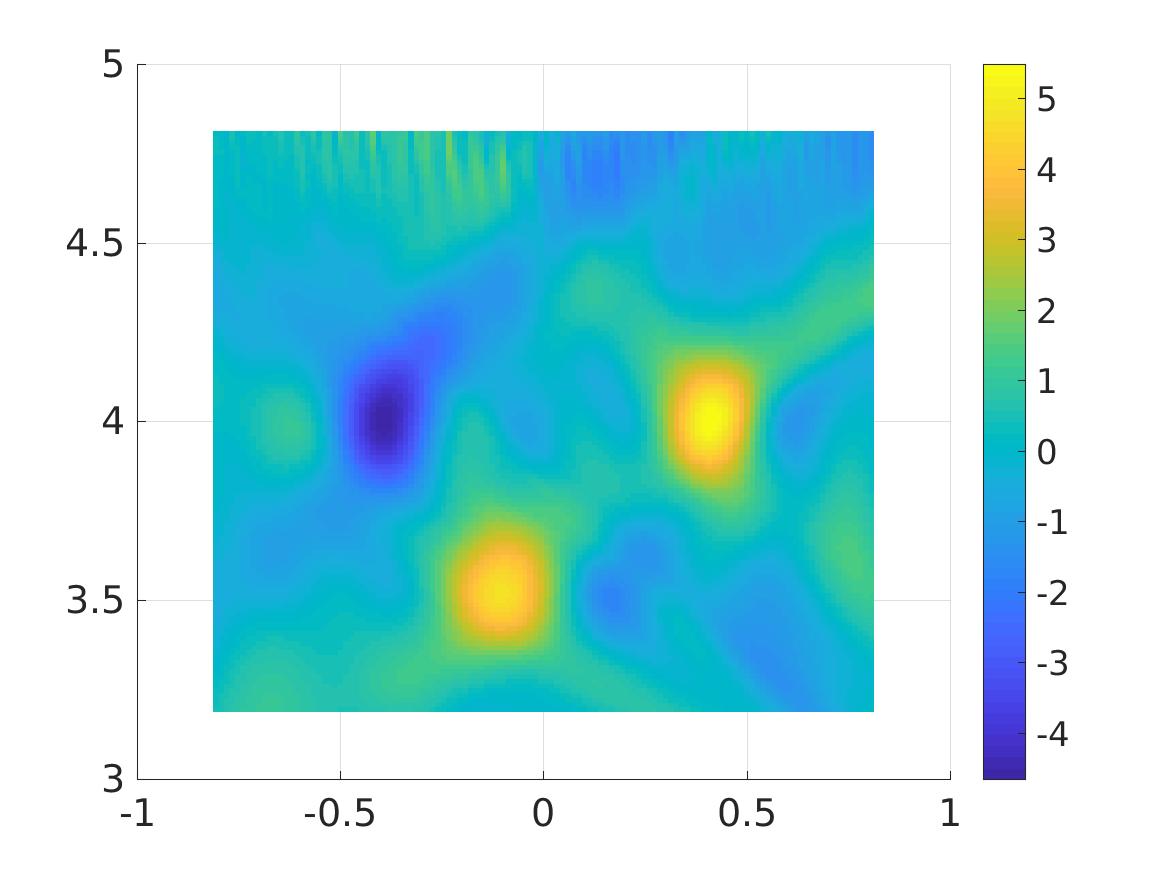}}  \hspace{.2cm}
		\subfloat[\label{fig Case 3 with 15 percent noise Clean}The function $f_{\rm comp}$ by our method in Section \ref{sec our approach}, noise level $15\%$, post processing of Section \ref{sec post}]{\includegraphics[width=.3\textwidth]{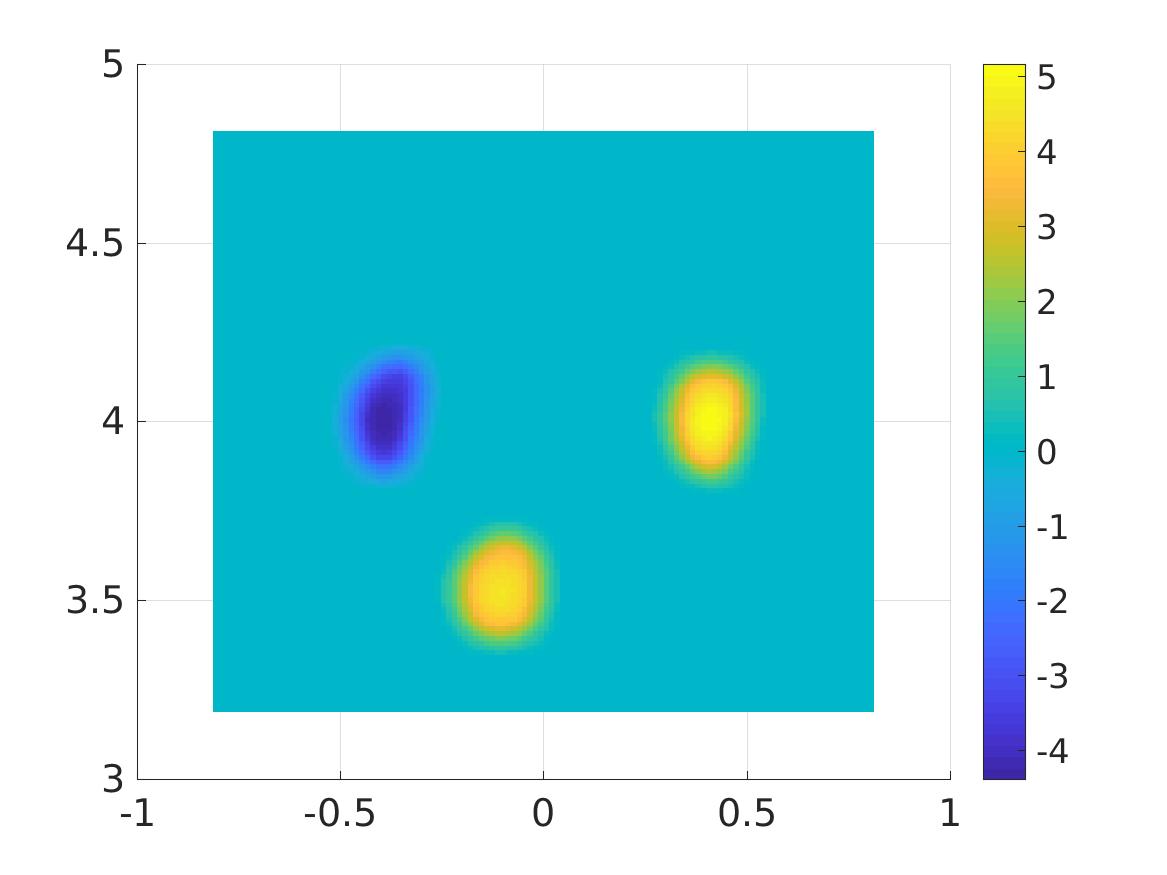}}  \hspace{.2cm}

\caption{\it  \label{fig case 3}
Test 2. The data and the reconstructions of the function $f^{\ast }$
in the case of three inclusions. On (a),(d)-(i)
inclusions  from left to right are numbered as 2,3 and 4. One can see from (e),(g),(i) that the
image quality provided by our method is better than that of the filtered
back projection method.}		
\end{center}
\end{figure}

\item \textbf{Test 3}. Next, we test a non smooth function and the inclusion
whose shape is not circular. Set $\Omega =(-1,1)\times (3.5,5.5).$ The
distance between the source line $\Gamma _{d}$ and the domain $\Omega $ is
now 3.5, which is greater than in previous two tests. The true function $%
f^{\ast }$ is 
\begin{equation*}
f^{\ast }(\mathbf{x})=\chi _{\{\mathbf{x}%
=(x,y):0.3<|x|+|y-4.5|<0.6,x>0.3,y>4.5\}},
\end{equation*}%
where $\chi $ is the characteristic function. The image of the true
inclusion looks like a letter $L$ rotated clockwise by $3\pi /4$ around the
center of $\Omega $. The true and reconstructed functions $f$ are displayed
in Figure \ref{fig case 6}.

\begin{figure}
	\begin{center}
		\subfloat[The true function $f^*$]{\includegraphics[width=.3\textwidth]{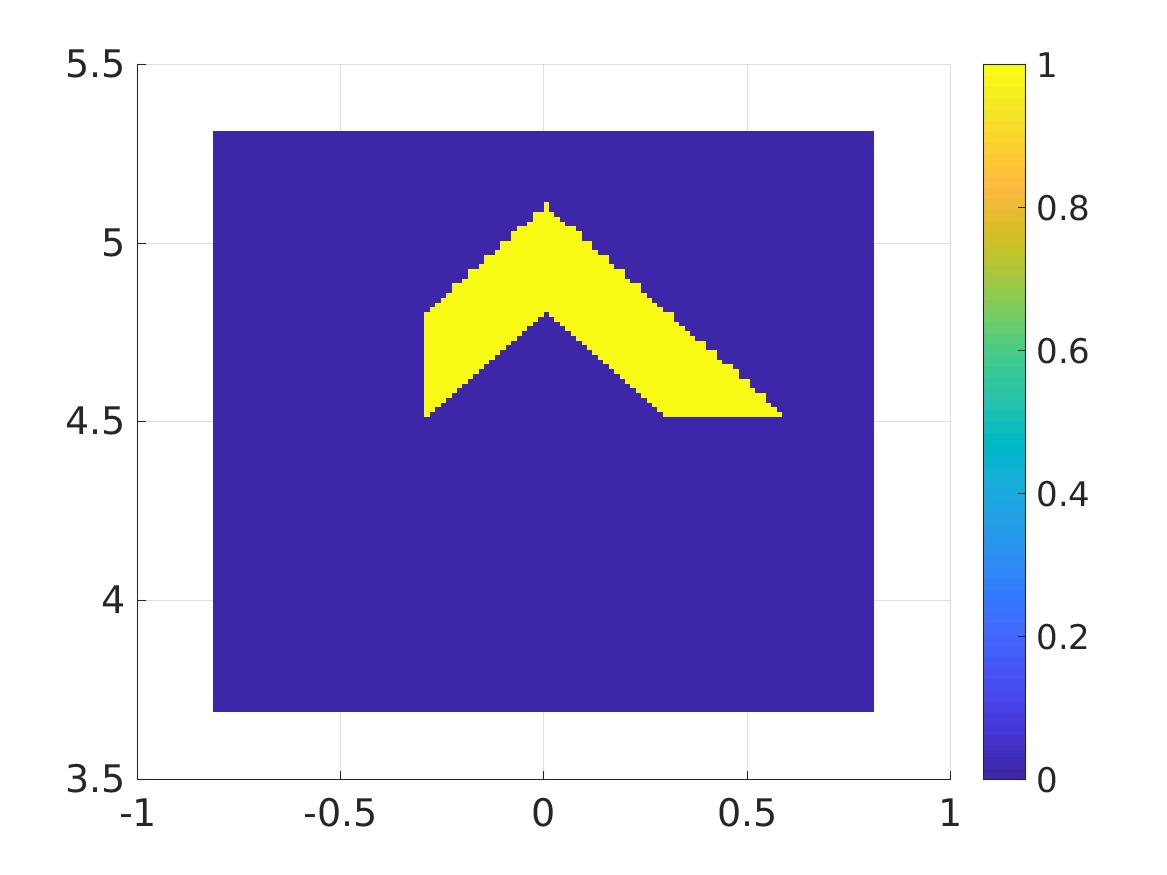}} \hspace{.2cm}
		\subfloat[The Radon transform of $f^*$ computed by the function ``radon" of Matlab]{\includegraphics[width=.3\textwidth]{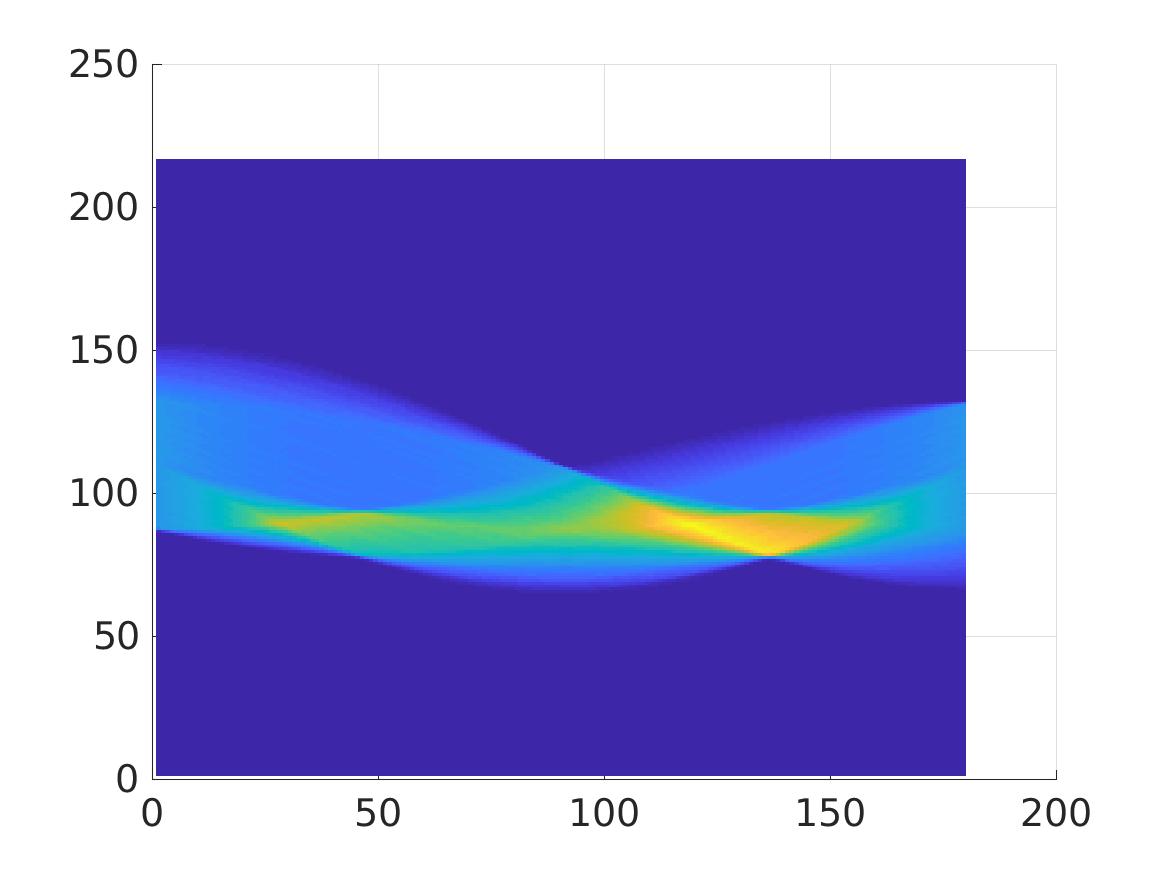}} \hspace{.2cm}
		\subfloat[\label{fig Case 6 incomplete} The incomplete tomographic data with $5\%$ noise]{\includegraphics[width=.3\textwidth]{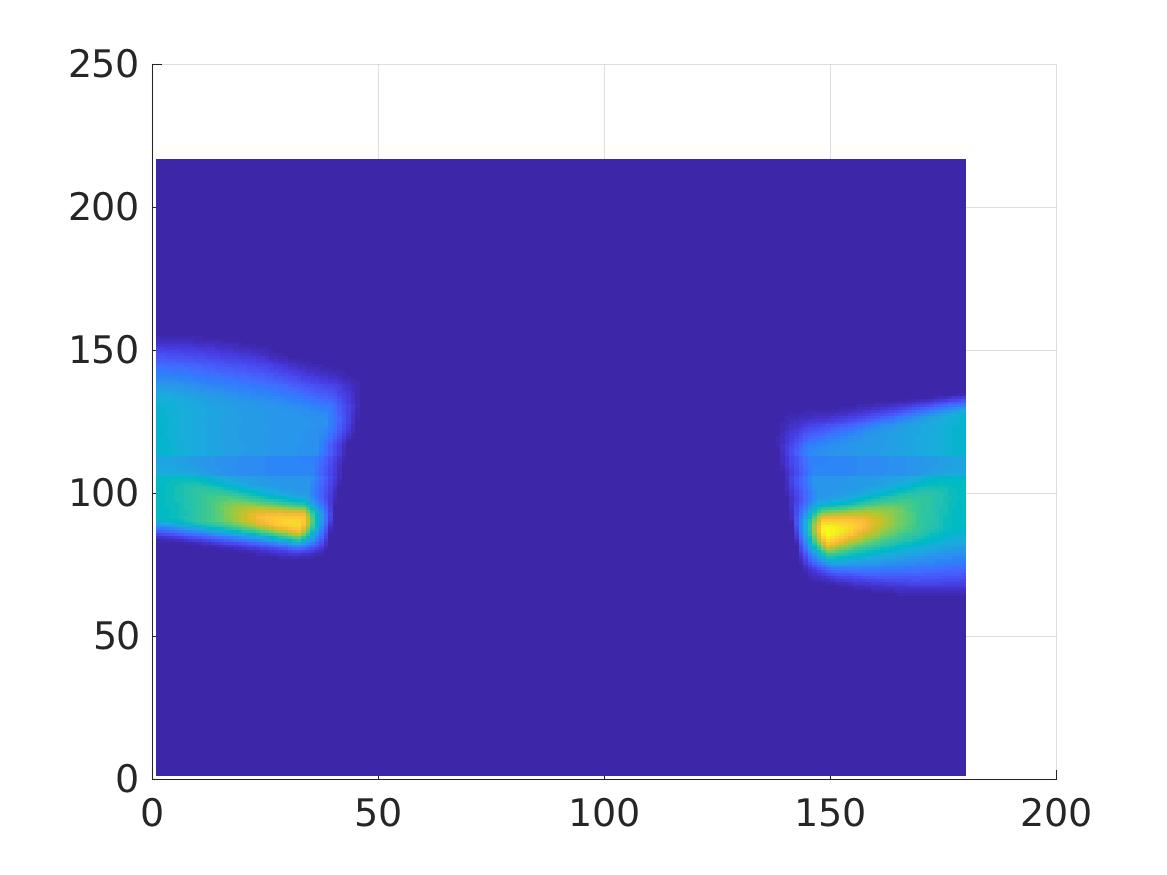}} 
		
		\subfloat[\label{fig Case 6 Radon No Clean}The function $f_{\rm comp}^{\rm iradon}$ computed by the filtered back projection algorithm, noise level $5\%$]{\includegraphics[width=.3\textwidth]{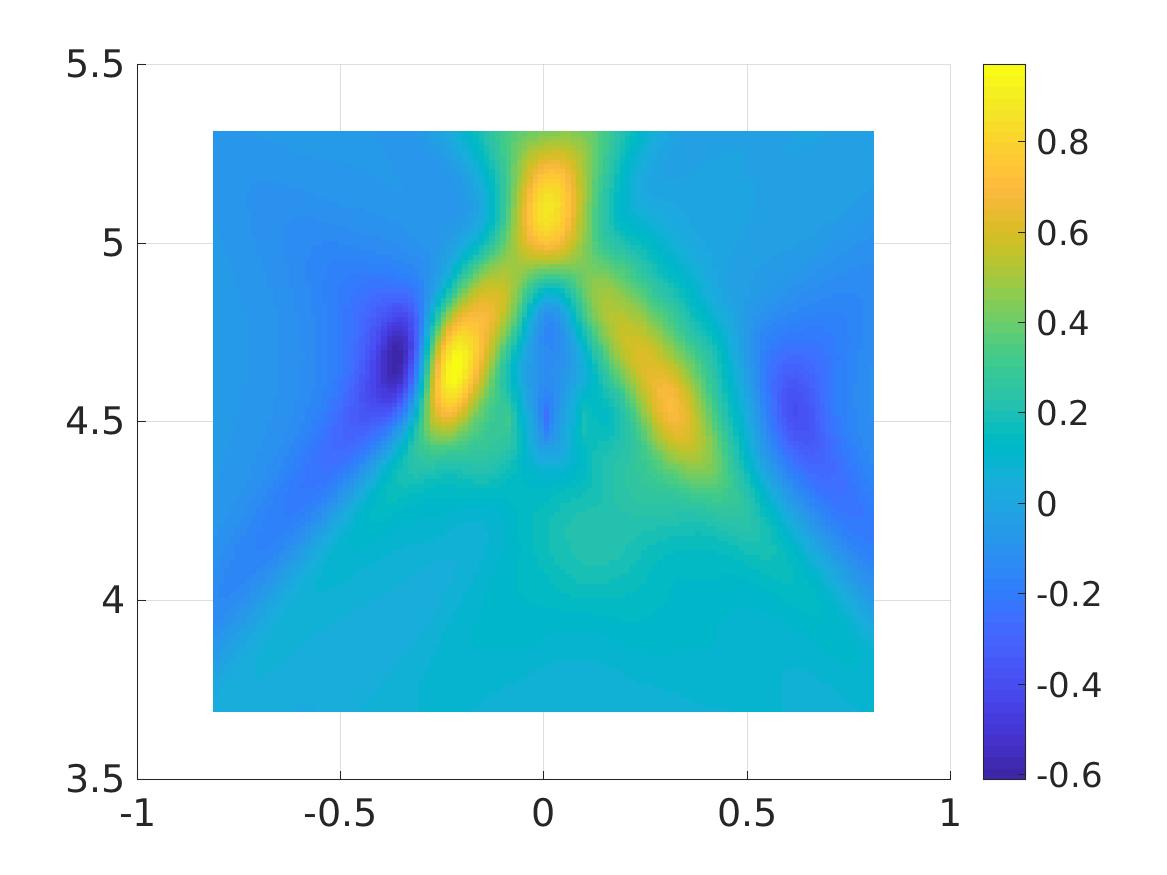}} \hspace{.2cm}
		\subfloat[\label{fig Case 6 Radon Clean}The function $f_{\rm comp}^{\rm iradon}$ computed by the filtered back projection algorithm, noise level $5\%$, together with the artifact remover]{\includegraphics[width=.3\textwidth]{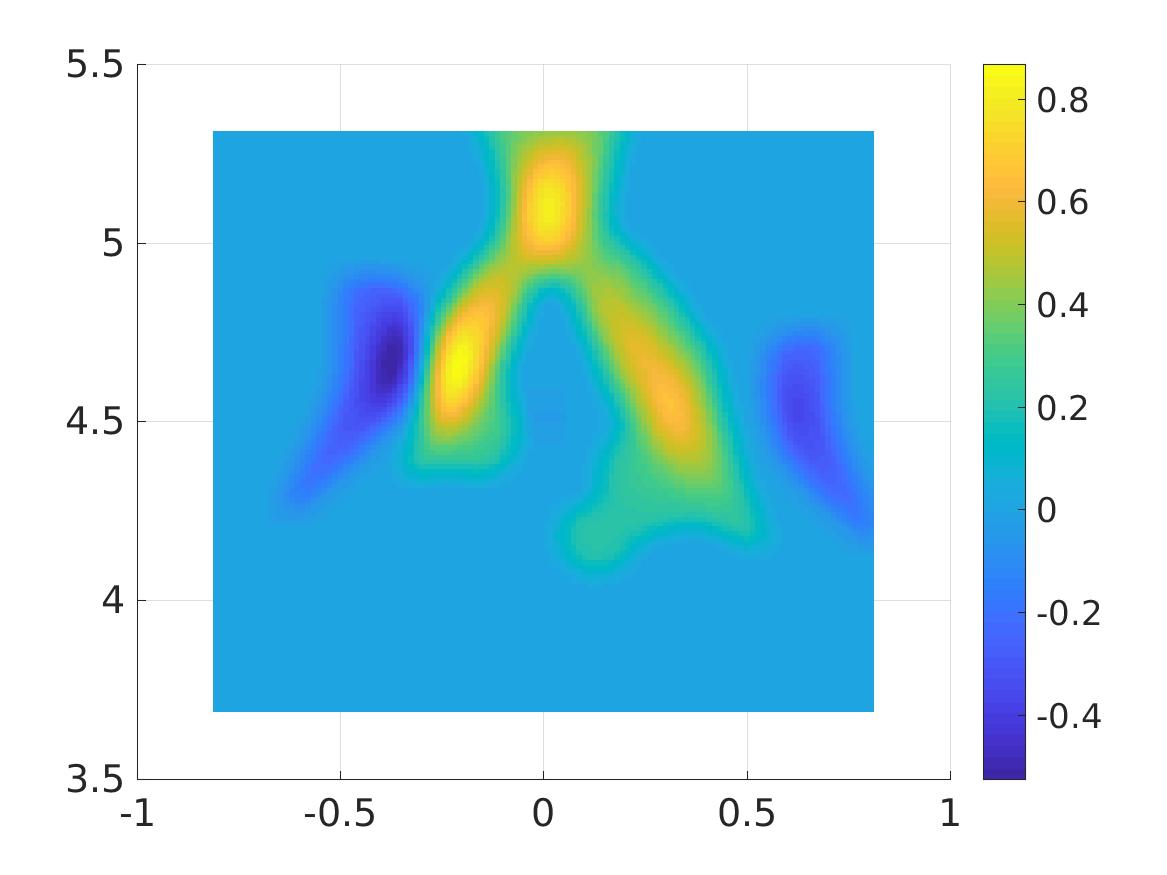}} \hspace{.2cm}
		\subfloat[\label{fig Case 6 with 5 percent noise}The function $f_{\rm comp}$ by our method in Section \ref{sec our approach}, noise level $5\%$]{\includegraphics[width=.3\textwidth]{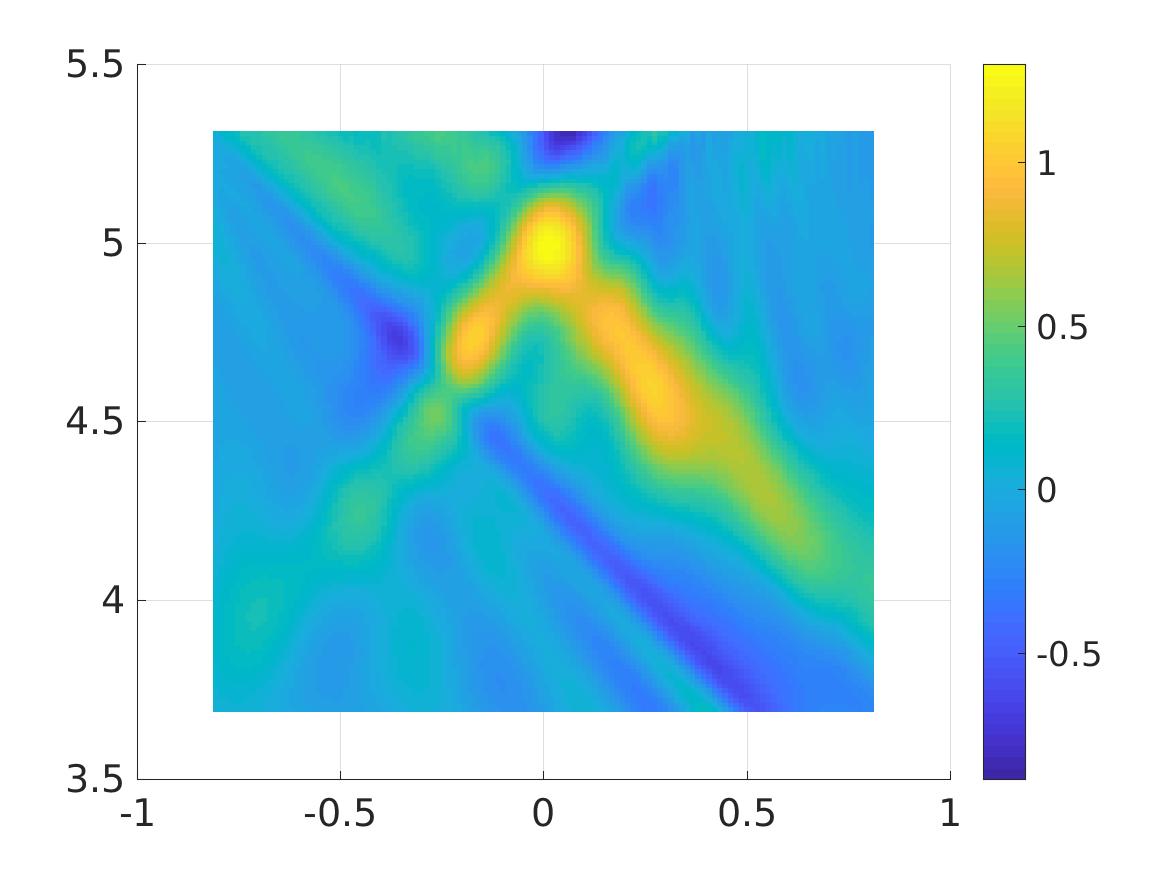}} 
		
		\subfloat[\label{fig Case 6 with 5 percent noise Clean} The function $f_{\rm comp}$ by our method in Section \ref{sec our approach}, noise level $5\%$, together with the artifact remover]{\includegraphics[width=.3\textwidth]{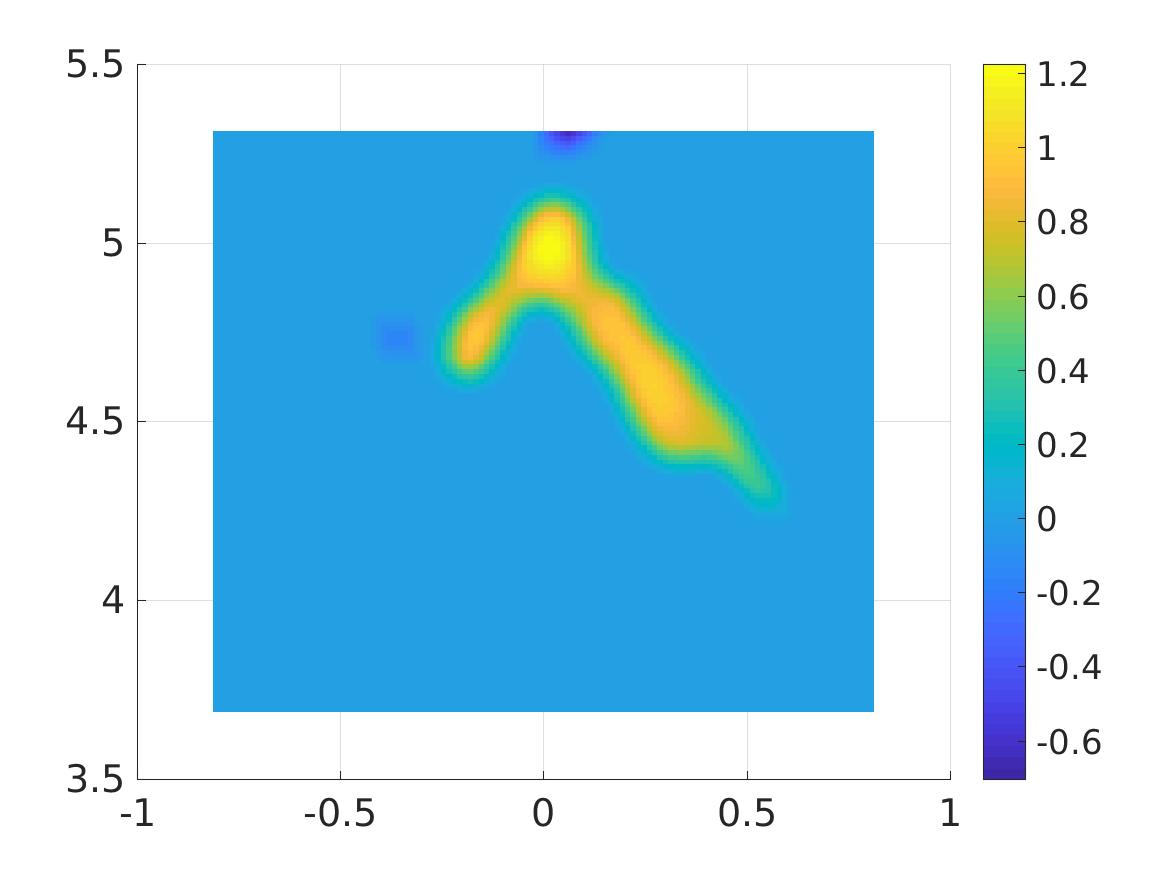}}  \hspace{.2cm}
		\subfloat[\label{fig Case 6 with 15 percent noise}The function $f_{\rm comp}$ by our method in Section \ref{sec our approach}, noise level $15\%$]{\includegraphics[width=.3\textwidth]{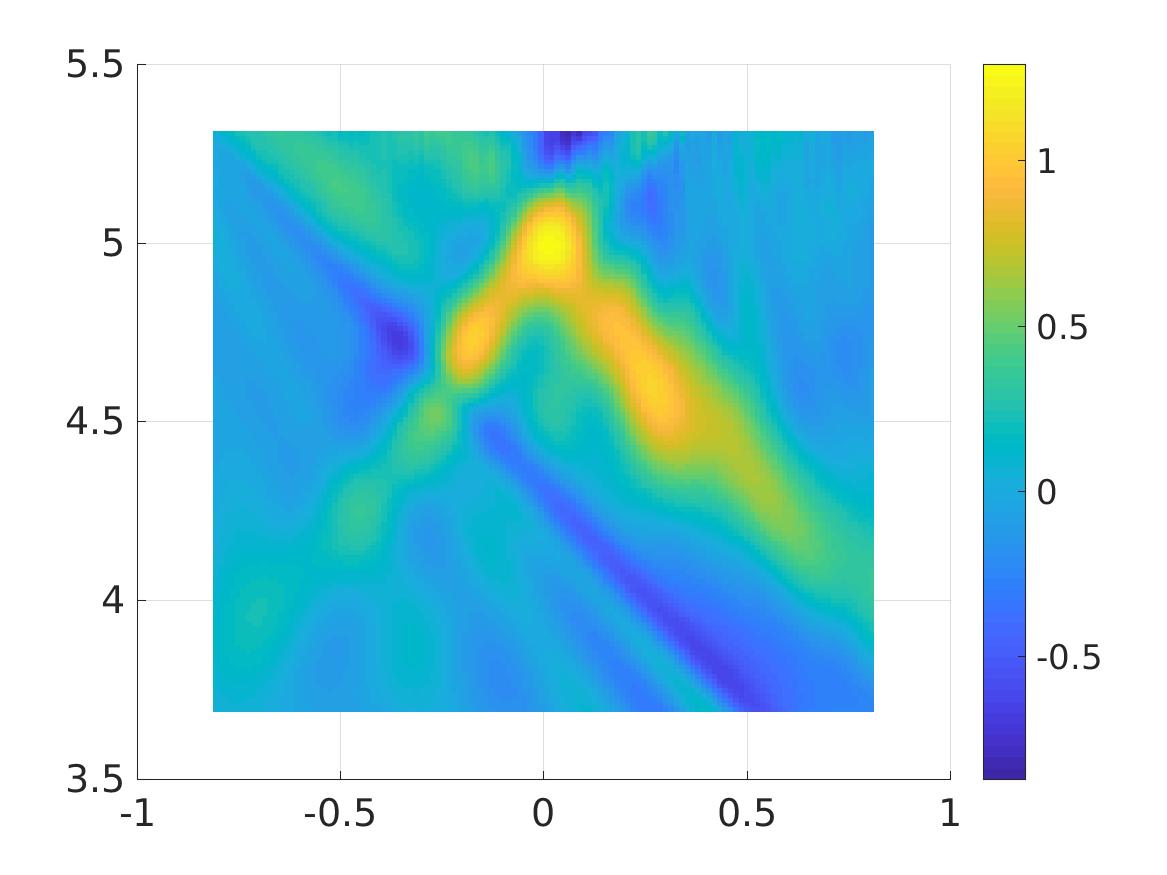}}  \hspace{.2cm}
		\subfloat[\label{fig Case 6 with 15 percent noise Clean}The function $f_{\rm comp}$ by our method in Section \ref{sec our approach}, noise level $15\%$, together with the artifact remover]{\includegraphics[width=.3\textwidth]{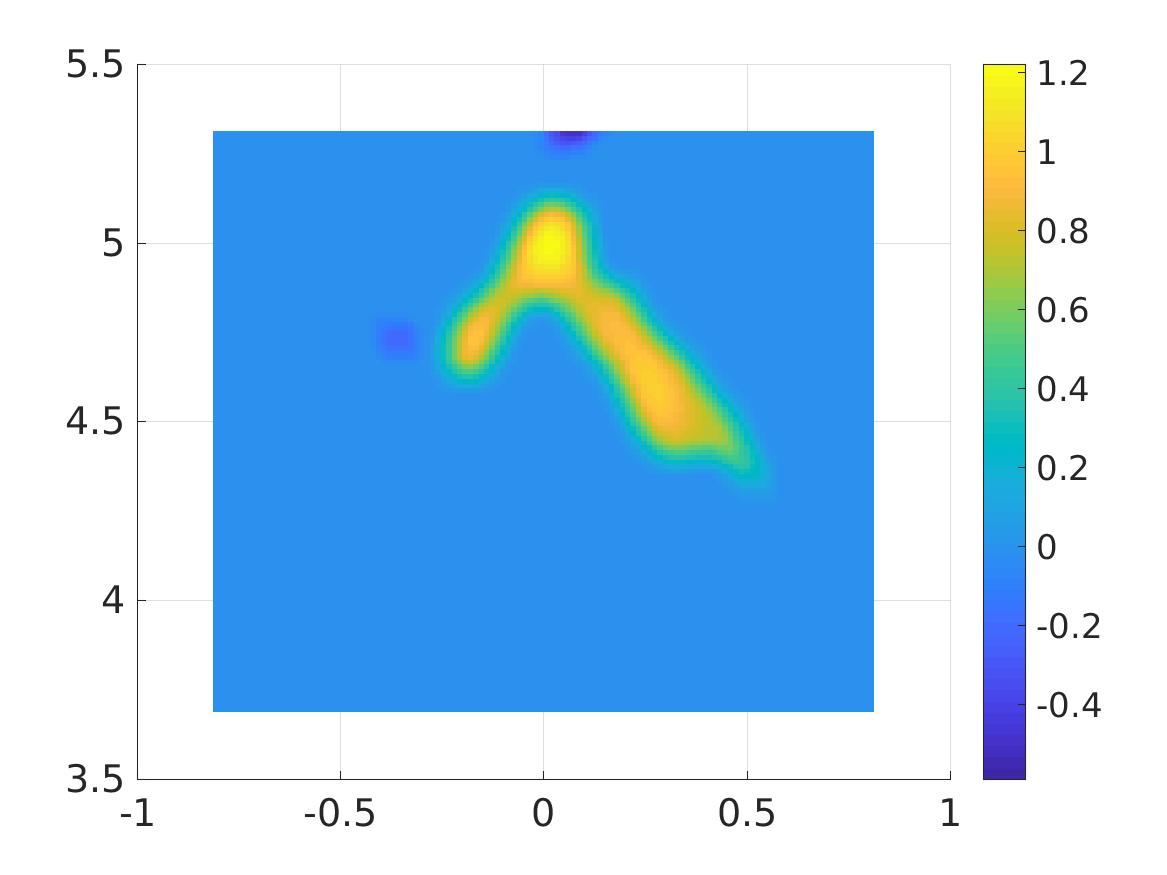}}  \hspace{.2cm}

\caption{\it  \label{fig case 6}
Test 3. The data and the reconstructions of the non smooth function $%
f^{\ast }$ in the case of an $L$-like shape. The shape is well
seen on images (g) and (i) which result from our method and it is not seem
well on (e), which results from the filtered back projection method.
Comparison of (e) with (g) and (i) indicates that the image quality provided
by our method is significantly better than that of the filtered back
projection method.}		
\end{center}
\end{figure}

\item \textbf{Test 4.} We next test our method with a non smooth function $%
f^{\ast }$ that is nonzero on a square rotated by $\pi /4$ around the center
of $\Omega .$ This square has two positive sides and two negative sides. In
particular, we want to see whether or not our method can detect a void
inside of a square.

The domain $\Omega $ is set to be $\Omega =(-1,1)\times (3.5,5.5),$ just as
in the previous numerical test. The function $f^{\ast }$ is given by 
\begin{multline*}
f^{\ast }(\mathbf{x})
=\chi _{\{\mathbf{x}=(x,y):0.3<|x|+|y-4.5|<0.6,y>4.5%
\}}
\\
-\chi _{\{\mathbf{x}=(x,y):0.3<|x|+|y-4.5|<0.6,y<4.5\}}.
\end{multline*}%
The true and reconstructed functions $f$ are displayed in Figure \ref{fig
case 7}.

\begin{figure}
	\begin{center}
		\subfloat[The true function $f^*$]{\includegraphics[width=.3\textwidth]{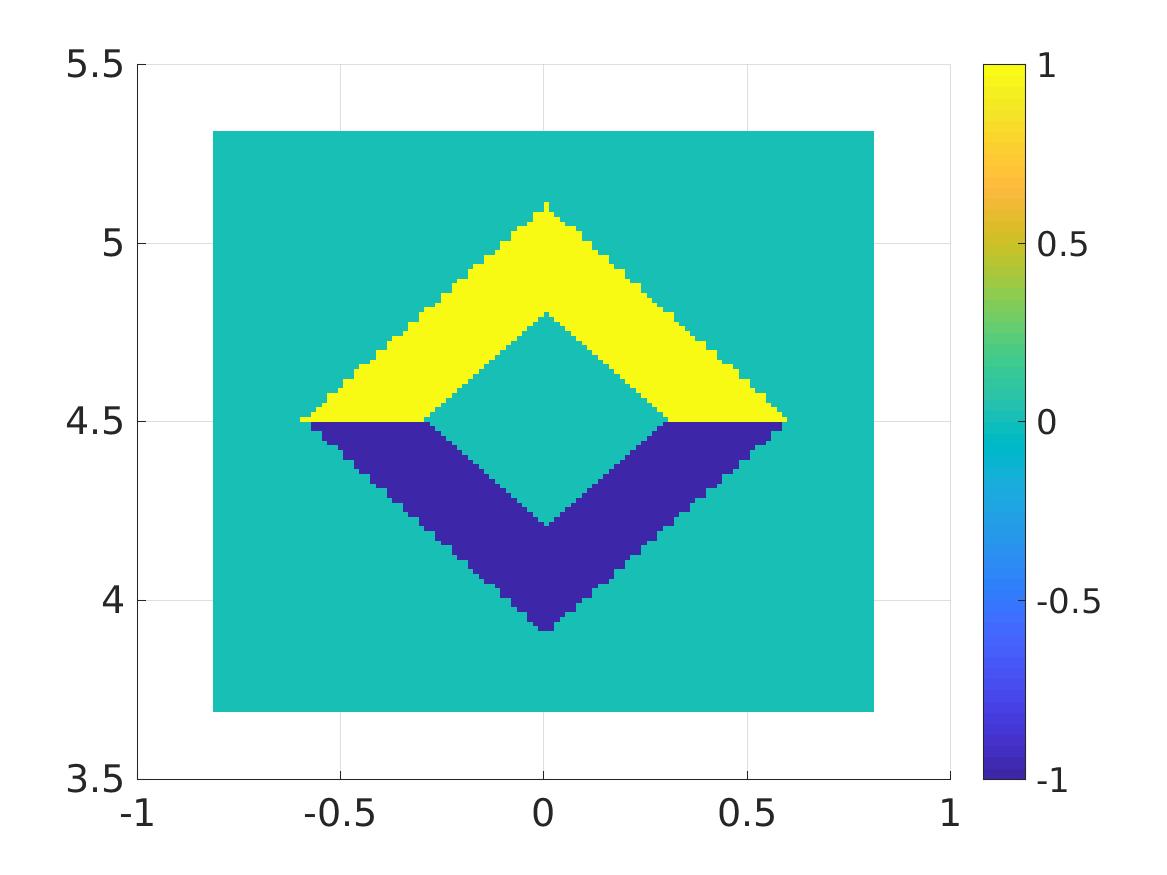}} \hspace{.2cm}
		\subfloat[\label{fig Case 7 full}The Radon transform of $f^*$ computed by the function ``radon" of Matlab]{\includegraphics[width=.3\textwidth]{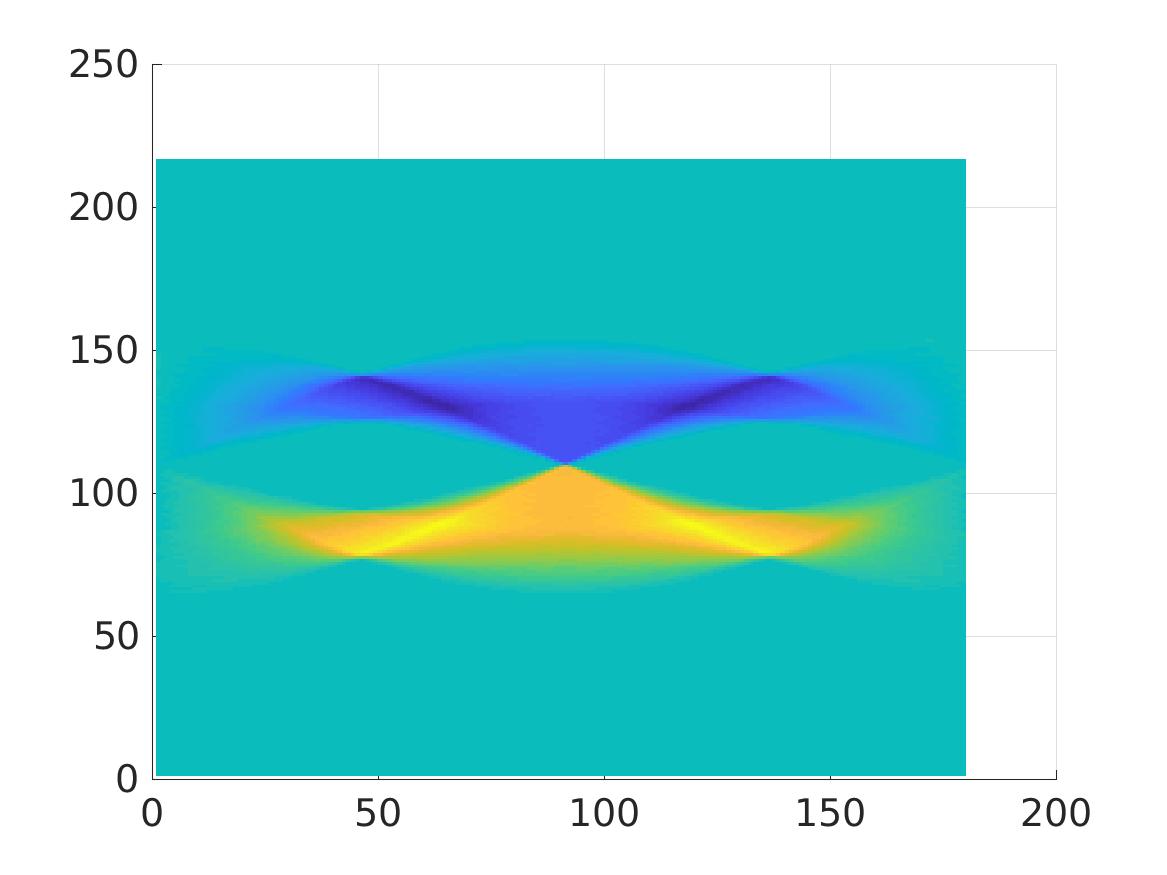}} \hspace{.2cm}
		\subfloat[\label{fig Case 7 incomplete} The incomplete tomographic data with $5\%$ noise]{\includegraphics[width=.3\textwidth]{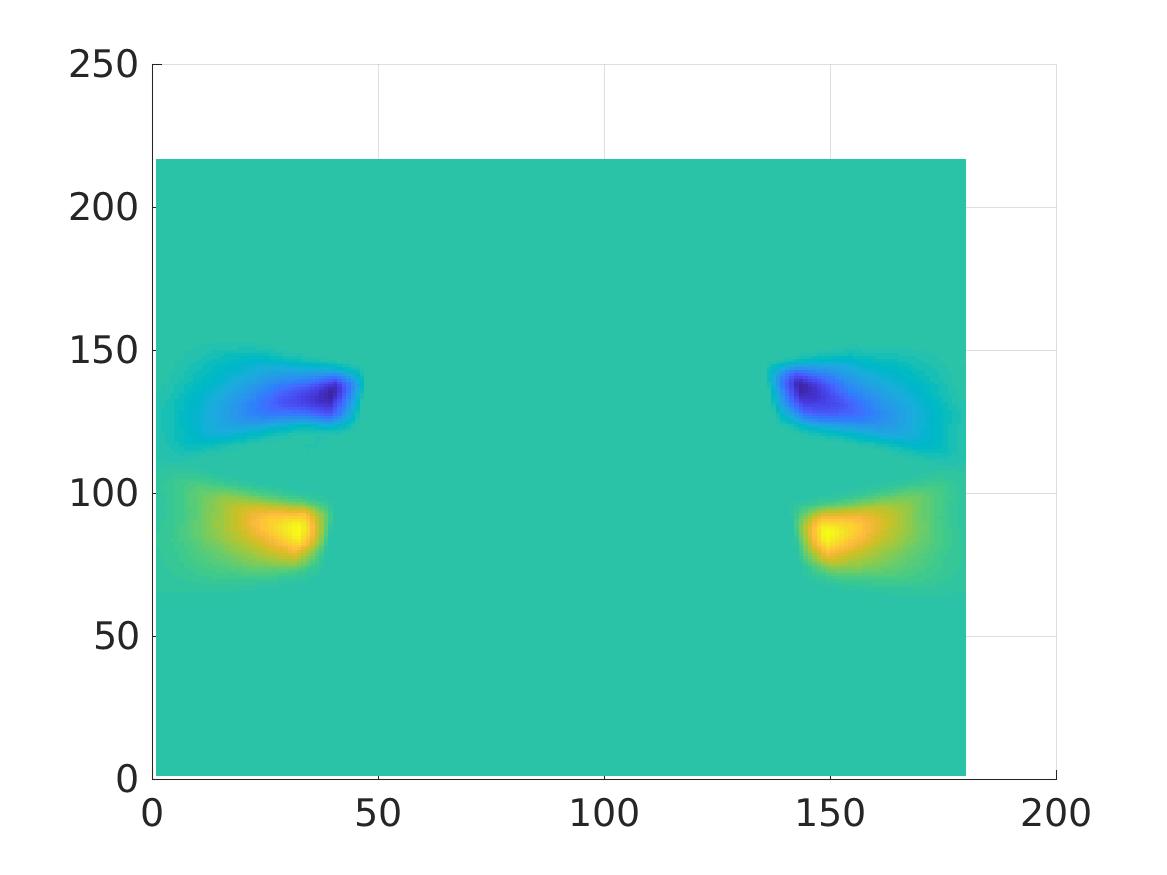}} 
		
		\subfloat[\label{fig Case 7 Radon No Clean}The function $f_{\rm comp}^{\rm iradon}$ computed by the filtered back projection algorithm, noise level $5\%$]{\includegraphics[width=.3\textwidth]{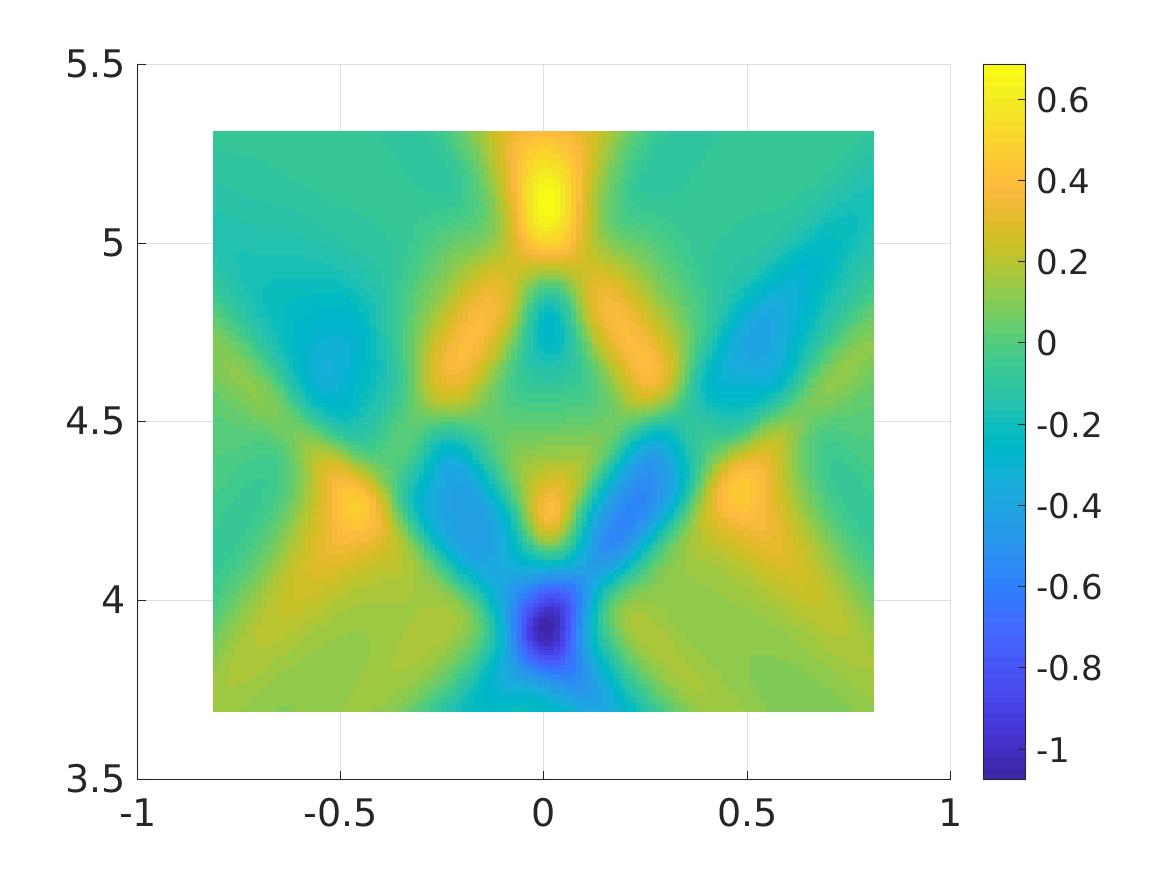}} \hspace{.2cm}
		\subfloat[\label{fig Case 7 Radon Clean}The function $f_{\rm comp}^{\rm iradon}$ computed by the filtered back projection algorithm, noise level $5\%$, together with the post processing of Section \ref{sec post}]{\includegraphics[width=.3\textwidth]{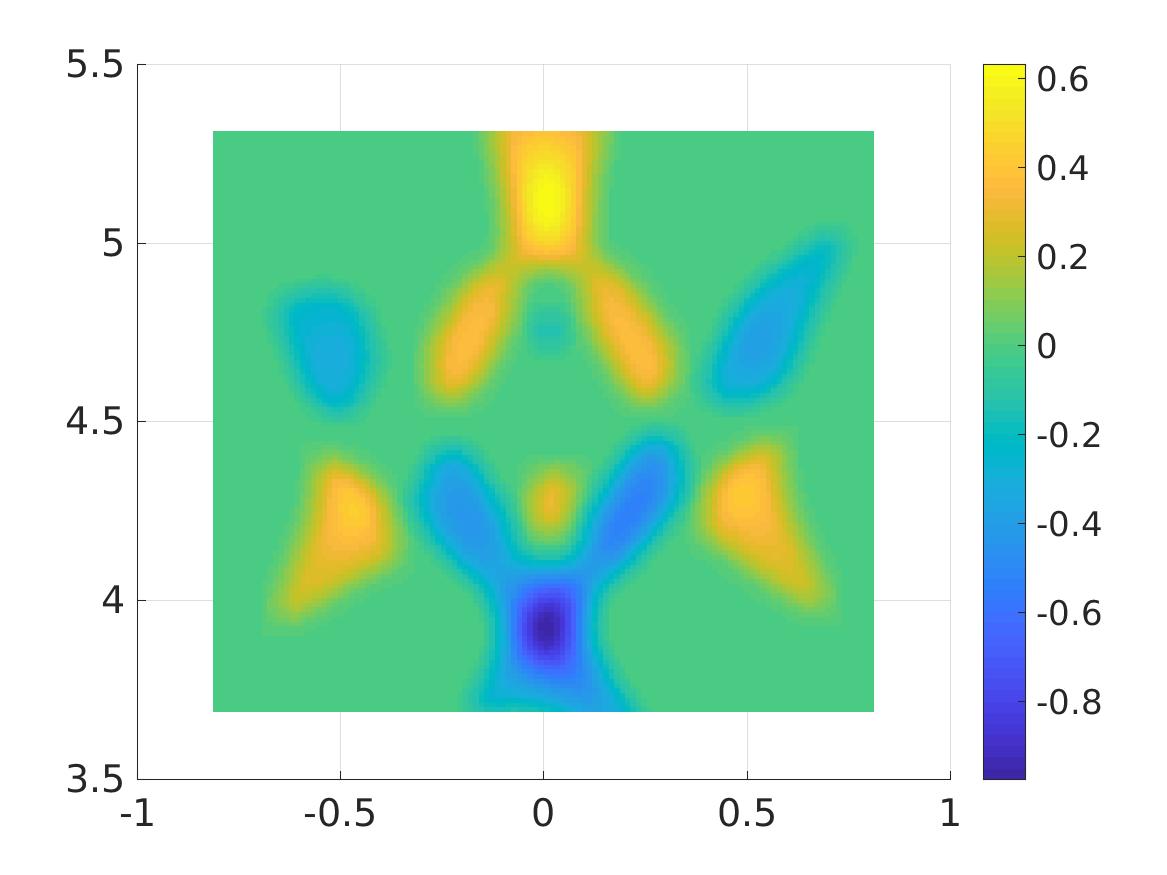}} \hspace{.2cm}
		\subfloat[\label{fig Case 7 with 5 percent noise}The function $f_{\rm comp}$ by our method in Section \ref{sec our approach}, noise level $5\%$]{\includegraphics[width=.3\textwidth]{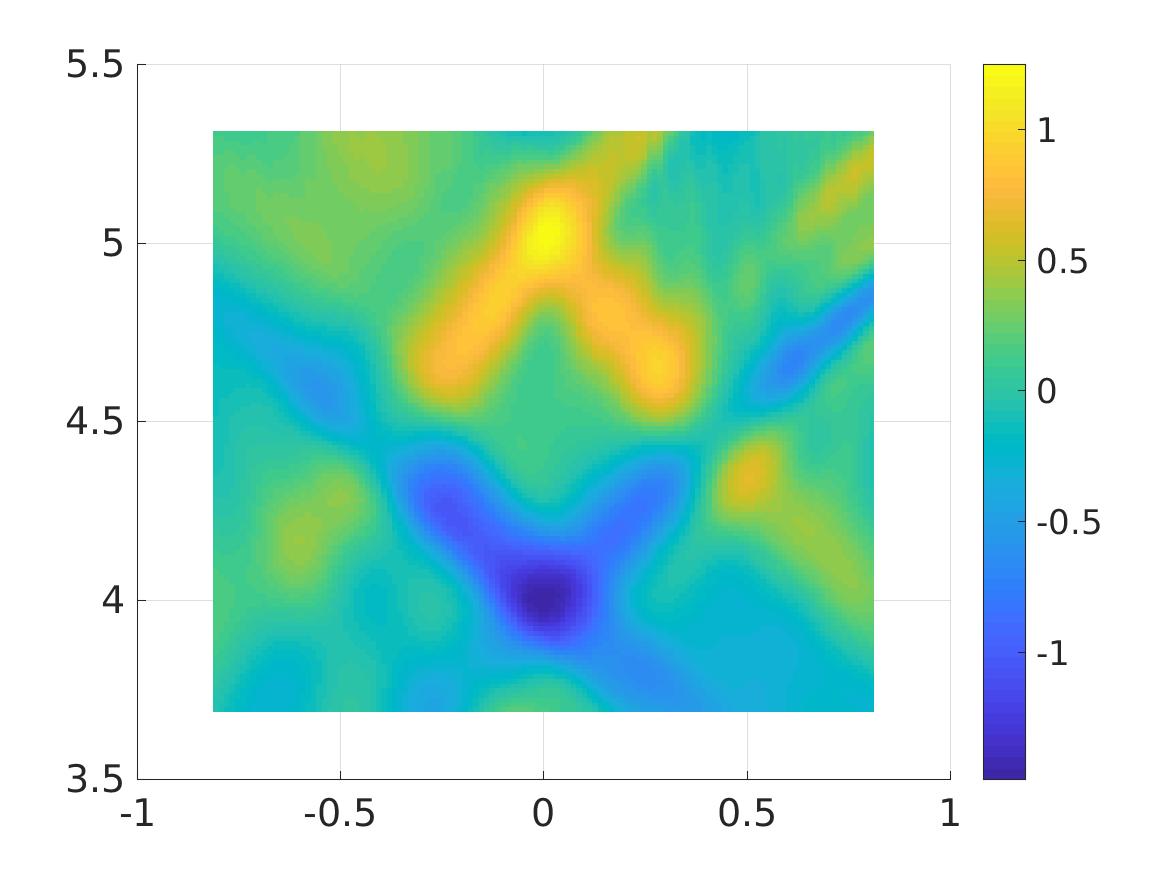}} 
		
		\subfloat[\label{fig Case 7 with 5 percent noise Clean} The function $f_{\rm comp}$ by our method in Section \ref{sec our approach}, noise level $5\%$, together withpost processing of Section \ref{sec post}]{\includegraphics[width=.3\textwidth]{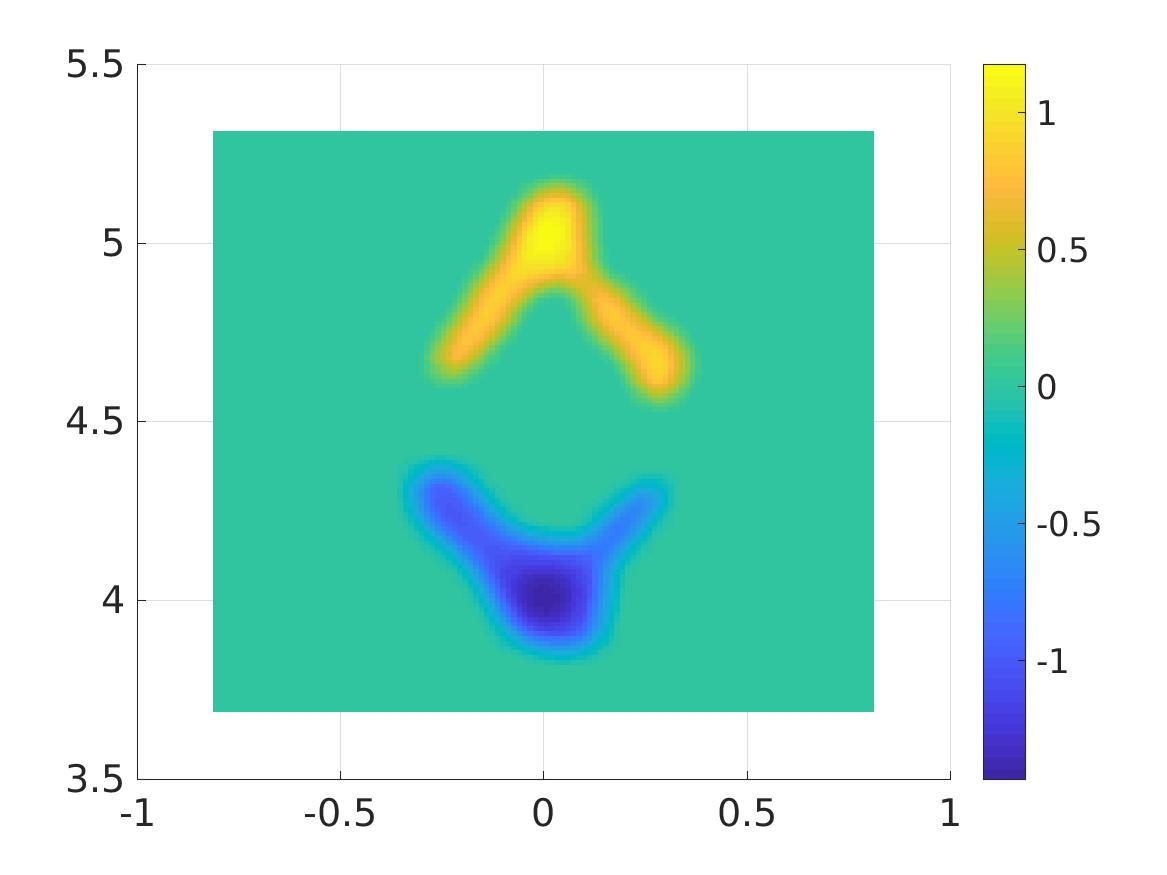}}  \hspace{.2cm}
		\subfloat[\label{fig Case 7 with 15 percent noise}The function $f_{\rm comp}$ by our method in Section \ref{sec our approach}, noise level $15\%$]{\includegraphics[width=.3\textwidth]{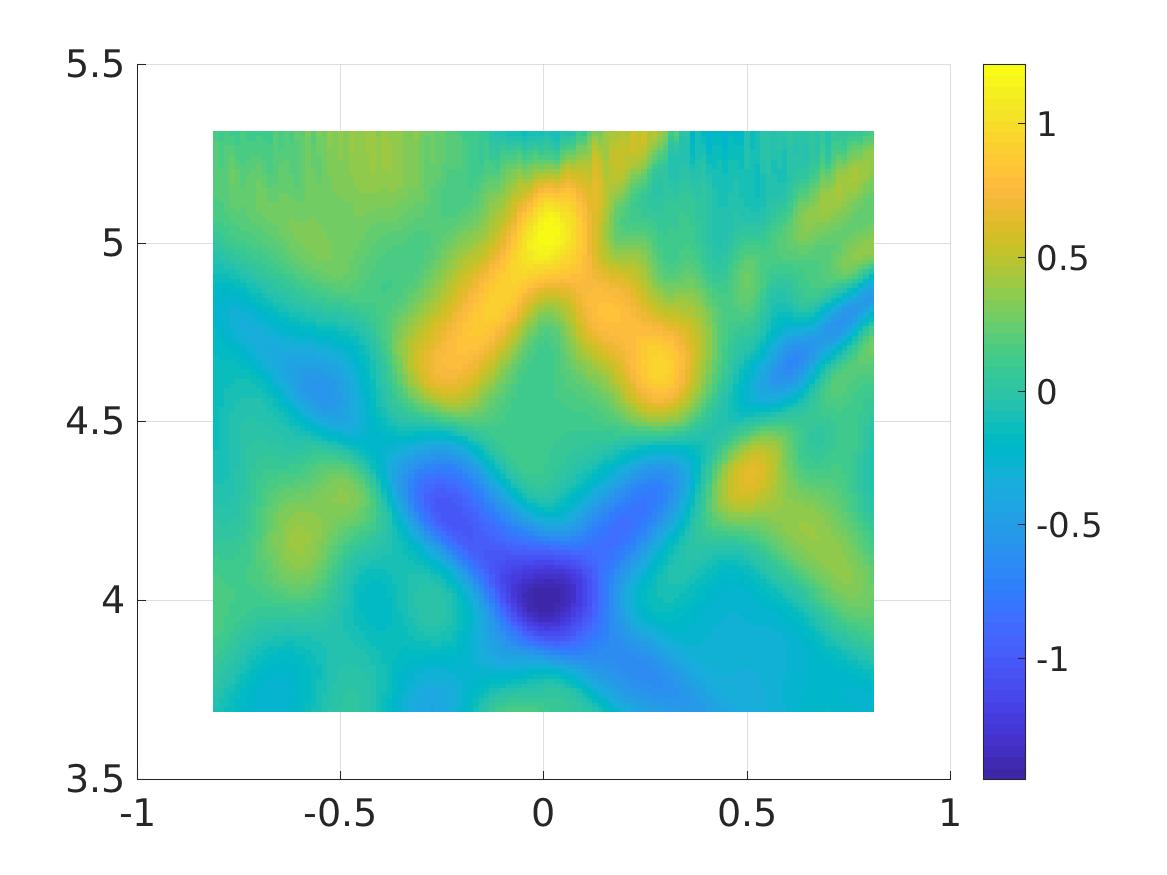}}  \hspace{.2cm}
		\subfloat[\label{fig Case 7 with 15 percent noise Clean}The function $f_{\rm comp}$ by our method in Section \ref{sec our approach}, noise level $15\%$, post processing of Section \ref{sec post}]{\includegraphics[width=.3\textwidth]{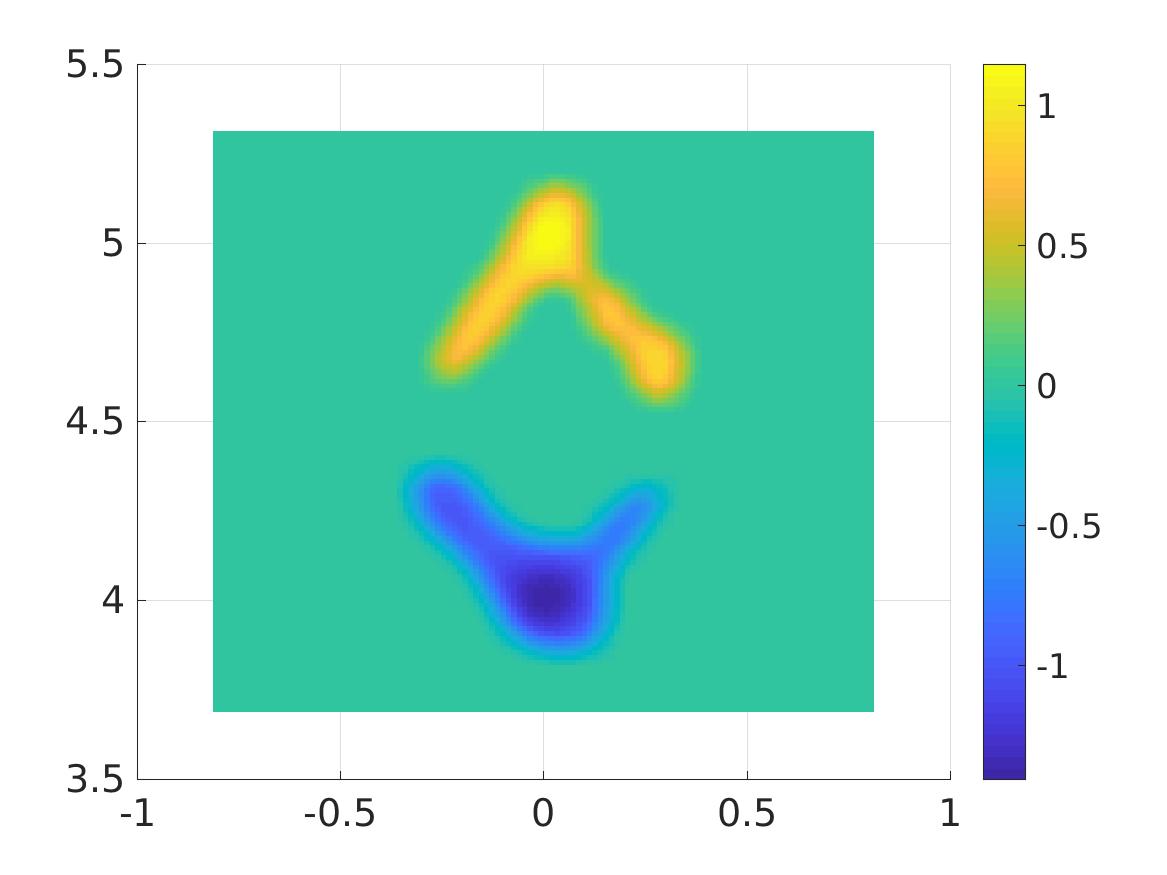}}  \hspace{.2cm}

\caption{\it  \label{fig case 7}
Test 4. The data and the reconstructions of the non smooth function $%
f^{\ast }$ in the case of a square shape. The shape is satisfactory
on (g) and (i) and the void is clearly seen on them, whereas (e) is less
clear. Comparing of (e) with (g) and (i), one can see that the image quality
provided by our method is significantly better than that of the filtered
back projection method.}		
\end{center}
\end{figure}

\end{enumerate}

One can see from these figures that our method is quite stable with respect
to the noise. In fact, the reconstructed errors and images do not change
much when the noise increases from $5\%$ to $15\%$. 

\begin{remark}[The comparison of artifacts]
Comparing Figures \ref{fig Case 1 Radon Clean}--\ref{fig Case 7 Radon Clean}
versus Figures \ref{fig Case 1 with 5 percent noise Clean}--\ref{fig Case 7
with 5 percent noise Clean} and Figures \ref{fig Case 1 with 15 percent
noise Clean}--\ref{fig Case 7 with 15 percent noise Clean}, we observe that
the unwanted artifacts involved in the results by the filtered back
projection method are much stronger than the ones in the numerical
reconstructions obtained by our method. More precisely, the
\textquotedblleft $20\%$ filter" in \eqref{eqn filter} cannot remove
unwanted artifacts in $f_{\mathrm{comp}}^{\mathrm{iradon}}$ while it works
well for the artifacts in $f_{\mathrm{comp}}$ obtained by using our method.
\end{remark}

It seems to be on the first glance that the longer the source line $\Gamma
_{d}$ in (\ref{2.10}) is, the wider is the angle to \textquotedblleft see"
the inclusions. However, in computation, there is a limiting length $%
2d_{\lim }$ for $\Gamma _{d}$ such that our method fails for $d>d_{\lim }.$
For example, for parameters $a$ and $b$ used in Tests 1 and 2, this limiting
length is $d_{\lim }=14$. To explain this length limitation, we observe that
a more detailed analysis of formulae (\ref{4.4}), (\ref{4.6}) and Lemma \ref%
{lem 3.2} shows that one should have in Lemma \ref{lem 3.3} $\left(
R+d\right) /a_{0}^{2}\ll 1.$ The fact that this inequality is not exactly
satisfied in Tests 1-4 can be viewed as another indication of the stability
of our technique. However, this inequality is violated at large for $d\geq
14,$ and this is why our method fails to work for such values of $d$.

\subsection{Reconstruction errors}

\label{sec errors}

As to the image quality, the visual analysis of Figures
2(e),(g),(i)-5(e),(g),(i) indicates that, at least in our four tests, our method provides better quality images than the filtered back
projection method. Furthermore, the difference of those qualities increases
in the favor of our method as the structures of inclusions become more
complicated.

\begin{figure}[tbp]
\begin{center}
\subfloat[Absolute errors in locations of reconstructed inclusions. In spite
of noise, the reconstructed locations of inclusions 2 and 4 are
exact.]{\includegraphics[width=.47\textwidth]{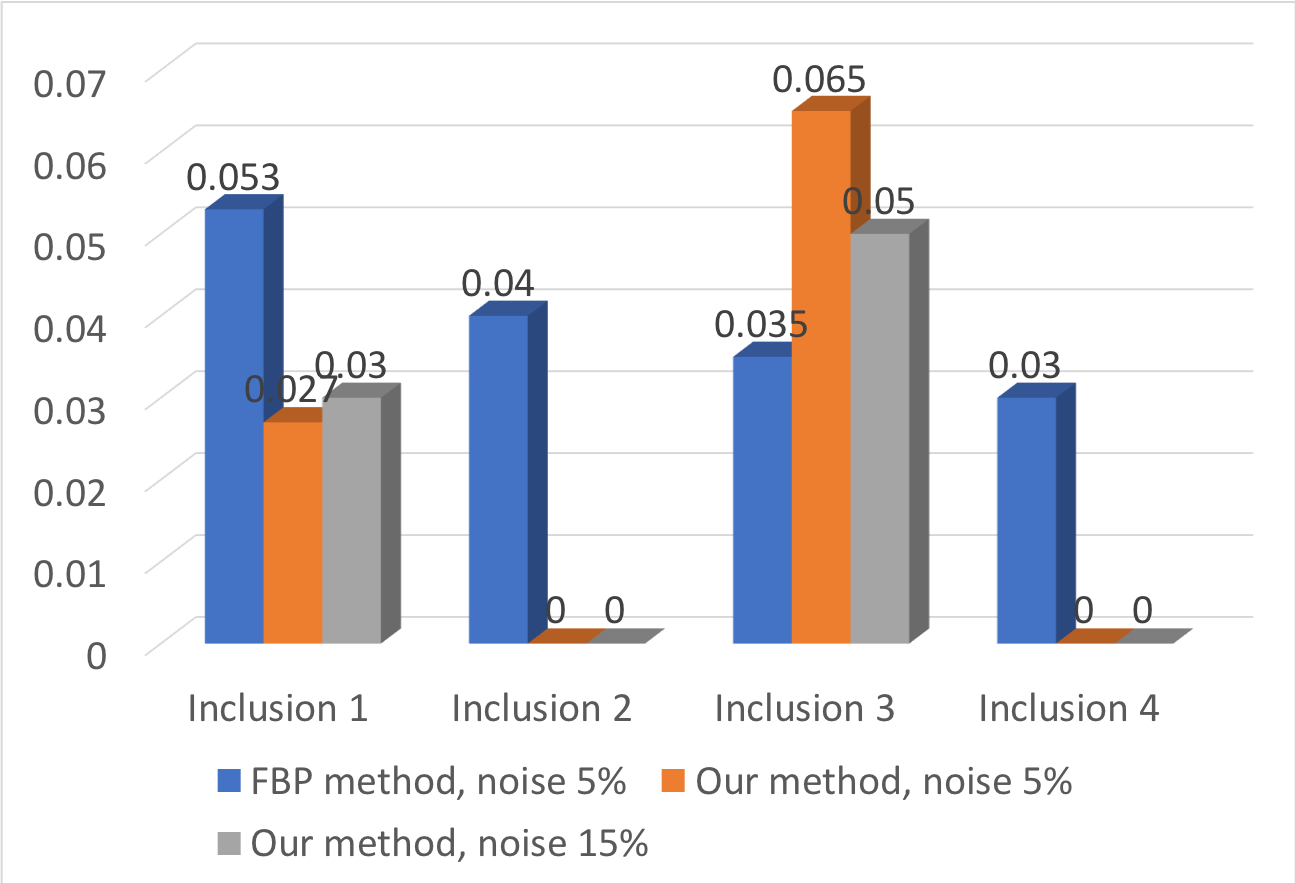}}  
\hspace{.2cm}  \subfloat[Relative errors (in \%) of the extreme value of the
reconstructed function $f_{\rm comp}$ in four inclusions in Tests 1 and
2.]{\includegraphics[width=.47\textwidth]{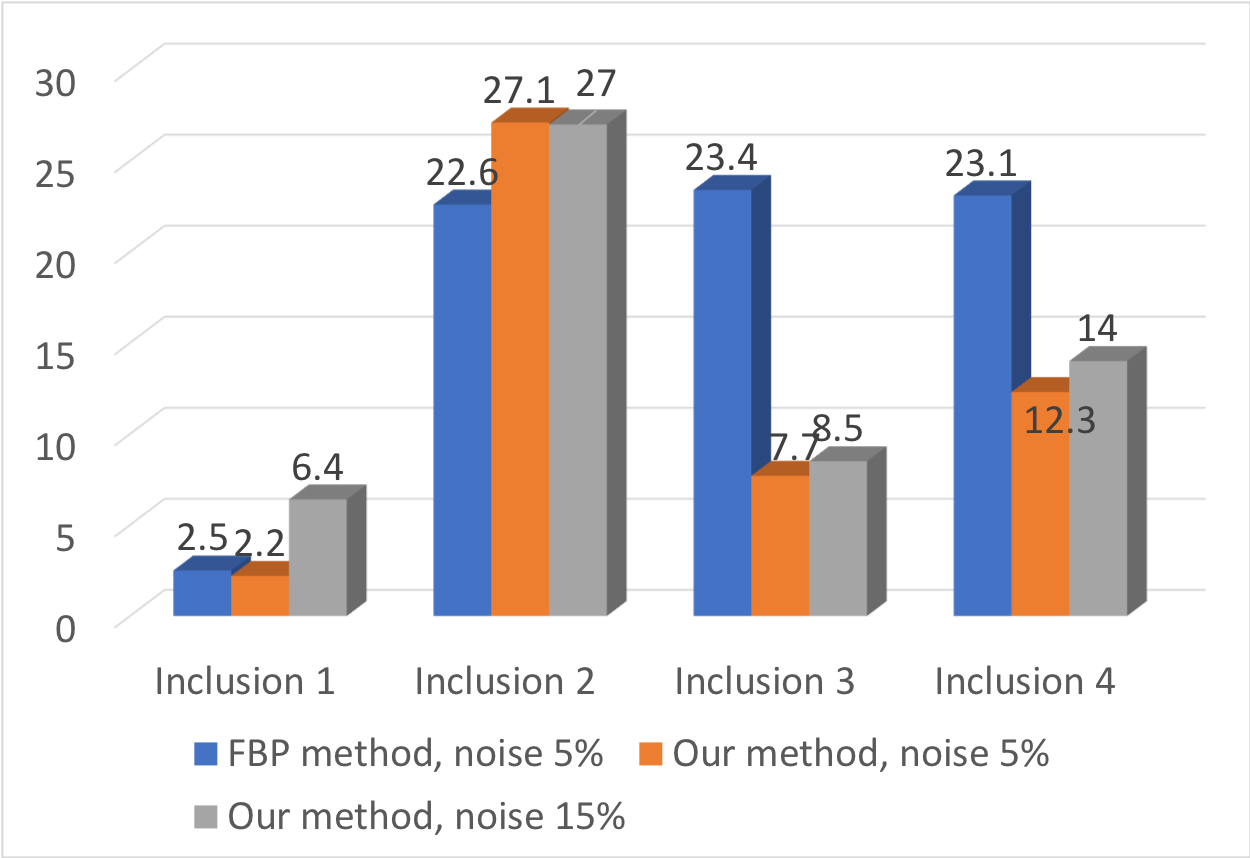}} 
\end{center}
\caption{\textit{Comparison of reconstruction errors of our method with
filtered back projection method for Tests 1 and 2, see their descriptions as
well as Figures \protect\ref{fig case 1} and \protect\ref{fig case 3} for
numbering of inclusions. a) Absolute errors in locations of points with
extreme values of the function $f$. b) Relative errors (in \%) of the
reconstructed extreme values of the function f inside the inclusions.}}
\label{fig error}
\end{figure}

We now discuss the reconstructed errors of the numerical solutions only in
the first two tests in which true function $f^{\ast }$ involves inclusions.
Satisfactory reconstructed values were obtained, see Table \ref{table 1} and
Figure \ref{fig error}. We do not present the error estimates for Tests 3
and 4 since it is not clear for us how to define the values of the
reconstructed functions for the kinds of non-convex inclusions in those two
tests. However, it can be seen from Figures \ref{fig Case 6 Radon Clean}, %
\ref{fig Case 6 with 5 percent noise Clean}, \ref{fig Case 6 with 15 percent
noise Clean}, \ref{fig Case 7 Radon Clean}, \ref{fig Case 7 with 5 percent
noise Clean} and \ref{fig Case 7 with 15 percent noise Clean} and the
enclosed color bars that the reconstructions of the function $f$ are
acceptable.

\begin{table}[tbp]
\caption{\textit{Correct and computed inclusions in Tests 1 and 2. Here, FBP
means filtered back projection, Nm means inclusion number (see descriptions
of Tests 1,2), $f_{\mathrm{true}}$ is the extreme value of the true function 
$f^*$ in the inclusion, loc$_{\mathrm{true}}$ means true location where the
extreme value of $f^*$ occurs, $f_{\mathrm{comp}}$ is the extreme value of
the computed function $f_{\mathrm{comp}}$ in the inclusion, and loc$_{%
\mathrm{true}}$ means true location where the extreme value of $f_{\mathrm{%
comp}}$ occurs. }}
\label{table 1}
\begin{center}
{\small
\begin{tabular}{|c|c|c|c|c|c|c|}
\hline
Inc. Nm & loc$_{\mathrm{true}}$ & $f_{\mathrm{true}}$ & Method & noise level
& loc$_{\mathrm{comp}}$ & $f_{\mathrm{comp}}$ \\ \hline
1 & (0.0, 2) & 1 & FBP method & 5\% & (0.053, 2.000) & 0.9751 \\ \hline
1 & (0.0, 2) & 1 & Our method & 5\% & (0.000, 1.973) & 0.9781 \\ \hline
1 & (0.0, 2) & 1 & Our method & 15\% & (0.013, 1.973) & 0.9361 \\ \hline
2 & (-0.4, 4) & -6 & FBP method & 5\% & (-0.400, 3.960) & -4.644, \\ \hline
2 & (-0.4, 4) & -6 & Our method & 5\% & (-0.4, 4) & -4.373 \\ \hline
2 & (-0.4, 4) & -6 & Our method & 15\% & (-0.4, 4) & -4.378 \\ \hline
3 & (-0.1, 3.5714) & 5 & FBP method & 5\% & (-0.067, 3.560) & 3.829 \\ \hline
3 & (-0.1, 3.5714) & 5 & Our method & 5\% & (-0.107, 3.507), & 4.615 \\ 
\hline
3 & (-0.1, 3.5714) & 5 & Our method & 15\% & (-0.107, 3.52) & 4.574 \\ \hline
4 & (0.4, 4) & 6 & FBP method & 5\% & (0.413, 4.027) & 4.617 \\ \hline
4 & (0.4, 4) & 6 & Our method & 5\% & (0.4, 4) & 5.261 \\ \hline
4 & (0.4, 4) & 6 & Our method & 15\% & (0.4, 4) & 5.16 \\ \hline
\end{tabular}
}
\end{center}
\end{table}

We now analyze Figure \ref{fig error}. In this figure \textquotedblleft
absolute errors in locations of reconstructed inclusions" means errors at
points where the reconstructed function $f$ achieves its extreme value, for
each inclusion of Tests 1,2. One can see from Figure \ref{fig error}(a)
that, in terms of locations, our method performs better than the filtered
back projection method for inclusions 1,2 and 4. And it performs worse for
inclusion number 3. As to Figure \ref{fig error}(b), one can observe that
our method provides more accurate extreme values for inclusions 3 and 4. In
the case of inclusion 1, the accuracy in calculating extreme values is about
the same for both methods for the case of 5\% noise. In the case of
inclusion 2, the accuracy in calculating the extreme value is better for
filtered back projection method.

\section{Concluding Remarks}

\label{sec concl}

While all current techniques of the inversion of the data for the X-ray
tomography are based on some inversion formulae, we have proposed a new
numerical method here, which does not intend to obtain an inversion formula.
Instead, it uses a well known transport PDE governing propagation of X-rays.
Our method works for a special case of a limited angle data, which might be
potentially applied to, e.g. checking out bulky luggage in airports and
checking out quality of walls in houses. Using the original idea of the
method of \cite{BukhKlib} as well as a recently introduced new orthonormal
basis in $L^{2}\left( -d,d\right) $ \cite{Klibanov:jiip2017}, we obtain a
system of coupled first order PDEs in which the target function $f$ is not
involved. The boundary value problem for this system is over determined.
Therefore, we solve this boundary value problem by the quasi-reversibility
method, which is perfectly suited for solutions of overdetermined boundary
value problems for PDEs. We prove a new Carleman estimate and use it then to
prove uniqueness and existence of the solution for the quasi-reversibility
method. Next, the same Carleman estimate enables us to establish convergence
rate of regularized solutions. We work with a semi discrete version of the
quasi-reversibility method, which is more realistic than its conventional
continuos version, see a survey in \cite{Ksurvey} for the continuos version. 

We have conducted numerical testing of this new method for noisy data,
including comparison with the filtered back projection method for Radon
transform. In the latter we have heuristically assigned zero to the missing
data, similarly with \cite{Borg}. We point out that this assignment cannot
be rigorously justified for the filtered back projection method, unlike our
method. We have observed that our method sustains 5\% and 15\% of noise and
resulting images are about the same.

The visual analysis indicates that, at least in the above Tests 1-4, images
resulting from our method have a better quality than those provided by the
filtered back projection method. Also, the more complicated the structures
of inclusions are, the larger in the favor of our method is the difference of qualities of those
images. It can be seen from Table \ref{table 1} that our method also
provides more accurate locations of imaged targets for Tests 1,2 for three
(3) out of four (4) inclusions. As to the extreme values of the function $f$
inside of inclusions, it can be seen from Table \ref{table 1} that our
method provides about the same accuracy as the filtered back projection
method for one inclusion (number 1), better accuracy for two (number 3,4)
and worse accuracy for one inclusion (number 2), also see Figure \ref{fig
error}. Comparison in numbers for Tests 3,4 is hard to provide due to the
complicated structures of inclusions in these tests.

Finally, we observe that, in the case of the attenuated X-ray transform \cite%
{Natterer:cmsiam2001}, the following analog of PDE (\ref{main eqn}) is valid 
\cite{HasanogluRomanov:s2017}: 
\begin{equation}
\frac{x-\alpha }{\sqrt{\left( x-\alpha \right) ^{2}+y^{2}}}u_{x}+\frac{y}{%
\sqrt{\left( x-\alpha \right) ^{2}+y^{2}}}u_{y}+c\left( x,y\right) u=f\left(
x,y\right) ,\quad (x,y)\in \Omega  \label{7.2}
\end{equation}%
with an appropriate function $c\left( x,y\right) $. This equation differs
from equation (\ref{main eqn}) by the term $c\left( x,y\right) u.$ Since
this is the lower order term in PDE (\ref{7.2}) and since Carleman estimates
are \textquotedblleft sensitive" only to the principal parts of PDE
operators and \textquotedblleft non sensitive" to their lower terms, then a
slight modification of our technique works for this case. Numerical studies
of this problem are outside of the scope of the current publication. We
refer to \cite{Nov} for an inversion formula for the attenuated X-ray
transform.


\begin{thebibliography}{99}
\bibitem{Linh:ip2016} \textsc{L.~L. Barannyk, J.~Frikel and L.~V. Nguyen}, 
\emph{On artifacts in limited data spherical radon transform: curved
observation surface}, Inverse Problems, 32 (2016), p.~015012.

\bibitem{Becacheelal:AIMS2015} \textsc{E.~B\'{e}cache, L.~Bourgeois,
L.~Franceschini, and J.~Dard\'{e}}, \emph{\ Application of mixed
formulations of quasi-reversibility to solve ill-posed problems for heat and
wave equations: The 1d case}, Inverse Problems \& Imaging, 9 (2015),
pp.~971--1002.

\bibitem{BK} \textsc{L. Beilina and M.V. Klibanov}, \emph{Approximate Global
Convergence and Adaptivity for Coefficient Inverse Problems}, Springer, New
York, 2012.

\bibitem{BY} \textsc{M. Bellassoued and M. Yamamoto}, \emph{Carleman
Estimates and Applications to Inverse Problems for Hyperbolic Systems},
Springer, Japan, 2017.

\bibitem{Borg} \textsc{L.~Borg, J.~Frikel, J.~S. J$\varnothing $rgensen, and
E.~T. Quinto}, \emph{\ Full characterization of reconstruction artifacts
from arbitrary incomplete x-ray ct data}, Arxiv:1707.03055v3, (2018).

\bibitem{Borg:Meas2017} \textsc{L.~Borg, J.~S. J$\varnothing $rgensen,
J.~Frikel, and J.~Sporring}, \emph{\ Reduction of variable-truncation
artifacts from beam occlusion during in situ {X}-ray tomography}, Meas. Sci.
Tech., 28 (2017), p.~19pp.

\bibitem{Borg:report2017} \textsc{L.~Borg, J.~S. J$\varnothing $rgensen, and
J.~Sporring}, \emph{Towards characterizing and reducing artifacts caused by
varying projection truncation}, tech. rep., Department of Computer Science,
University of Copenhagen, 2017/1.

\bibitem{Bourgeois:ip2006} \textsc{L.~Borg, J.~S. J$\varnothing $rgensen, and
J.~Sporring}, \emph{Convergence rates for the quasi-reversibility method to solve
the {C}auchy problem for {L}aplace's equation}, Inverse Problems, 22 (2006),
pp.~413--430.

\bibitem{BourgeoisDarde:ip2010} \textsc{L.~Bourgeois and J.~Dard\'{e}}, 
\emph{A duality-based method of quasi-reversibility to solve the {C}auchy
problem in the presence of noisy data}, Inverse Problems, 26 (2010),
p.~095016.

\bibitem{BukhKlib} \textsc{A.~Bukhgeim and M.~Klibanov}, \emph{Uniqueness in
the large of a class of multidimensional inverse problems}, Soviet Math.
Doklady, 17 (1981), pp.~244-247.

\bibitem{ClasonKlibanov:sjsc2007} \textsc{C.~Clason and M.~V. Klibanov}, 
\emph{The quasi-reversibility method for thermoacoustic tomography in a
heterogeneous medium}, SIAM J. Sci. Comput., 30 (2007), pp.~1--23.

\bibitem{Dadre:ipi2016} \textsc{J.~Dard\'e}, \emph{Iterated
quasi-reversibility method applied to elliptic and parabolic data completion
problems}, Inverse Problems and Imaging, 10 (2016), pp.~379--407.

\bibitem{Frikel:ip2013} \textsc{J.~Frikel and E.~T. Quinto}, \emph{%
Characterization and reduction of artifacts in limited angle tomography},
Inverse Problems, 29 (2013), p.~125007.

\bibitem{HasanogluRomanov:s2017} \textsc{A.~H. Hasano\u{g}lu and V.~G.
Romanov}, \emph{Introduction to Inverse Problems for Differential Equations}%
, Springer, Cham, 2017.

\bibitem{Kab1} \textsc{S.I. Kabanikhin, A. D. Satybaev and M.A. Shishlenin}, 
\emph{Direct Methods of Solving Inverse Hyperbolic Problems}, VSP, The
Netherlands, 2005.

\bibitem{Kab2} \textsc{S.I. Kabanikhin and M.A. Shishlenin}, \emph{Numerical
algorithm for two-dimensional inverse acoustic problem based on
Gel'fand--Levitan--Krein equation}, J. Inverse and Ill-Posed Problems, 18
(2011), 979-995.

\bibitem{Kab3} \textsc{S.I. Kabanikhin, K.K. Sabelfeld, N.S. Novikov and
M.A. Shishlenin, }\emph{Numerical solution of the multidimensional
Gelfand--Levitan equation,} J. Inverse and Ill-Posed Problems, 23 (2015),
439-450.\textsc{\ }

\bibitem{KS} \textsc{M.V. Klibanov and F. Santosa, }\emph{A computational
quasi-reversibility method for Cauchy problems for Laplace's equation}, SIAM
J. Appl. Math. 51 (1991), pp. 1653--1675.

\bibitem{KT} \textsc{M.V. Klibanov and A. Timonov, }\emph{Carleman Estimates
for Coefficient Inverse Problems and Numerical Applications, }VSP, Utrecht,
2004.

\bibitem{Klibanov:jiipp2013} \textsc{M.~V. Klibanov}, \emph{Carleman
estimates for global uniqueness, stability and numerical methods for
coefficient inverse problems}, J. Inverse and Ill-Posed Problems, 21 (2013),
pp.~477--560.

\bibitem{Ksurvey} \textsc{M.V. Klibanov}, \emph{Carleman estimates for the
regularization of ill-posed Cauchy problems}, Applied Numerical Mathematics,
94 (2015), pp. 46--74.

\bibitem{KlTh} \textsc{M. V. Klibanov and N. T. Th\`{a}nh}, Recovering
dielectric constants of explosives via a globally strictly convex cost
functional, SIAM J. Appl. Math., 75 (2015), 518-537.

\bibitem{Klibanov:jiip2017} \textsc{M.~V. Klibanov}, \emph{Convexification
of restricted {D}irichlet to {N}eumann map}, J. Inverse and Ill-Posed
Problems, 25 (2017), pp.~669--685.

\bibitem{EIT} \textsc{M. V. Klibanov, J. Li, and W. Zhang}, Electrical
impedance tomography with restricted dirichlet-to-neumann map data, 2018,
arXiv:1803.11193.

\bibitem{KlibKol} \textsc{M. V. Klibanov, A.E. Kolesov, A. Sullivan and L.
Nguyen, }\emph{A new version of the convexification method for a 1-D
coefficient inverse problem with experimental data}, Inverse Problems,
accepted for publication, 2018, a preprint is available at
https://doi.org/10.1088/1361-6420/aadbc6.

\bibitem{LRS} \textsc{M.M. Lavrentiev, V.G. Romanov and S.P. Shishatskii}, 
\emph{Ill-Posed Problems of Mathematical Physics and Analysis}, American
Mathematical Society, Providence, RI, 1986.

\bibitem{Louis:nm1985} \textsc{A.~K. Louis}, \emph{Incomplete data problems
in x-ray computerized tomography {I.} {S}ingular value decomposition of the
limited angle transform}, Numer. Math., 48 (1986), pp.~251--262.

\bibitem{LattesLions:e1969} \textsc{R.~Latt\`es and J.~L. Lions}, \emph{The
Method of Quasireversibility: Applications to Partial Differential Equations}%
, Elsevier, New York, 1969.

\bibitem{Natterer:cmsiam2001} \textsc{N.~Natterer}, \emph{The mathematics of
computerized tomography}, Classics in Mathematics. Society for Industrial
and Applied Mathematics, New York, 2001.

\bibitem{Linh:SIAM2015} \textsc{L. V. Nguyen}, \emph{On artifacts in limited
data spherical {R}adon transform: flat observation surfaces}, SIAM Journal
on Mathematical Analysis, 47 (2015), pp.~2984--3004.



\bibitem{Linh:jfaa2017} \textsc{L. V. Nguyen}, \emph{On the strength of streak artifacts in filtered back-projection
reconstructions for limited angle weighted x-ray transform}, J. Fourier
Anal. Appl., 23 (2017), pp.~712--728.

\bibitem{Loc:Arxiv2018} \textsc{L. H. Nguyen}, \emph{An inverse source problem for hyperbolic equations and the Lipschitz-like convergence of the quasi-reversibility method}, preprint, Arxiv : 1806.03921.

\bibitem{Nov} \textsc{R.G. Novikov}, \emph{An inversion formula for the
attenuated X-ray transformation}, Ark. Mat., 40 (2002), pp. 145-167.

\bibitem{PanSidkyVannier:ip2009} \textsc{X.~Pan, E.~Y. Sidky, and M. Vannier}%
, \emph{Why do commercial {CT} scanners still employ traditional, filtered
back- projection for image reconstruction?}, Inverse Problems, 25 (2009),
p.~123009.

\bibitem{Radon1917} \textsc{J.~Radon}, \emph{\"{U}ber die {B}estimmung von {F%
}unktionen durch ihre {I}ntegralwerte l\"{a}ngs gewisser {M}annigfaltigkeiten%
}, Berichte S\"{a}chsische Akademie der Wissenschaften, Leipzig,
Mathematisch-Physikalische Klasse, 69 (1917), pp.~262--277.

\bibitem{Radon:IEEE1986} \textsc{J.~Radon}, \emph{On the determination of functions from their integral values
along certain manifolds}, IEEE Transactions on Medical Imaging. \textrm{%
Translated by P.C. Parks from the original German text}, 5 (1986),
pp.~170--176.

\bibitem{T} A\textsc{.N. Tikhonov, A.V. Goncharsky, V.V. Stepanov and A.G.
Yagola}, \emph{Numerical Methods for the Solution of Ill-Posed Problems},
Kluwer, London, 1995.
\end{thebibliography}
\end{document}